\documentclass[notitlepage,leqno,10pt]{article}

\usepackage{latexsym}
\usepackage{amsmath}
\usepackage{amsfonts}
\usepackage{amssymb}
\usepackage{hyperref}
\renewcommand{\a }{\alpha }
\renewcommand{\b }{\beta }
\renewcommand{\d}{\delta }
\newcommand{\D }{\Delta }

\newcommand{\e }{\varepsilon }
\newcommand{\g }{\gamma}

\newcommand{\G }{\Gamma }
\renewcommand{\l }{\lambda }
\renewcommand{\L }{\Lambda }

\newcommand{\n }{\nabla }
\newcommand{\var }{\varphi }

\newcommand{\s }{\sigma }

\renewcommand{\t }{\tau }
\renewcommand{\th }{\theta }
\renewcommand{\o }{\omega }

\newcommand{\z }{\zeta }
\newcommand{\ov}{\overline}
\newcommand{\pa}{\partial}
\newcommand{\tPhi}{\tilde{\Phi}}

\newcommand{\be}{\begin{equation}}
\newcommand{\ee}{\end{equation}}
\newenvironment{pf}{\noindent{\sc Proof}.\enspace}{\rule{2mm}{2mm}\medskip}
\newenvironment{pfn}{\noindent{\sc Proof}}{\rule{2mm}{2mm}\medskip}

\newcommand{\R}{\mathbb{R}}

\newcommand{\N}{\mathbb{N}}
\newcommand{\C}{\mathbb{C}}
\newcommand{\no}{\noindent}
\newcommand{\oy }{\overline{y}}

\newcommand{\os}{\ov{s}}

\author{Fethi MAHMOUDI$^{\rm a}$ \and Andrea MALCHIODI$^{\rm a}$}

\date{}

\title{Solutions to the nonlinear Schr\"odinger equation \\
carrying momentum along a curve. \\ Part II: proof of the existence result }

\begin{document}

\newtheorem{lem}{Lemma}[section]
\newtheorem{pro}[lem]{Proposition}
\newtheorem{thm}[lem]{Theorem}
\newtheorem{rem}[lem]{Remark}
\newtheorem{cor}[lem]{Corollary}
\newtheorem{df}[lem]{Definition}

\maketitle

\begin{center}

$^{\rm a}${\small SISSA, Sector of Mathematical Analysis \\ Via
Beirut 2-4, 34014 Trieste, Italy}

\end{center}

\footnotetext[1]{E-mail addresses: mahmoudi@ssissa.it
(F.Mahmoudi), malchiod@sissa.it (A. Malchiodi)}

\noindent {\sc abstract}. We prove existence of a special class of
solutions to the (elliptic) Nonlinear Schr\"odinger Equation $- \e^2
\D \psi + V(x) \psi = |\psi|^{p-1} \psi$ on a manifold or in the
Euclidean space. Here $V$ represents the potential, $p$ is an
exponent greater than $1$ and $\e$ a small parameter corresponding
to the Planck constant. As $\e$ tends to zero (namely in the
semiclassical limit) we prove existence of complex-valued solutions
which concentrate along closed curves, and whose phase in highly
oscillatory. Physically, these solutions carry quantum-mechanical
momentum along the limit curves. In the first part of this work we
identified the limit set and constructed approximate solutions,
while here we give the complete proof of our main existence result
Theorem \ref{t:main}.
\begin{center}

\bigskip\bigskip


\bigskip

\centerline{\bf AMS subject classification: 34B18, 35B25, 35B34,
35J20, 35J60}

\end{center}

\section{Introduction}\label{s:i}

In this paper we continue our study of \cite{mmm1}, concerning
concentration phenomena for solutions of the singularly-perturbed
elliptic problem
\begin{equation}\label{eq:pe}
    - \e^2 \D_g \psi + V(x) \psi = |\psi|^{p-1} \psi \qquad
    \hbox{ on } M,
\end{equation}
where $M$ is an $n$-dimensional compact manifold (or the flat
Euclidean space $\R^n$),  $V$ a smooth positive function on $M$
satisfying the properties
\begin{equation}\label{eq:VV}
    0 < V_1 \leq V \leq V_2; \qquad \qquad \|V\|_{C^3} \leq V_3,
\end{equation}
$\psi$ a complex-valued function, $\e > 0$ a small parameter and $p$
is an exponent greater than $1$. Here $\D_g$ stands for the
Laplace-Beltrami operator on $(M,g)$.

Solutions to \eqref{eq:pe} represent {\em standing waves} of the
Nonlinear Schr\"odinger Equation, and here we are interested in
the {\em semiclassical limit}, namely the asymptotics of solutions
when the parameter $\e$ (representing the Planck constant) tends
to zero. Typically, if concentration occurs near some point $x_0
\in M$, such solutions behave like $\psi_\e(x) \simeq u\left(
\frac{dist(x,x_0)}{\e} \right)$, where $dist(\cdot, \cdot)$
denotes the distance on $M$ and where $u$ solves the equation
\begin{equation}\label{eq:vx0}
    - \D u + V(x_0) u = u^p \qquad \quad \hbox{ in } \R^n.
\end{equation}
We refer the reader to the introduction of \cite{mmm1} for the
motivation of the study and for a brief description of the
existing results in the literature about this topic. There are
several works concerning standing waves concentrating at a single
(or multiple) points of $M$, and for which the corresponding
solutions of \eqref{eq:vx0} decay to zero at infinity. On the
other hand, only very recently it has been proven existence of
solutions concentrating at higher dimensional sets, like curves or
manifolds. In all these results (except for \cite{dap}), the
profile is given by (real) solutions to \eqref{eq:vx0} which are
independent of some of the variables and hence do not tend to zero
at infinity: if concentration occurs near a $k$-dimensional set,
then the profile in the directions orthogonal to the limit set
will be given by a soliton in $\R^{n-k}$.

\

\noindent In this paper we are going to construct a different type
of solutions. These still concentrate along curves in $M$, but their
phase is highly oscillatory along the limit set. More precisely, we
consider standing waves (namely the solution of \eqref{eq:vx0})
whose profile  has the following expression
\begin{equation}\label{eq:proflimit}
    \phi(x',x_n) = e^{-i \hat{f} x_n} \hat{U}(x'),
    \qquad \quad x' = (x_1, \dots, x_{n-1}),
\end{equation}
where $\hat{f}$ is some constant and $\hat{U}(x')$ a real
function. With this choice of $\phi$, if concentration occurs near
some point $x_0$, then the function $\hat{U}$ satisfies the
equation
\begin{equation}\label{eq:ovU}
    - \D \hat{U} + \left( \hat{f}^2 + V(x_0) \right) \hat{U}
    = |\hat{U}|^{p-1} \hat{U} \qquad \quad \hbox{ in } \R^{n-1},
\end{equation}
and decays to zero at infinity. Solutions to \eqref{eq:ovU} can be
found by considering the radial function $U : \R^{n-1} \to \R$ which
solves
\begin{equation}\label{eq:UU}
    - \D U + U = U^p \qquad \quad \hbox{ in } \R^{n-1}.
\end{equation}
It is known that $U$ (and its derivatives) behaves at infinity like
\begin{equation}\label{eq:decoU}
    U(r) \simeq e^{-r} r^{- \frac{n-2}{2}} \qquad \qquad \hbox{ as } r \to + \infty.
\end{equation}
Using the scaling
\begin{equation}\label{eq:ovhovk}
    \hat{U}(x') = \hat{h} U(\hat{k} x'), \qquad \qquad \hat{h} = \left(
  \hat{f}^2 + V(x_0) \right)^{\frac{1}{p-1}}, \quad \hat{k} = \left(
  \hat{f}^2 + V(x_0) \right)^{\frac 12},
\end{equation}
in \eqref{eq:proflimit} the constant $\hat{f}$ can be chosen
arbitrarily, and then $\hat{h}, \hat{k}$ are determined according to
the last formula, depending on $V(x_0)$. Indeed $\hat{f}$ represents
the speed of the phase oscillation, and is physically related to the
velocity of the quantum-mechanical particle associated to the wave
function. If concentration occurs near some closed curve $\g =
\g(\ov{s})$ in $M$, and if we allow the parameter $\hat{f}$ to
depend on the variable $s$, then the solution $\psi$ will be of the
form
\begin{equation}\label{eq:prof}
    \psi(\ov{s},\z) \simeq e^{- i \frac{f(\ov{s})}{\e}}
    h(\ov{s}) U\left(\frac{k(\ov{s}) y}{\e}\right),
\end{equation}
where $\ov{s}$ stands for the arc-length parameter of $\g$, and $y$
for a system geodesic coordinates normal to $\g$. Here the functions
$h(\ov{s})$ and $k(\ov{s})$ are chosen so that
\begin{equation}\label{eq:ovhovks}
    h(\ov{s}) = \left(
  (f'(\ov{s}))^2 + V(\ov{s}) \right)^{\frac{1}{p-1}}, \qquad \qquad
  k(\ov{s}) = \left( (f'(\ov{s}))^2 + V(\ov{s}) \right)^{\frac 12}.
\end{equation}
basically replacing $\hat{f}$ with $f'(\ov{s})$ in
\eqref{eq:ovhovk}.

If $\g$ is given, it was shown in \cite{mmm1} (using formal
expansions in $\e$) that the corresponding function $f$ should
satisfy the following condition
\begin{equation}\label{eq:f'Cintr}
    f'(\ov{s}) \simeq \mathcal{A} h^{\s}(\ov{s}) \qquad
    \quad \hbox{ with } \quad \s = \frac{(n-1)(p-1)}{2}-2,
\end{equation}
where $\mathcal{A}$ is an arbitrary constant. At this point, only
the limit curve $\g$ should be determined. Since we require the
function $f$ to be periodic, if we consider variations of $\g$
(for all of which \eqref{eq:f'Cintr} holds true) then it is
natural to work in the restricted class
$$
  \G := \left\{ \g : \R \to M \hbox{ periodic} \; : \;
  \mathcal{A} \int_\g h(\ov{s})^\s d \ov{s} = \int_\g f'(\ov{s})^\s
  d \ov{s} = \hbox{constant} \right\},
$$
where, as before, $\ov{s}$ stands for the arc-length parameter. It
is shown in \cite{mmm1} that the candidate limit curves are critical
points of the functional $\g \mapsto \int_\g h^{\th}(\ov{s}) d
\ov{s}$, where
\begin{equation*}
    \th = p + 1 - \frac 12 (p-1)(n-1).
\end{equation*}
With a direct computation one can prove that the extremality
condition is the following
\begin{equation}\label{eq:eulerintr}
  \n^N V = \left( \frac{p-1}{\th}
  h^{p-1} - 2 \mathcal{A}^2 h^{2\s} \right){\bf H},
\end{equation}
where $\n^N V$ represents the normal gradient of $V$ and ${\bf H}$
is the  curvature vector of $\g$.

Similarly, via some long but straightforward calculation, one can
find a natural non-degeneracy condition for stationary points, which
is expressed by the invertibility of the operator $\mathfrak{J}$,
acting on normal sections $\mathcal{V}$ to $\g$, which in components
is given by
\begin{eqnarray}\label{eq:2ndvarfin4} \nonumber
  (\mathfrak{J} \mathcal{V})^m & = & - \left( h^\th - \frac{2 \mathcal{A}^2
  \th}{p-1} h^\s \right) \ddot{\mathcal{V}}^m - \th
  \left( h^{\th-1} - \frac{2 \mathcal{A}^2 \s}{p-1} h^{\s-1}
  \right) h' \dot{\mathcal{V}}^m
  + \frac{\th}{p-1} h^{-\s} ((\n^N)^2 V)[\mathcal{V},E_m]
  \\  & + & \frac 12 \left( h^\th - \frac{2
  \mathcal{A}^2 \th}{p-1} h^\s \right) \left(
  \sum_j \left(  \partial^2_{jm} g_{11} \right) {\mathcal{V}}^j \right)
  - 2 \mathcal{A} \mathcal{A}'_1 \frac{(\th-\s)
  h^{p-1}}{\left[ (p-1) h^{\th} - 2 \s \mathcal{A}^2 h^{\s}
  \right]} H^m \\
  & + & H^m \langle {\bf H}, \mathcal{V} \rangle
  \left[ \frac{- (p-1) \left( 3 + \frac{\s}{\th} \right)
  h^{2\th} - \frac{16 \s \th \mathcal{A}^4}{p-1} h^{2\s} + 2 \mathcal{A}^2 (5\s +
  3\th) h^{\th+\s}}{(p-1) h^{\th} - 2 \mathcal{A}^2 \s h^{\s}}
  \right], \nonumber
\end{eqnarray}
$m = 2, \dots, n$. We refer to Section 2 in \cite{mmm1} for the
notation used in this formula. We point out that, since
\eqref{eq:f'Cintr} determines only the derivative of the phase, to
obtain periodicity we need to introduce a nonlocal term, denoted
here with $\mathcal{A}'_1$. Letting $L(\g)$ be the length of the
curve $\g$, our main result is the following.

\begin{thm}\label{t:main}
Let $M$ be a compact $n$-dimensional manifold and let $V : M \to \R$
be a smooth positive function (or let $M = \R^n$ and let $V$ satisfy
\eqref{eq:VV}) and $1 < p < \frac{n+1}{n-3}$. Let $\g$ be a simple
closed curve in $M$: then there exists a positive constant
$\mathcal{A}_0$, depending on $V|_\g$ and $p$ for which the
following holds. If $0 \leq \mathcal{A} < \mathcal{A}_0$, if $\g$
satisfies \eqref{eq:eulerintr} and the operator in
\eqref{eq:2ndvarfin4} is invertible on normal sections of $\g$,
there is a sequence $\e_k \to 0$ such that problem $(NLS_{\e_k})$
possesses solutions $\psi_{\e_k}$ having the asymptotics in
\eqref{eq:prof}, with $f$ satisfying \eqref{eq:f'Cintr}.
\end{thm}

\noindent As a consequence of this theorem, see Corollary 1.3 in
\cite{mmm1}, we prove a conjecture posed in \cite{amn1} for the case
of one-dimensional limit sets. We also improve the result in
\cite{dap}, in the sense that we characterize explicitly the limit
set and we do not require any symmetry on the potential $V$: indeed
in \cite{dap} $V$ is assumed cylindrically symmetric in $\R^3$, and
solutions are found via separation of variables. The restriction on
the exponent $p$ is natural since it is a necessary condition for
the solutions of \eqref{eq:UU} to vanish at infinity by the {\em
Pohozaev's identity}. The smallness condition on the constant
$\mathcal{A}$ and the fact that concentration is not proved for all
the values of (small) $\e$ are discussed below in the introduction.
For the latter issue and for the main difficulties caused by
removing the symmetries see also the introduction of \cite{mmm1}.

\

\no The main goal of Part I, \cite{mmm1}, was to show that the
condition \eqref{eq:eulerintr} and the non-degeneracy of the
operator in \eqref{eq:2ndvarfin4}, arising from the {\em reduced
functional} $\g \mapsto \int_\g h^\th(\ov{s}) d \ov{s}$, appear
naturally when considering \eqref{eq:pe}, and in particular when we
try to solve it formally with an expansion in power series of $\e$.
To explain this fact, it is convenient to scale problem
\eqref{eq:pe} in the following way
\begin{equation}\label{eq:new}
  - \D_{g_\e} \psi + V(\e x) \psi = |\psi|^{p-1} \psi \qquad \hbox{ in
  } M_\e,
\end{equation}
where $M_\e$ denotes the manifold $ M$ endowed with the scaled
metric $g_\e = \frac{1}{\e^2} g$ (with an abuse of notation we might
often write $M_\e=\frac1\e M$). We are now looking for a solution
concentrated near the dilated curve $\g_\e := \frac 1 \e \g$. We let
$s$ be the arc-length parameter of $\g_\e$, so that $\ov{s} = \e s$,
and we let $(E_j)_{j=2, \dots, n}$ denote an orthonormal frame in $N
\g$ (the normal bundle of $\g$) transported parallely with respect
to  the normal connection, see Section 2 in \cite{mmm1}: we also let
$(y_j)_j$ be a corresponding set of normal coordinates. Since we
want to allow some flexibility both in the choice of the phase and
of the curve of concentration, we define $\tilde{f}_0(\ov{s}) =
f(\ov{s}) + \e f_1(\ov{s})$, and we set $z_j = y_j -
\Phi_j(\ov{s})$, where $(\Phi_j)_{j=2,\dots,n}$ are the components
(with respect to the above coordinates $y$) of a section $\Phi$ in
$N \g$. Then, with a formal expansion of $\psi$ in powers of $\e$ up
to the second order, in the coordinates $(s,z)$ near $\g_\e$, we set
  \begin{equation*}
  \psi_{2,\e}(s,z) = e^{- i \frac{\widetilde{f}_0(\e s)}{\e}}
  \left\{ h(\e s) U \left( k( \e
  s) z \right) + \e \left[ w_r + i w_i \right] +
  \e^2 \left[ v_r + i v_i \right]\right\}, \qquad \quad s
  \in [0, L/\e], z \in \R^{n-1},
\end{equation*}
($L = L(\g)$) for some corrections  $w_r, w_i, v_r, v_i$ (which have
to be determined) to the above approximate solutions. We saw in
Section 3 of \cite{mmm1} that these terms solve equations of the
form $\mathcal{L}_r w_r = \mathcal{F}_r$, $\mathcal{L}_i w_i =
\mathcal{F}_i$, $\mathcal{L}_r v_r = \tilde{\mathcal{F}}_r$,
$\mathcal{L}_i v_i = \tilde{\mathcal{F}}_i$ where
\begin{equation}\label{eq:LrLi}
    \left\{
      \begin{array}{ll}
        \mathcal{L}_r v = - \D_z v + V(\ov{s}) v - p h(\ov{s})^{p-1}
    U(k(\ov{s}) z)^{p-1} v; & \\
        \mathcal{L}_i v = - \D_z
    v + V(\ov{s}) v - h(\ov{s})^{p-1} U(k(\ov{s}) z)^{p-1} v &
      \end{array}
    \right.
      \qquad
    \quad \hbox{ in } \R^{n-1},
\end{equation}
and where $\mathcal{F}_r$, $\mathcal{F}_i$, $\tilde{\mathcal{F}}_r$,
$\tilde{\mathcal{F}}_i$ are given data which depend on $V$, $\g$,
$\ov{s}$, $\mathcal{A}$, $\Phi$ and $f_1$. The operators
$\mathcal{L}_r$ and $\mathcal{L}_i$ are Fredholm (and symmetric)
from $H^2(\R^{n-1})$ into $L^2(\R^{n-1})$, and the above equations
for the corrections can be solved provided the right-hand sides are
orthogonal to the kernels of these operators. As explained in
\cite{mmm1}, the condition \eqref{eq:eulerintr} and the
non-degeneracy of the operator $\mathfrak{J}$ allow us to determine
$w_r, w_i$ and $v_r, v_i$ respectively, namely to solve
\eqref{eq:new} at order $\e$ first, and then at order $\e^2$.

\

\no To make the above arguments rigorous, we can start with
 an approximate solution ${\psi}_{0,\e}$ behaving like ${\psi}_{0,\e} \simeq e^{-
i \frac{f(\e s)}{\e}} h(\e s) U(k(\e s) z)$, and try to find a true
solution of the form $e^{- i \frac{{\tilde f}(\e s)}{\e}} [h(\e s)
U(k(\e s) z) + \tilde{w}]$, with $\tilde{w}$ suitably small and
$\tilde{f}$ close to $f$, via some local inversion arguments. From a
linearization of the equation near ${\psi}_{0,\e}$, the operator
$L_\e$ acting on $\tilde{w}$ in the coordinates $(s,z)$ is then the
following
\begin{equation}\label{eq:lintildew}
    L_\e \tilde{w} : = - \pa^2_{ss} \tilde{w} - \D_z \tilde{w} +
    V(\e x) - |{\psi}_{0,\e}|^{p-1} \tilde{w} - (p-1) |{\psi}_{0,\e}|^{p-3}
    {\psi}_{0,\e} \Re({\psi}_{0,\e} \ov{\tilde{w}}).
\end{equation}
Here $\Re$ denotes the real part. Decomposing first $\tilde{w}$ into
its real and imaginary parts, and then in Fourier modes with respect
to the variable $\e s$, we can write $\tilde{w} = \tilde{w}_r + i
\tilde{w}_i = \sum_j \sin(j \e s) \tilde{w}_{r,j}(z) + i \sum_j
\sin(j \e s) \tilde{w}_{i,j}(z)$ (forgetting for simplicity about
the cosine functions). If we take (as a model problem) $V \equiv 1$,
then the operators (in the $z$ variables) acting on the real and
imaginary components are respectively $\mathcal{L}_r + \e^2 j^2$ and
$\mathcal{L}_i + \e^2 j^2$. It is well-known, see for example
\cite{kwo}, that $\mathcal{L}_r$ has a single negative eigenvalue, a
kernel with multiplicity $n-1$ spanned by the functions $\pa_l
U(k(\ov{s}) z)$, $l = 2, \dots, n$ (the generators of the normal
translations), while all the remaining eigenvalues are positive. The
operator $\mathcal{L}_i$ instead has one zero eigenvalue with
eigenfunction $U(k(\ov{s}) z)$ (the generator of complex rotations)
and all the remaining eigenvalues positive.

As a consequence, the kernels of $\mathcal{L}_r$ and $\mathcal{L}_i$
produce a sequence of eigenvalues for $L_\e$ which behave {\em
qualitatively} like $\e^2 j^2$, and for small values of $j$ these
become resonant. With an accurate expansion of these eigenvalues,
one finds that the non-degeneracy assumption on
\eqref{eq:2ndvarfin4} prevents each of them to vanish: anyway, a
direct application of the implicit function theorem is not possible
since a further resonance phenomenon occurs. This arises from the
fact that $\mathcal{L}_r$  possesses a negative eigenvalue as well,
which generates an extra sequence of eigenvalues of $L_\e$,
qualitatively of the form $-1 + \e^2 j^2$, $j \in \N$. This
resonance is typical of concentration for \eqref{eq:pe} along sets
of positive dimension, and the only hope to get invertibility is to
choose the values of $\e$ appropriately. Indeed, differently from
the previous sequence of eigenvalues, this new one causes a
divergence of the Morse index when $\e$ tends to zero, and the
presence of a kernel for some epsilon's is unavoidable. The
eigenfunctions corresponding to these eigenvalues will have faster
and faster oscillations along the limit curve $\g$.

This phenomenon is also present when one looks for solutions of the
singularly perturbed problem $- \e^2 \D u + u = u^p$ in bounded
domains of $\R^n$, when Neumann boundary conditions are imposed. In
the papers \cite{mm}, \cite{mal}, \cite{malm}, \cite{malm2}
concentration along sets of dimension $k = 1, \dots, n-1$ has been
proved, and analogous spectral properties hold true. By the Weyl's
asymptotic formula, if solutions concentrate along a set of
dimension $d$, the counterpart of the latter sequence of eigenvalues
behaves like $- 1 + \e^2 j^{\frac 2d}$, and the average distance
between those close to zero is of order $\e^d$. The resonance
phenomenon was taken care of using a theorem by T. Kato, see
\cite{ka}, page 445, which allows to differentiate eigenvalues with
respect to $\e$. In the aforementioned papers it was shown that when
varying the parameter $\e$ the spectral gaps near zero almost do not
shrink, and invertibility can be obtained for a large family of
epsilon's.

However when the concentration set is one-dimensional the spectral
gaps of the resonant eigenvalues (with fast-oscillating
eigenfunctions) are relatively large, of order $\e$, and the profile
of the corresponding eigenfunctions can be analyzed by means of a
scalar function on $[0,L]$ (see below in the introduction and in
Subsection \ref{sss:pr3}). This fact indeed allows sometimes to
bypass Kato's theorem and to use a more direct approach, employed in
\cite{mp} to study existence of constant mean curvature surfaces of
cylindrical type embedded in manifolds, and in \cite{dkw} for
studying solutions of \eqref{eq:pe} in $\R^2$. We can partially take
advantage of these techniques, see the comments in Section
\ref{s:LS}, but some new difficulties arise due to the fast phase
oscillations in \eqref{eq:prof}. We describe them below, together
with the strategy of the proof.

\

\noindent By the above discussion, we expect to find three possible
resonances: two of them for small values of the index $j$ (with
eigenvectors roughly of the form $e^{-i\frac{f(\e s)}{\e}}\pa_l
U(k(\ov{s}) z)\sin(\e js)$, $l = 2, \dots, n$, and $ie^{-i\frac{f(\e
s)}{\e}} U(k(\ov{s}) z)\sin(\e js)$ respectively) and a third one
for $j$ of order $\frac 1 \e$, precisely when $- \e^2 j^2$ coincides
with the first eigenvalue of $\mathcal{L}_r$.

To understand this behavior, we first study the spectrum of a {\em
model} operator similar to \eqref{eq:lintildew},  where we assume $V
\equiv \hat{V} > 0$ and $\psi_{0,\e}$ to coincide with the function
in \eqref{eq:proflimit}. For this case we characterize completely
the spectrum of the operator and the properties of the
eigenfunctions, see Subsection \ref{ss:appker} and in particular
Proposition \ref{p:alphamu}. The condition on the smallness of
$\mathcal{A}$ appears precisely here (and only here), and is used to
show that the resonant eigenvalues are only of the forms described
above. Removing the smallness assumption might indeed lead to
further resonance phenomena, see Remark \ref{r:Asmall2} for further
comments.

We next consider the case of non-constant potential $V$: since this
has a slow dependence in $s$ along $\g_\e$, one might guess that the
approximate kernel of $L_\e$ (see \eqref{eq:lintildew}) might be
obtained from that for constant $V$, introducing also a slow
dependence in $s$ of the profile of these functions. With this
criterion, given a small positive parameter $\d$, we introduce a set
$K_\d$ (see \eqref{eq:Kd} and the previous formulas) consisting of
{\em candidate} approximate eigenfunctions on $L_\e$, once
multiplied by the phase factor $e^{- i \frac{f(\e s)}{\e}}$. More
comments on the specific construction of this set can be found in
Subsection \ref{ss:appker}, especially before \eqref{eq:Kd}.

In Proposition \ref{p:inv} we show that this guess is indeed
correct: in fact, we prove that the operator $L_\e$ is invertible
provided we restrict ourselves to the subset $\ov{H}_\e$ of
functions which are orthogonal to $e^{- i \frac{f(\e s)}{\e}} K_\d$.
This property allows us to solve the equation up to a lagrange
multiplier in $K_\d$, see Proposition \ref{p:truefndec}. For
technical reasons, we prove invertibility of $L_\e$ in suitable
weighted norms, which are convenient to deal with functions decaying
exponentially away from $\g_\e$. As done in \cite{dkw}, \cite{mm}
and \cite{mal}, this decay allows us to shift the problem from the
whole manifold $M_\e$ to the normal bundle $N \g_\e$ via a
localization method, see Subsection \ref{subsec:loc}.

Compared to the other results in the literature which deal with this
kind of resonance phenomena, the approximate kernel here depends
{\em genuinely} on the variable $s$ (in \cite{mm}, \cite{mal},
\cite{malm}, \cite{malm2},  \cite{mp} the problem is basically
homogeneous along the limit set, while in \cite{dkw} it can be made
such through a change of variables). To deal with this feature,
which mostly causes difficulties in Proposition \ref{p:inv}, we
localize the problem in the variable $s$ as well. Multiplying by a
cutoff function depending on $s$, we show that orthogonality to
$K_\d$ implies approximate orthogonality to the set $\hat{K}_\d$,
see \eqref{eq:hatKd} and the previous formulas, which is the
counterpart of $K_\d$ for a potential freezed at some point in
$\g_\e$: once this is shown, we use the spectral analysis of
Proposition \ref{p:alphamu}.

\

\noindent Section \ref{s:as} is devoted to choose a family of
approximate solutions to \eqref{eq:new}: since we have many small
eigenvalues appearing, it is natural trying to look for functions
which solve \eqref{eq:new} as accurately as possible. Our final goal
is to annihilate the Lagrange multiplier in Proposition
\ref{p:truefndec}, and to do this we choose approximate solutions
$\tilde{\Psi}_{2,\e}$ (in the notation of Section \ref{s:as}) which
depend on suitable parameters: a normal section $\Phi$, a phase
factor $f_2$ and a real function $\b$. The latter parameters
correspond to different components of $K_\d$, and are related to the
kernels of $\mathcal{L}_r (+ \e^2 j^2)$ and $\mathcal{L}_i (+ \e^2
j^2)$, see the above comments. The function $\b$ in particular is
highly oscillatory, and takes care of the resonances due to the fast
Fourier modes.

Differently from \cite{mmm1} (see in particular  Section 4 there),
where the expansions were only formal, we need here to derive
rigorous estimates on the error terms, and to study in particular
their Lipschitz dependence on the data $\Phi, f_2$ and $\b$.
Proposition \ref{p:errest} collects the final expression for $-
\D_{g_\e} \psi + V(\e x) \psi - |\psi|^{p-1} \psi$ on the
approximate solutions $\tilde{\Psi}_{2,\e}$: the error terms
$\tilde{\mathfrak{A}}$'s are listed (and estimated) before in that
section, together with their Lipschitz dependence on the parameters.

\

\no Finally, after performing a Lyapunov-Schmidt reduction onto the
set $K_\d$, see Proposition \ref{p:exhatphi}, we study the
bifurcation equation  in order to annihilate the Lagrange
multiplier. In doing this we use crucially the computations in Part
I, \cite{mmm1}, together with the error estimates in Section
\ref{s:as}. In particular, for $\Phi$ and $f_2$ we find as main
terms respectively the operator $\mathfrak{J}$ in
\eqref{eq:2ndvarfin4} and the one in the left-hand side of
\eqref{eq:f1111}, both appearing when we performed formal
expansions: these operators are both invertible by our assumptions,
and therefore we are able to determine $\Phi$ and $f_2$ without
difficulties.

The operator acting on $\b$ instead is more delicate, since it is
\underline{qualitatively} of the form
\begin{equation}\label{eq:eqaintr}
    - \e^2 \b''(\ov{s}) + \l(\ov{s}) \b \qquad \qquad \hbox{ on }
    [0,L],
\end{equation}
with periodic boundary conditions, where $\l$ is a negative
function. This operator is precisely the one related to the peculiar
resonances described above. In particular it is resonant on
frequencies of order $\frac 1\e$, and this requires to choose a norm
for $\b$ which is weighted in the Fourier modes, see
\eqref{eq:normsharp} and Subsection \ref{sss:pr3}. For operators
like that in \eqref{eq:eqaintr} there is in general a sequence of
epsilon's for which a non-trivial kernel exists. Using Kato's
theorem though, as in \cite{mal}, \cite{malm}, \cite{malm2},
\cite{mm} and \cite{mmp}, we provide estimates on the derivatives of
the eigenvalues with respect to $\e$, showing that for several
values of this parameter the operator acting on $\b$ is invertible.
In this operation also the value of the constant $\mathcal{A}$, see
\eqref{eq:f'Cintr}, has to be suitably modified (depending on $\e$),
in order to preserve the periodicity of our functions. Once we have
this, we apply the contraction mapping theorem to solve the
bifurcation equation as well.

\

\noindent The results in this paper and in \cite{mmm1} are briefly
summarized in the note \cite{mmmnote}.

\

\begin{center}
{\bf Notation and conventions}
\end{center}

\no  Dealing with coordinates, capital letters like $A, B, \dots$
will vary between $1$ and $n$ while indices like $j, l, \dots$
will run between $2$ and $n$. The symbol $i$ will stands for the
imaginary unit.

\

\no For summations, we  use the standard convention of summing
terms where repeated indices appear.

\

\no We will often work with coordinates $(\ov{s}, \ov{y}_1, \dots,
\ov{y}_{n-1})$ near the curve $\g$, where $\ov{s}$ is the arc-length
parameter of $\g$ and the $\ov{y}_j$'s are Fermi normal coordinates
(see the next section). The dilated curve $\g_\e:=\frac1\e\g$ will
be parameterized by $s=\frac1\e\ov{s}$, and we also denote by $y$
the scaled normal  coordinates. The length of $\g$ is denoted by
$L$. When dealing with functions (or normal sections) depending on
the variables $s$ or $\ov{s}$, the notation {\em prime} will always
denote the derivative with respect to $\ov{s}$. When we
differentiate with respect to $s$, we usually adopt the symbol
$\pa_s$.

\

\no For simplicity, a constant $C$ is allowed to vary from one
formula to another, also within the same line.

\

\no For a real positive variable $r$ and an integer $m$, $O(r^m)$
(resp. $o(r^m)$) will denote a complex-valued  function for which
$\left| \frac{O(r^m)}{r^m} \right|$ remains bounded (resp. $\left|
\frac{o(r^m)}{r^m} \right|$ tends to zero) when $r$ tends to zero.
We might also write $o_\e(1)$ for a quantity which tends to zero
as $\e$ tends to zero.

\

\no Sometimes we shall need to work with integer indices $j$ which
belong to sets depending on $\e$. e.g. $\{0, 1, \dots, [1/\e]\}$,
where the latter square brackets stand for the integer part. For
convenience, we will often omit the to add the square brackets,
assuming that this convention is understood.

\

\section{Lyapunov-Schmidt reduction of the problem}\label{s:LS}

In this section we show how to reduce problem \eqref{eq:new} to a
system of three ordinary (integro-)differential equations on $\R /
[0,L]$. We first introduce a metric on the normal bundle $N \g_\e$
of $\g_\e$ and  then study operators which mimic the properties of
the linearization of \eqref{eq:new} near an approximate solution.
Next, we turn to the reduction procedure: this follows basically
from a localization method, since the functions we are dealing with
have an exponential decay away from $\g_\e$. We introduce a set
$K_\d$ consisting of approximate (resonant) eigenfunctions of the
linearized operator $L_\e$: calling $\ov{H}_\e$ the orthogonal
complement of this set (which has to be multiplied by a phase factor
close to $e^{- i \frac{f(\e s)}{\e}}$) we show in Proposition
\ref{p:truefndec} that $L_\e$ is invertible on the projection onto
this set, once suitable weighted norms are introduced.

\subsection{A metric structure on $N \g_\e$}\label{ss:model}

In this subsection we define a metric $\hat{g}_\e$ on $N \g_\e$, the
normal bundle to $\g_\e$, and then introduce some basic tools which
are useful for working in local coordinates on this set.

\

\no First of all, we choose a local orthonormal frame $(E_i)_i$ in
$N \g$ and, using the notation of Subsection 2.2 in \cite{mm}, we
set $\n^N_{\pa_{\ov{s}}} E_j = \b^l_j(\pa_{\ov{s}}) E_l$, $j,l = 1,
\dots, n-1$. If we impose that the $E_j$'s are transported parallely
via the normal connection $\n^N$, as in Subsection 2.1 of
\cite{mmm1}, we find that $\b^l_j(\pa_{\ov{s}}) \equiv 0$ for all
$j,l$. As a consequence, see formula $(18)$ in \cite{mm}, we have
that if $(V^j)_j$, $j = 1, \dots, n-1$ is a normal section to $\g$,
then the components of the normal Laplacian $\D^N V$ are simply
given by
\begin{equation}\label{eq:normlapl}
    (\D^N V)^j = \D_\g (V^j) = \pa^2_{\ov{s} \ov{s}} V^j,
    \qquad \quad j = 1, \dots, n-1.
\end{equation}

We next define a metric $\hat{g}$ on $N \g$ as follows. Given $v \in
N \g$, a tangent vector $W \in T_v N \g$ can be identified with the
velocity of a curve $w(t)$ in $N \g$ which is equal to $v$ at time
$t = 0$. The metric $\hat{g}$ on $N \g$ acts on an arbitrary couple
$(W, \tilde{W}) \in (T_v N \g)^2$ in the following way (see
\cite{do2}, page 79)
$$
  \hat{g}(W,\tilde{W}) = g \left( \pi_* W, \pi_* \tilde{W} \right) + \left\langle
  \frac{D^N w}{d t}|_{t=0}, \frac{D^N \tilde{w}}{d t}|_{t=0} \right\rangle_N.
$$
In this formula $\pi$ denotes the natural projection from $N \g$
onto $\g$, $\frac{D^N w}{d t}$ the (normal) covariant derivative
of the vector field $w(t)$ along the curve $\pi \; w(t)$, and
$\tilde{w}(t)$ stands for a curve in $N \g$ with initial value $v$
and initial velocity equal to $\tilde{W}$.

Following the notation in Subsection 2.1 of \cite{mmm1} we have
that, if $w(t) = w^j(t) E_j(t)$, then
$$
  \frac{D^N w}{d t} = \frac{d w^j(t)}{dt} E_j(t).
$$
Therefore, if we choose a system of coordinates $(\os,\oy)$ on $N
\g$ defined by
$$
  (\os, \oy) \in \R \times \R^{n-1} \qquad \mapsto \qquad \oy_j E_j(\g(\os)),
$$
we get that
$$
  \hat{g}_{11}(\os,\oy) = g_{11}(\os) + \oy_l
  \oy_j \left\langle \n^N_{E_1} E_l, \n^N_{E_1} E_j
  \right\rangle_N = g_{11}(\os) \equiv 1,
$$
and
$$
  \hat{g}_{1 \ov{l}}(\os,\oy) \equiv 0; \qquad \qquad
  \hat{g}_{\ov{l} \ov{j}}(\os,\oy) = \d_{\ov{l} \ov{j}},
$$
where we have set $\pa_{\ov{l}} = \frac{\pa}{\pa \oy_l}$. We
notice also that the following co-area type formula holds, for any
smooth compactly supported function $f : N \g  \to \R$
\begin{equation}\label{eq:coarea}
    \int_{N \g} f dV_{\hat{g}} = \int_\g \left( \int_{N \g(\os)}
    f(\oy) d \oy \right) d \os.
\end{equation}
This follows immediately from the fact that $\det \hat{g} = \det
g$, and by our choice of $(\ov{s},\ov{y})$.

Since in the above coordinates the metric $\hat{g}$ is diagonal,
the Laplacian of any (real or complex-valued) function $\phi$
defined on $N \g$ with respect to this metric is
\begin{equation*}
    \D_{\hat{g}} \phi = \pa^2_{\os \os} \phi + \pa^2_{\ov{j} \ov{j}}
    \phi \qquad \quad \hbox{ in } N \g.
\end{equation*}
We endow next $N \g_\e$ with a natural metric, inherited by
$\hat{g}$ through a scaling. If $T_\e$ denotes the dilation $x
\mapsto \e x$, we define a metric $\hat{g}_\e$ on $N \g_\e$ simply
by
$$
  \hat{g}_\e = \frac{1}{\e^2} [(T_\e)_* \hat{g}].
$$
In particular, choosing coordinates $(s,y)$ on $N \g_\e$ via the
scaling $(\os,\oy) = \e (s,y)$, one easily checks that the
components of $\hat{g}_\e$ are given by
$$
  (\hat{g}_\e)_{11}(s,y) = g_{11}(\os) \equiv 1; \qquad \qquad
  (\hat{g}_\e)_{1l}(s,y) \equiv 0; \qquad \qquad
  (\hat{g}_\e)_{lj}(s,y) = \d_{lj}.
$$
Therefore, if $\psi$ is a smooth function in $N \g_\e$, it follows
that in the above coordinates $(s,y)$
\begin{equation*}
    \D_{\hat{g}_\e} \psi = \pa^2_{s s} \psi + \pa^2_{j j} \psi,
    \qquad \quad \hbox{ in } N \g_\e.
\end{equation*}
In the case $\psi(s, y) = e^{-i\hat{f}s}u(s, y)$, for
$\hat{f}=\mathcal{A}\hat{h}^\s$ (see \eqref{eq:ovhovk}) and for $u$
real, we have clearly that
\begin{equation*}
    \D_{\hat{g}_\e} \psi = e^{-i\hat{f}s}\pa^2_{ss} u
    -2i\hat{f}e^{-i\hat{f}s}\pa_su-\hat{f}^2e^{-i\hat{f}s}u
    +e^{-i\hat{f}s}\pa^2_{j j} u.
\end{equation*}
Similarly to \eqref{eq:coarea} one easily finds
\begin{equation}\label{eq:coarea2}
    \int_{N \g_\e} f dV_{\hat{g}_\e} = \int_{\g_\e} \left( \int_{N \g_\e(s)}
    f(y) d y \right) d s.
\end{equation}

\subsection{Localizing the problem to a subset of the normal bundle
$N \g_\e$}\label{subsec:loc}

We next exploit the exponential decay of solutions (or approximate
solutions) away from $\g_\e$ to reduce \eqref{eq:new} from the whole
scaled  manifold $M_\e$ to the normal bundle $N \g_\e$: this step of
the proof follows closely a procedure in \cite{dkw}. We first define
a smooth non-increasing cutoff function $\ov{\eta} : \R \to \R$
satisfying
\begin{equation*}
     \left\{%
\begin{array}{ll}
    \ov{\eta}(t) = 1, & \hbox{ for } t \leq 0; \\
    \ov{\eta}(t) = 0, & \hbox{ for } t \geq 1; \\
    \ov{\eta}(t) \in [0,1], & \hbox{ for every } t \in \R. \\
\end{array}%
\right.
\end{equation*}
Next, if $(s,y)$ are the coordinates introduced above in $N \g_\e$,
and if $\Phi(\e s)$ is a section of $N \g$, using the notation of
Subsection 3.1 in \cite{mmm1} we define
$$
 z = y - \Phi(\e s).
$$
We will assume throughout the paper that $\Phi$ satisfies the
following bounds
\begin{equation}\label{eq:bdPhi}
    \|\Phi\|_{\infty} + \|\Phi'\|_{\infty} + \e \|\Phi''\|_{\infty} \leq C \e
\end{equation}
for some fixed constant $C > 0$. Next, for a small $\ov{\d}
> 0$ and for a smooth function $K(\e s) > 0$, both to be determined below
and letting $h, k : [0,L] \to \R$
be as in the introduction, we set
\begin{equation}\label{eq:defetae}
    \tilde{\psi}_{0,\e} = \ov{\eta}_\e(s,z) \psi_{0,\e} := \ov{\eta}\left( K(\e s)
    \left( |z| - \frac{\e^{-\ov{\d}}}{K(\e s)} \right) \right)
   e^{- i \frac{\tilde{f}(\e s)}{\e}} h(\e s) U(k(\e s) z),
\end{equation}
where $\tilde{f}$ (to be defined later) is close to the function $f$
(also defined in the introduction). For $\t \in (0,1)$, we let $S_\e
: C^{2,\t}(M_\e) \to C^\t(M_\e)$ be the operator
\begin{equation}\label{eq:defse}
    S_\e (\psi) = - \D_{g_\e} \psi + V(\e x) \psi -
  |\psi|^{p-1} \psi \qquad \quad \hbox{ in } M_\e.
\end{equation}
If we let $\tilde{\psi}_\e$ denote an approximate solution of
\eqref{eq:new} (we will take later $\tilde{\psi}_\e$ equal to
$\tilde{\psi}_{0,\e}$, with some small correction), then setting
$\psi = \tilde{\psi}_\e + \tilde{\phi}$, we have $S_\e(\psi) = 0$ if
and only if
$$
  L_\e (\tilde{\phi}) = S_\e(\tilde{\psi}_\e) + N_\e(\tilde{\phi})
\qquad \quad \hbox{ in } M_\e,
$$
where $L_\e(\tilde{\phi})$ stands for the linear correction in
$\tilde{\phi}$, namely
\begin{equation}\label{eq:defLe}
    L_\e(\tilde{\phi}) = - \D_{g_\e} \tilde{\phi} + V(\e x) \tilde{\phi} -
  |\tilde{\psi}_\e|^{p-1} \tilde{\phi} - (p-1) |\tilde{\psi}_\e|^{p-3}
  \tilde{\psi}_\e \Re (\tilde{\psi}_\e \ov{\tilde{\phi}}) \qquad \quad \hbox{ in } M_\e,
\end{equation}
and where the nonlinear operator $N_\e(\tilde{\phi})$ is defined
as
\begin{equation}\label{eq:defNe}
    N_\e(\tilde{\phi}) = |\tilde{\psi}_\e + \tilde{\phi}|^{p-1} (\tilde{\psi}_\e
  + \tilde{\phi}) - |\tilde{\psi}_\e|^{p-1} \tilde{\psi}_\e - |\tilde{\psi}_\e|^{p-1}
  \tilde{\phi} - (p-1) |\tilde{\psi}_\e|^{p-3} \tilde{\psi}_\e
  \Re (\tilde{\psi}_\e \ov{\tilde{\phi}}).
\end{equation}
Then, in the coordinates $(s,z)$ we can write
$$
  \tilde{\phi} = \ov{\eta}_\e(z) \phi + \varphi
$$
where, with an abuse of notation, we assume $\phi$ defined on $N
\g_\e$ (through the exponential map normal to $\g_\e$) and where
the correction $\varphi$ is defined on the whole $M_\e$. In this
way we need to solve the equation
\begin{equation}\label{eq:sol}
    L_\e (\ov{\eta}_\e(z) {\phi}) + L_\e (\varphi) = S_\e(\tilde{\psi}_\e) +
    N_\e (\ov{\eta}_\e(z) {\phi} + \varphi) \qquad \hbox{ in } M_\e.
\end{equation}
We will require $\phi$ to be supported in a cylindrical-shaped
region in $N \g_\e$ centered around the zero section. For technical
reasons,  convenient for proving the results in the next subsection,
we define
\begin{equation}\label{eq:tildeDe}
    \tilde{D}_\e = \left\{ (s,z) \in N \g_\e \; : \; |z| \leq
\frac{\e^{-\ov{\d}}+1}{K(\e s)}  \right\},
\end{equation}
and then the subspace of functions in $N \g_\e$
\begin{equation*}
    H_{\tilde{D}_\e} = \left\{ u \in L^2(N \g_\e; \C) \; : \;
    u \hbox{ is supported in } \tilde{D}_\e \right\}.
\end{equation*}
Using elementary computations, we see that \eqref{eq:sol} is
satisfied if (tautologically) the following two conditions are
imposed
\begin{eqnarray}\label{eq:eqiner} \nonumber
  L_\e(\phi) & = & \left[ S_\e(\tilde{\psi}_\e)
    + N_\e(\ov{\eta}_\e(z) \phi + \varphi) \right] \\ & + &
     |\tilde{\psi}_\e|^{p-1} \varphi + (p-1)
    |\tilde{\psi}_\e|^{p-3} \tilde{\psi}_\e \Re(\tilde{\psi}_\e \ov{\varphi})
    \quad \hbox{ in } \tilde{D}_\e, \qquad \phi \in H_{\tilde{D}_\e};
\end{eqnarray}
\begin{equation}\label{eq:eqphi}
  \mathcal{L}_{\tilde{\psi}_\e} \varphi
  = (1 - \ov{\eta}_\e(z)) \left[ S_\e(\tilde{\psi}_\e) + N_\e(\ov{\eta}_\e(z) \phi
   + \varphi) \right] + 2 \, \n_{g_\e} \ov{\eta}_\e(z) \cdot \n_{g_\e} \phi +
   \D_{g_\e} \ov{\eta}_\e(z) \phi \qquad \hbox{ in } M_\e,
\end{equation}
where
\begin{equation}\label{eq:mathclphie}
    \mathcal{L}_{\tilde{\psi}_\e} \varphi = - \D_{g_\e} \varphi + V(\e x)
  \varphi - (1 - \ov{\eta}_\e(z)) \left[ |\tilde{\psi}_\e|^{p-1} \varphi +
  (p-1) |\tilde{\psi}_\e|^{p-2} \tilde{\psi}_\e \Re (\tilde{\psi}_\e
  \ov{\varphi}) \right].
\end{equation}

\

\noindent We have next an existence result for equation
\eqref{eq:eqphi}: in order to state it we need to introduce some
notation. For a regular periodic function $\mathfrak{p} : [0,L] \to
\R$, for $m \in \N$ and $\t \in (0,1)$ we define the weighted norms
\begin{equation}\label{eq:wnmak1}
    \|\varphi\|_{C^{m,\t}_{\mathfrak{p}}} = \sup_{x \in \tilde{D}_\e} \left[ e^{\mathfrak{p}(\e
  s)|z|} \|\varphi\|_{C^{m,\t}(B_1(x))} \right], \quad \qquad x = (s,z).
\end{equation}
We also recall the definition of $k(\e s)$ in \eqref{eq:ovhovks}.

\begin{pro}\label{p:exvarphi}
Let $k_2(\ov{s}) < k_1(\ov{s}) < k_0(\ov{s}), K(\ov{s})$ be smooth
positive $L$-periodic functions in $\ov{s}$, and $\t \in (0,1)$.
Then, if $V(\ov{s}), K^2(\ov{s})
> k_2^2(\ov{s})$ and if $\|\tilde{\psi}_\e\|_{C^ \t_{k_0}},
\|S_\e(\tilde{\psi}_\e)\|_{C^ \t_{k_0}}
\leq 1$, there exists a positive constant $C$ depending on $\ov{\d},
\t, k, k_0, k_1$ and $k_2$ such that given any $\phi$ with
$\|\phi\|_{C^{1,\t}_{k_1}} \leq 1$ problem \eqref{eq:eqphi} has a
unique solution $\varphi(\phi)$ whose restriction to $\tilde{D}_\e$
satisfies
\begin{equation}\label{eq:qqqqq1}
    \|\varphi(\phi)\|_{C^\t_{-k_2}} \leq C \left( e^{- \e^{-\ov{\d}} \inf \frac{k_2+k_0}{K}
} \|S_\e(\tilde{\psi}_\e)\|_{C^ \t_{k_0}} + e^{- \e^{-\ov{\d}} \inf
\frac{k_1+k_2}{K}} \|\phi\|_{C^{1,\t}_{k_1}} \right).
\end{equation}
Moreover, if $\tilde{\psi}^1_\e, \tilde{\psi}^2_\e$ satisfy
$\|S_\e(\tilde{\psi}^j_\e)\|_{C^ \t_{k_0}} \leq 1$, $j=1, 2$,  if
$\|\phi_j\|_{C^{1,\t}_{k_1}} \leq 1$, $j = 1, 2,$ and if
$\varphi_j(\phi_j)$, $j = 1, 2$, are the corresponding solutions,
for the restrictions to $\tilde{D}_\e$ we also have
\begin{equation}\label{eq:qqqqq2}
    \|\varphi(\phi_1) - \varphi(\phi_2)\|_{C^\t_{-k_2}} \leq C
\left( e^{- \e^{-\ov{\d}} \inf \frac{k_2+k_0}{K}} \|S_\e(\tilde{\psi}^1_\e)
- S_\e(\tilde{\psi}^2_\e)\|_{C^ \t_{k_0}} + e^{- \e^{-\ov{\d}} \inf
\frac{k_1+k_2}{K}} \|\phi_1 - \phi_2\|_{C^{1,\t}_{k_1}} \right).
\end{equation}
\end{pro}

\begin{rem}\label{r:posexp} (a) The choice of the norm in
\eqref{eq:wnmak1} is done for considering functions which grow at
most like $e^{-\mathfrak{p}(\e s) |z|}$, and in particular functions
which decay at infinity if $\mathfrak{p}$ is positive. In the
left-hand side of \eqref{eq:qqqqq1} we have a negative exponent,
representing the fact that $\varphi$ can grow as $|z|$ increases.
However (we will take later $k_0, k_1, k_2, K$ very close), the
coefficients in the right-hand side are so tiny that $\varphi$ is
everywhere small in $\tilde{D}_\e$, and indeed with an even smaller
bound for $|z|$ close to zero. This reflects the fact that the
support of the right-hand side in \eqref{eq:eqphi} is
$\frac{\e^{-\ov{\d}}}{K(\e s)} < |z| < \frac{\e^{-\ov{\d}}+1}{K(\e
s)} $, so $\varphi$ should decay away from this set.

(b) We introduced the functions $k_0, k_1$ and $k_2$ for technical
reasons, since we want to allow some flexibility for the
(exponential) decay rate in $|z|$.
\end{rem}

\begin{pf} We prove the result only when the manifold $M$ in
\eqref{eq:pe} is compact. For the modifications needed for $M =
\R^n$ see Remark \ref{r:imprregvp} $(b)$.

Consider a smooth non-decreasing cutoff function $\chi : [0,1] \to
[0,1]$ satisfying
$$
  \left\{
    \begin{array}{ll}
      \chi(t) = 0 & \hbox{ for } t \leq \frac 14; \\
      \chi(t) = t & \hbox{ for } t \geq \frac 34; \\
      0 \leq\chi'(t) \leq 4  & \hbox{ for all } t; \\
      0 \leq\chi''(t) \leq 16  & \hbox{ for all } t.
    \end{array}
  \right.
$$
Next, given a large constant $\mathcal{B}$ (to be specified later)
depending only on $V$ and $k_2$, we define
$\tilde{\chi}(\ov{s},|z|)$ as
$$
  \tilde{\chi}(\ov{s},|z|) = \left\{
    \begin{array}{ll}
       \mathcal{B} \chi\left( \frac{|z|}{\mathcal{B}} \right)
& \hbox{ for } |z| \leq \mathcal{B}; \\
      |z| & \hbox{ for } \mathcal{B} \leq |z| \leq \frac{\e^{-\ov{\d}}}{k_2(\ov{s})} -1; \\
      |z| - \frac{\e^{-\ov{\d}}}{k_2(\ov{s})} - \frac 12 - 2 \chi \left(
       |z| - \frac{\e^{-\ov{\d}}}{k_2(\ov{s})} - \frac 12 \right) & \hbox{ for }
      \frac{\e^{-\ov{\d}}}{k_2(\ov{s})} -1 \leq |z| \leq \frac{\e^{-\ov{\d}}}{k_2(\ov{s})}
      + 1; \\ 2 \frac{\e^{-\ov{\d}}}{k_2(\ov{s})} - |z| & \hbox{ for }
\frac{\e^{-\ov{\d}}}{k_2(\ov{s})} + 1 \leq |z| \leq 2 \frac{\e^{-\ov{\d}}}{k_2(\ov{s})} - \mathcal{B}; \\
      \mathcal{B} \chi \big( \frac{2 \frac{\e^{-\ov{\d}}}{k_2(\ov{s})} - |z|}{\mathcal{B}}
\big) & \hbox{ for }  2 \frac{\e^{-\ov{\d}}}{k_2(\ov{s})} - \mathcal{B} \leq |z|
\leq 2 \frac{\e^{-\ov{\d}}}{k_2(\ov{s})};
     \\ 0 & \hbox{ for } |z| \geq 2 \frac{\e^{-\ov{\d}}}{k_2(\ov{s})}.
    \end{array}
  \right.
$$
By our choice of $\chi$, the function $\tilde{\chi}$ satisfies the
following inequalities (where, here, the gradient and the Laplacian
are taken with respect to the  Euclidean metric)
\begin{equation*}
    |\n_z \tilde{\chi}| \leq 1; \qquad \qquad \D_z \tilde{\chi} \leq
\frac{16 + 4 (n-2)}{\mathcal{B}}.
\end{equation*}
Using the above coordinates $(s,z)$, we define next the barrier
function $\mathfrak{u} : M_\e \to \R$ as
$$
  \mathfrak{u}(s,z) = e^{k_2(\e s) \tilde{\chi}(\e s, z)} \qquad \qquad
\hbox{ for } |z| \leq 2 \frac{\e^{-\ov{\d}}}{k_2(\e s)},
$$
and we extend $\mathfrak{u}$ identically equal to $1$ elsewhere. By
our choice of $\tilde{\chi}$, this function is indeed smooth and
strictly positive on the whole $M_\e$. We consider next the linear
equation (motivated by \eqref{eq:mathclphie})
$$
  \mathcal{L}_{\tilde{\psi}_\e} \varphi = \vartheta \qquad \quad \hbox{ on }
M_\e,
$$
where $\vartheta : M_\e \to \R$ is H\"older continuous (with
$supp(\vartheta) \subset\subset \tilde{D}_\e$, see
\eqref{eq:equivalence} below). Since the operator
$\mathcal{L}_{\tilde{\psi}_\e}$ is uniformly elliptic, the latter
equation is (uniquely) solvable, and we would like next to derive
some pointwise estimates on its solutions. To this aim we define
$$
  v(x) = \frac{\varphi(x)}{\mathfrak{u}(x)}, \qquad \quad x \in M_\e.
$$
With this notation, we have that
$$
  \mathfrak{u} \mathcal{L}_{\tilde{\psi}_\e} v - v \D_{g_\e} \mathfrak{u} - 2 \n_{g_\e}
  v \cdot \n_{g_\e} \mathfrak{u} = \vartheta \qquad \quad \hbox{ on } M_\e.
$$
Using the expression of the metric coefficients in the coordinates
$(s,z)$, see Lemma 2.1 in \cite{mmm1}, \eqref{eq:bdPhi} and the
properties of the cutoff function $\tilde{\chi}$, one easily checks
that
$$
  \D_{g_\e} \mathfrak{u} \left\{
           \begin{array}{ll}
             \leq ( k_2(\ov{s})^2 + o_\mathcal{B}(1) + o_\e(1) ) \mathfrak{u}  & \hbox{ for } |z| \leq
    \frac{\e^{-\ov{\d}}}{k_2(\ov{s})}, \\
             = 0 & \hbox{ elsewhere},
           \end{array}
         \right.
$$
where $o_\e(1) \to 0$ as $\e \to 0$ and $o_\mathcal{B}(1) \to 0$ as
$\mathcal{B} \to + \infty$. Therefore we obtain that the function
$v$ satisfies
$$
  \left\{
    \begin{array}{ll}
      \left| (\mathcal{L}_{\tilde{\psi}_\e} - k_2^2(\ov{s}) + o_\mathcal{B}(1) + o_\e(1)) v
\right| \leq \frac{|\vartheta|}{\mathfrak{u}} & \hbox{ for } |z| \leq 2 \frac{\e^{-\ov{\d}}}{k_2(\ov{s})},  \\
     \mathcal{L}_{\tilde{\psi}_\e} v = \frac{\vartheta}{\mathfrak{u}}  & \hbox{
elsewhere}.
    \end{array}
  \right.
$$
Since we assumed $V(\ov{s}) > k_2(\ov{s})^2$, we obtain that $V -
k_2(\ov{s})^2 + o_\mathcal{B}(1) + o_\e(1)$ is strictly positive
(provided $\mathcal{B}$ is sufficiently large and $\e$ sufficiently
small) for $|z| \leq 2 \frac{\e^{-\ov{\d}}}{k_2(\ov{s})}$, and hence
the function $v$ satisfies a uniformly elliptic equation with a
non-negative coefficient in the zero-th order term with right-hand
side given by $\frac{\vartheta}{\mathfrak{u}}$. Therefore from the
maximum principle we derive the estimate
$$
  \max_{M_\e} |v| \leq C \max_{M_\e} \frac{|\vartheta|}{\mathfrak{u}},
$$
where $C$ depends on the uniform lower bound of the above
coefficient. The latter estimate clearly implies
$$
  |\varphi(x)| \leq C \mathfrak{u}(x) \max_{M_\e} \frac{\vartheta}{\mathfrak{u}}, \qquad \quad
  \hbox{ for every } x \in M_\e.
$$
We define next the weighted norm
$$
  \|\varphi\|_{m,\t,\mathfrak{u}} := \sup_{x \in M_\e} \left\| \frac{\varphi}{\mathfrak{u}}
\right\|_{C^{m,\t}(B_1(x))},
$$
which is equivalent (with constants depending on $\mathcal{B}$ only)
to $\| \cdot \|_{C^{m,\t}_{-k_2}}$   on the set $\tilde{D}_\e$.
Using the explicit form of the function $\mathfrak{u}$ and standard
elliptic regularity estimates one can improve the latter inequality
to
\begin{equation}\label{eq:estfinlin}
    \|\varphi\|_{2,\t,\mathfrak{u}} \leq C \|\vartheta\|_{0,\t,\mathfrak{u}}.
\end{equation}
The proof of the proposition will now follow from this linear
estimate and the contraction mapping theorem: in fact,  defining
$$
  G_{\phi,\e}(\varphi) = (1 - \ov{\eta}_\e(s,z)) \left[ S_\e(\tilde{\psi}_\e)
 + N_\e(\ov{\eta}_\e(z) \phi + \varphi) \right] + 2 \, \n_{g_\e}
 \ov{\eta}_\e(z) \cdot \n_{g_\e} \phi + \D_{g_\e} \ov{\eta}_\e(z) \phi,
$$
equation \eqref{eq:eqphi} is equivalent to
\begin{equation}\label{eq:equivalence}
    \varphi = \mathcal{L}_{\tilde{\psi}_\e}^{-1} G_{\phi,\e}(\varphi).
\end{equation}
First of all, notice that $\mathcal{L}_{\tilde{\psi}_\e}$ is
invertible since we are assuming $M$ (and hence $M_\e$) to be
compact, see the beginning of the proof. Secondly, to apply
\eqref{eq:estfinlin}, we need to estimate
$\|G_{\phi,\e}(\varphi)\|_{0,\t,\mathfrak{u}}$, together with its
Lipschitz dependence in $\varphi$: our goal indeed is to apply the
contraction mapping theorem.

Let us consider for instance the term $(1 - \ov{\eta}_\e(s,z))
S_\e(\tilde{\psi}_\e)$. Using the fact that $(1-\ov{\eta}_\e)$ is
zero for $|z| \leq \frac{\e^{-\ov{\d}}}{K(\e s)}$, that
$S_\e(\tilde{\psi}_\e)$ is zero for $|z| \geq
\frac{\e^{-\ov{\d}}+1}{K(\e s)} $ and that $k_2 < K$, we obtain
\begin{equation}\label{eq:estfirst}
    \left\| (1 - \ov{\eta}_\e) S_\e(\tilde{\psi}_\e) \right\|_{0,\t,\mathfrak{u}} \leq
  C \sup_{\{\frac{\e^{-\ov{\d}}}{K(\e s)} \leq |z| \leq
\frac{\e^{-\ov{\d}}+1}{K(\e s)} \}} \left\|
\frac{S_\e(\tilde{\psi}_\e)}{\mathfrak{u}} \right\|_{C^\t(B_1(z))}
\leq C e^{- \e^{-\ov{\d}} \inf \frac{k_0 + k_2}{K}}
\|S_\e(\tilde{\psi}_\e)\|_{C^ \t_{k_0}}.
\end{equation}
Now to estimate the remaining terms of $G_{\phi,\e}$ we notice that
$$
  |N_\e(\ov{\eta}_\e(z) \phi + \varphi)| \leq \left\{
                                                \begin{array}{ll}
                                                  C |\tilde{\psi}_\e|^{p-2}
  |\ov{\eta}_\e(z) \phi + \varphi|^2 & \hbox{ if } |\ov{\eta}_\e(z) \phi + \varphi| \leq
  |\tilde{\psi}_\e|; \\
                                                  |\ov{\eta}_\e(z) \phi + \varphi|^p &
  \hbox{ otherwise.}
                                                \end{array}
                                              \right.
$$
Since $p > 1$, we can find a number $\zeta \in (0,1)$ such that $p -
2 + 1 - \zeta > 0$, so the last formula implies
\begin{equation}\label{eq:Ne1}
    |N_\e(\ov{\eta}_\e(z) \phi + \varphi)| \leq C \left(
    |\tilde{\psi}_\e|^{p-1-\zeta} (|\ov{\eta}_\e(z) \phi|^\zeta + |\varphi|^\zeta)
   (|\ov{\eta}_\e(z) \phi| + |\varphi|) + |\ov{\eta}_\e(z) \phi|^p + |\varphi|^p \right).
\end{equation}
Using the fact that $\|\tilde{\psi}_\e\|_{C^ \t_{k_0}},
\|\phi\|_{C^{1,\t}_{k_1}} \leq 1$ and reasoning as for
\eqref{eq:estfirst}, after some computations we deduce (assuming
$\|\varphi\|_\infty \leq 1$, which will be verified later)
\begin{eqnarray}\label{eq:contrvp1} \nonumber
  \|G_{\phi,\e}(\varphi)\|_{0,\t,\mathfrak{u}}  & \leq & C \left( e^{- \e^{-\ov{\d}} \inf
 \frac{k_2+k_0}{K}} \|S_\e(\tilde{\psi}_\e)\|_{C^\t_{k_0}} + e^{- \e^{-\ov{\d}} \inf
 \frac{k_2+k_1}{K}} \|\phi\|_{C^{1,\t}_{k_1}} \right. \\ & + & \left.
 e^{- \e^{-\ov{\d}} \inf \frac{p k_1 + k_2}{K}}
\|\phi\|_{C^{0,\t}_{k_1}} + e^{- \e^{-\ov{\d}} \inf \frac{(p-1)k_0+k_1+k_2}{K}}
\|\phi\|_{C^{0,\t}_{k_1}} \right) \\
  & + & C \left(  e^{- \e^{-\ov{\d}} (p-1-\zeta) \inf \frac{k_0}{K}} +
\|\varphi\|_{\infty}^{p-1} \right) \|\varphi\|_{0,\t,\mathfrak{u}}. \nonumber
\end{eqnarray}
Similarly, for two functions $\varphi_1, \varphi_2$ with
$\|\varphi_1\|_\infty, \|\varphi_2\|_\infty \leq 1$ and with finite
$\|\cdot\|_{0,\t,\mathfrak{u}}$  norm we have
\begin{equation}\label{eq:contrvp2}
  \|G_{\phi,\e}(\varphi_1)-G_{\phi,\e}(\varphi_2)\|_{0,\t,\mathfrak{u}} \leq
  \left(  e^{- \e^{-\ov{\d}} (p-1-\zeta) \inf \frac{k_0}{K}} +
\|\varphi_1\|_{\infty}^{p-1} + \|\varphi_2\|_{\infty}^{p-1}\right) \|\varphi_1 -
\varphi_2\|_{0,\t,\mathfrak{u}}.
\end{equation}
We now consider the map $\varphi \mapsto G_{\phi,\e}(\varphi)$ in
the set
$$
  \mathfrak{B} = \left\{ \varphi \; : \; \|\varphi\|_{0,\t,\mathfrak{u}} \leq C_1
  \left( e^{- \e^{-\ov{\d}} \inf \frac{k_2+k_0}{K}} \|S_\e(\tilde{\psi}_\e)\|_{C^\t_{k_0}}
+ e^{-\e^{-\ov{\d}} \inf \frac{k_2+k_1}{K}} \|\phi\|_{C^{1,\t}_{k_1}} \right)
\right\},
$$
where $C_1$ is a sufficiently large positive constant: notice that
if $\varphi \in \mathfrak{B}$ then $\|\varphi\|_\infty = o_\e(1)$.
From \eqref{eq:contrvp1}, \eqref{eq:contrvp2}  it then follows that
this map is a contraction from $\mathfrak{B}$ into itself, endowed
with the above norm, and therefore a solution $\varphi$ exists as a
fixed point of $G_{\phi,\e}$. The fact that $k_2 < K$ implies that
the norm $\|\cdot\|_{C^\t_{-k_2}}$ is equivalent to
$\|\cdot\|_{0,\t,\mathfrak{u}}$ in $\tilde{D}_\e$ (see also the
comments in Remark \ref{r:posexp}), so we obtain  \eqref{eq:qqqqq1}.
A similar reasoning, still based on regularity theory and elementary
inequalities, also yields \eqref{eq:qqqqq2}.
\end{pf}

\begin{rem}\label{r:imprregvp} (a) From elliptic regularity
theory it follows that in \eqref{eq:qqqqq1}-\eqref{eq:qqqqq2} the
norm $\|\cdot\|_{C^\t_{-k_2}}$ can be replaced by the stronger one
$\|\cdot\|_{C^{2,\t}_{-k_2}}$, yielding an estimate in the
 norm $\|\cdot\|_{C^{\t'}_{-k_2}}$ for any $\t' \in
(\t,1)$.

(b) In the case $M = \R^n$ the above proof needs to be slightly
modified: in fact the invertibility of
$\mathcal{L}_{\tilde{\psi}_\e}$ will be guaranteed provided we work
in an appropriate class of functions $Y$ decaying exponentially at
infinity. To guarantee this condition, we can vary the form of the
barrier function $\mathfrak{u}$ in order that it both remains a
super-solution of $\mathcal{L}_{\tilde{\psi}_\e} = 0$ and decays
exponentially to zero at infinity. This is indeed possible using the
uniform positive lower bound on $V$, see \eqref{eq:VV}: we omit the
details of this construction.
\end{rem}

\

\noindent As a consequence of Proposition \ref{p:exvarphi}, we
obtain that solvability of \eqref{eq:pe} is equivalent to that of
\eqref{eq:eqiner}.

\begin{pro}\label{p:sec2final} Suppose the  assumptions of
Proposition \ref{p:exvarphi} hold,  and consider the corresponding
$\varphi = \varphi(\phi)$. Then $\psi = \tilde{\psi}_\e +
\ov{\eta}_\e(z) \phi + \varphi(\phi)$ solves \eqref{eq:new} if and
only if $\phi \in H_{\tilde{D}_\e}$ satisfies
\begin{equation}\label{eq:LetildeSe}
    L_\e(\phi) = \tilde{S}_\e(\phi) \qquad \qquad \hbox{ in } \tilde{D}_\e,
\end{equation}
with \begin{equation}\label{eq:tildeSe}
    \tilde{S}_\e(\phi) = S_\e(\tilde{\psi}_\e)
    +  N_\e(\ov{\eta}_\e(z) \phi + \varphi(\phi)) +
  |\tilde{\psi}_\e|^{p-1} \varphi(\phi) + (p-1)  |\tilde{\psi}_\e|^{p-3}
  \tilde{\psi}_\e \Re (\tilde{\psi}_\e \ov{\varphi(\phi)}) \quad \hbox{ in }
\tilde{D}_\e,
\end{equation}
where $\ov{\eta}_\e, S_\e$ and $N_\e$ are given in
\eqref{eq:defetae}, \eqref{eq:defse} and \eqref{eq:defNe}
respectively.
\end{pro}

\subsection{Construction of an approximate kernel for $L_\e$}\label{ss:appker}

We perform here some preliminary analysis useful to understand the
spectral properties of $L_\e$.  More precisely in this subsection we
consider a {\em model} case, when the domain $\tilde{D}_\e$ (see
\eqref{eq:tildeDe}) is replaced by $[0,L/\e] \times \R^{n-1}$ and
the profile of approximate solutions is independent of the variable
$s$ (only the phase varies, periodically in $s$). As in formula
\eqref{eq:ovhovk}, we consider positive constants $\hat{V}$,
$\hat{h}$, $\hat{k}$ satisfying
\begin{equation}\label{eq:ovhovk2}
    \hat{h} = \left( \hat{f}^2 + \hat{V} \right)^{\frac{1}{p-1}};
    \qquad \qquad \hat{k} = \left( \hat{f}^2 + \hat{V} \right)^{\frac 12}.
\end{equation}
Our goal is to study the  following eigenvalue problem, which models
our linearized equation
\begin{equation*}
  \hat{L}_\e u = \l u \quad \hbox{ in }
      [0,L/\e] \times \R^{n-1};
\end{equation*}
\begin{equation}\label{eq:hatLe}
    \hat{L}_\e u = - \D_{\hat{g}_\e} u + \hat{V} u - \hat{h}^{p-1}
      U(\hat{k} y)^{p-1} u - (p-1) \hat{h}^{p-1} U(\hat{k} y)^{p-1} e^{- i \hat{f} s}
      \Re(e^{- i \hat{f} s} \ov{u}),
\end{equation}
and in particular we would like to characterize the small
eigenvalues and the corresponding eigenfunctions. First of all we
can write $u$ as
$$
  u = e^{- i \hat{f} s} (u_r + i u_i),
$$
for some real $u_r$ and $u_i$. With this notation, we are reduced to
study the coupled system
\begin{equation*}
  \left\{
    \begin{array}{ll}
      - \D_{\hat{g}_\e} u_r + (\hat{V} + \hat{f}^2) u_r - p \hat{h}^{p-1}
      U(\hat{k} y)^{p-1} u_r - 2 \hat{f} \frac{\pa u_i}{\pa s} = \l u_r & \hbox{ in }
     [0,L/\e] \times \R^{n-1}; \\
      - \D_{\hat{g}_\e} u_i + (\hat{V} + \hat{f}^2) u_i - \hat{h}^{p-1}
      U(\hat{k} y)^{p-1} u_i + 2 \hat{f} \frac{\pa u_r}{\pa s} = \l u_i
& \hbox{ in }
     [0,L/\e] \times \R^{n-1}.
    \end{array}
  \right.
  \end{equation*}
Making the change of variables $y \mapsto \hat{k} y$ and using
\eqref{eq:ovhovk2}, we are reduced to
\begin{equation}\label{eq:coupledsys}
  \left\{
    \begin{array}{ll}
      - \frac{1}{\hat{k}^2} \frac{\pa^2 u_r}{\pa s^2}- \D_{y} u_r +  u_r - p
      U(y)^{p-1} u_r - \frac{2 \hat{f}}{\hat{k}^2} \frac{\pa u_i}{\pa s} = \frac{\l}{\hat{k}^2} u_r
    & \hbox{ in } [0,L/\e] \times \R^{n-1};  \\
      - \frac{1}{\hat{k}^2} \frac{\pa^2 u_i}{\pa s^2}- \D_{y} u_i + u_i -
      U(y)^{p-1} u_i + \frac{2 \hat{f}}{\hat{k}^2} \frac{\pa u_r}{\pa s} = \frac{\l}{\hat{k}^2} u_i
     & \hbox{ in } [0,L/\e] \times \R^{n-1}.
    \end{array}
  \right.
\end{equation}
It is now convenient to use a Fourier decomposition in $s$ of $u_r$
and $u_i$, writing
$$
  u_r = \sum_j \left(\cos \left( \frac{2\pi \e j s}{L} \right) u_{r,c,j}(y) +
  \sin \left( \frac{2\pi \e j s}{L} \right) u_{r,s,j}(y) \right), \qquad s \in
  [0,L/\e], y \in \R^{n-1};
$$
$$
  u_i = \sum_j \left(\cos \left( \frac{2\pi \e j s}{L} \right) u_{i,c,j}(y)
  +\sin \left( \frac{2\pi \e j s}{L} \right)
u_{i,s,j}(y) \right), \qquad s \in [0,L/\e], y \in \R^{n-1}.
$$
In this way the functions $u_{r,c,j}, u_{r,s,j}, u_{i,c,j},
u_{i,s,j}$ satisfy the following systems of equations
$$
  \left\{
    \begin{array}{ll}
      - \D_{y} u_{r,c,j} + \left(1 + \frac{4 \pi^2 \e^2 j^2}{L^2 \hat{k}^2}
      \right) u_{r,c,j} - p U(y)^{p-1} u_{r,c,j} - \frac{4 \pi \hat{f}
      \e j}{L \hat{k}^2} u_{i,s,j} = \frac{\l}{\hat{k}^2} u_{r,c,j} & \hbox{ in } \R^{n-1}; \\
      - \D_{y} u_{i,s,j} + \left(1 + \frac{4 \pi^2 \e^2 j^2}{L^2 \hat{k}^2} \right) u_{i,s,j} -
      U(y)^{p-1} u_{i,s,j} - \frac{4 \pi \hat{f} \e j}{L \hat{k}^2} u_{r,c,j} =
      \frac{\l}{\hat{k}^2} u_{i,s,j} & \hbox{ in } \R^{n-1},
    \end{array}
  \right.
$$
$$
  \left\{
    \begin{array}{ll}
      - \D_{y} u_{r,s,j} + \left(1 + \frac{4 \pi^2 \e^2 j^2}{L^2 \hat{k}^2} \right)
      u_{r,s,j} - p U(y)^{p-1} u_{r,s,j} + \frac{4 \pi \hat{f} \e j}{L \hat{k}^2}
      u_{i,c,j} = \frac{\l}{\hat{k}^2} u_{r,s,j} & \hbox{ in } \R^{n-1}; \\
      - \D_{y} u_{i,c,j} + \left(1 + \frac{4 \pi^2 \e^2 j^2}{L^2 \hat{k}^2} \right) u_{i,c,j} -
      U(y)^{p-1} u_{i,c,j} + \frac{4 \pi \hat{f} \e j}{L \hat{k}^2} u_{r,s,j} =
      \frac{\l}{\hat{k}^2} u_{i,c,j} & \hbox{ in } \R^{n-1}.
    \end{array}
  \right.
$$
If we set $\frac{2 \pi \e j}{L \hat{k}} = \a$,  $\frac{2
\hat{f}}{\hat{k}} = \mu$ and $\tilde{\l} = \frac{\l}{\hat{k}^2}$
then the latter two systems are equivalent to the following one
\begin{equation}\label{eq:systuv}
  \left\{
    \begin{array}{ll}
      - \D_{y} u + (1 + \a^2) u - p U(y)^{p-1} u + \mu \a v = \tilde{\l} u & \hbox{ in }
    \R^{n-1}; \\
      - \D_{y} v + (1 + \a^2) v - U(y)^{p-1} v + \mu \a u = \tilde{\l} v & \hbox{ in }
    \R^{n-1}.
    \end{array}
  \right.
\end{equation}
The equivalence with the second system is obvious: for the first one
it is sufficient to switch the sign of the second component. We
characterize the spectrum of the last system in the next
proposition: the value of $\mu$ is fixed, while $\a$ is allowed to
vary. We remark that it is irrelevant for our purposes to take $\a$
positive or negative, since we can still switch the sign of one of
the two components.

\begin{pro}\label{p:alphamu}
Let $\eta_\a, \s_\a$ and $\t_\a$ denote the first three eigenvalues
of \eqref{eq:systuv}. Then there exists $\mu_0 > 0$ such that for
$\mu \in [0,\mu_0]$ the following properties hold

\begin{description}

  \item[a] there exists $\a_0 > 0$ such that
$\eta_\a$ is simple, increasing and differentiable in $\a$ for
$\a \in [0,\a_0]$, $\frac{\pa \eta_\a}{\pa \a}
> 0$ for $\a \in (0,\a_0]$, $\eta_0 < 0$ and $\eta_{\a_0} > 0$;

  \item[b] the eigenvalue $\s_\a$ is zero for $\a =
  0$ with multiplicity $n$,  it satisfies $\frac{\partial
\s_\a}{\partial \a} > 0$ for $\a$ small positive and stays
uniformly bounded away from zero if $\a$ stays bounded away from
zero;

  \item[c] $\t_\a$ is strictly positive and stays uniformly
  bounded away from zero for all $\a$'s;

  \item[d] the eigenfunction $u_\a$ corresponding to $\eta_\a$ is simple, radial in
  $y$, radially decreasing and depends smoothly on $\a$; for $\a
  = 0$ the eigenfunction of \eqref{eq:systuv} corresponding to
  $\eta_0 < 0$ is of the form $(\tilde{Z},0)$ with $\tilde{Z}$
  radial and radially decreasing, while those corresponding to
  $\s_0 = 0$ are linear combinations of $(\n_{y_j} U, 0)$, $j =
1,
  \dots, n-1$, and $(0,U)$;

  \item[e] let $\ov{\a}$ be the unique $\a$ for which $\eta_{\ov{\a}} = 0$
  (see {\bf a}): then the corresponding eigenfunction is of the
  form $(Z,W)$ for some radial functions $Z, W$ satisfying the
  following decay $|Z| + |W| \leq C e^{-(1+\hat{\eta})|y|}$ for
  some constants $C, \hat{\eta} > 0$.
\end{description}
\end{pro}

\begin{pf} This result is known for $\mu = 0$, see e.g. Proposition 4.2 in
\cite{mm} and Proposition 2.9 in \cite{mal}.

For $\mu \neq 0$ sufficiently small the functions $\a \mapsto
\eta_\a$, $\a \mapsto \s_\a$ and $\a \mapsto \t_\a$ will be
$C^1$-close to those corresponding to $\mu = 0$: therefore, to prove
{\bf a}-{\bf d} it is sufficient to show that $\eta_\a, \s_\a$ are
twice differentiable in $\a$ for $\a$ small, that $\frac{\pa
\eta_\a}{\pa \a} = \frac{\pa \s_\a}{\pa \a} =0$, and that
$\frac{\pa^2 \eta_\a}{\pa \a^2}, \frac{\pa^2 \s_\a}{\pa \a^2} > 0$.

We prove this statement only formally, but a rigorous proof can be
easily derived. Differentiating
\begin{equation}\label{eq:systuva}
  \left\{
    \begin{array}{ll}
      - \D_{y} u_\a + (1 + \a^2) u_\a - p U(y)^{p-1} u_\a + \mu \a v_\a
      = \eta_\a u_\a & \hbox{ in }
    \R^{n-1}; \\
      - \D_{y} v_\a + (1 + \a^2) v_\a - U(y)^{p-1} v_\a + \mu \a u_\a =
      \eta_\a v_\a & \hbox{ in }
    \R^{n-1}
    \end{array}
  \right.
\end{equation}
with respect to $\a$ we find
\begin{equation}\label{eq:systuvadiff}
  \left\{
    \begin{array}{ll}
      - \D_{y} \frac{\pa u_\a}{\pa \a} + (1 + \a^2) \frac{\pa u_\a}{\pa \a} -
      p U(y)^{p-1} \frac{\pa u_\a}{\pa \a} + \mu \a \frac{\pa v_\a}{\pa \a}
      + 2 \a u_\a + \mu v_\a = \eta_\a \frac{\pa u_\a}{\pa \a} + \frac{\pa
      \eta_\a}{\pa \a} u_\a & \hbox{ in } \R^{n-1}; \\
      - \D_{y} \frac{\pa v_\a}{\pa \a} + (1 + \a^2) \frac{\pa v_\a}{\pa \a} -
      U(y)^{p-1} \frac{\pa v_\a}{\pa \a} + \mu \a \frac{\pa u_\a}{\pa \a} + 2 \a v_\a
      + \mu u_\a = \eta_\a \frac{\pa v_\a}{\pa \a} + \frac{\pa \eta_\a}{\pa \a} v_\a & \hbox{ in }
    \R^{n-1}.
    \end{array}
  \right.
\end{equation}
To compute $\frac{\pa \eta_\a}{\pa \a}$ at $\a = 0$ it is sufficient
to multiply the first equation by $u_\a$, the second by $v_\a$, to
take the sum and integrate: if we choose $\frac{\pa u_\a}{\pa \a}$
and $\frac{\pa v_\a}{\pa \a}$ so that $\int_{\R^{n-1}} u_\a
\frac{\pa u_\a}{\pa \a} + v_\a \frac{\pa v_\a}{\pa \a} = 0$
(choosing for example $\int (u_\a^2 + v_\a^2) = 1$ for all $\a$'s),
then with an integration by parts we find that
$$
  \frac{\pa \eta_\a}{\pa \a}|_{\a=0} \int_{\R^{n-1}} u_\a^2 + v_\a^2 = 2 \mu
  \int_{\R^{n-1}} u_0 v_0.
$$
Using the fact that $v_0 = 0$, see {\bf d}, we then obtain
$\frac{\pa \eta_\a}{\pa \a}|_{\a=0} = 0$. The same argument applies
for evaluating $\frac{\pa \s_\a}{\pa \a}|_{\a=0}$, since the
eigenfunctions corresponding to $\s_0 = 0$ always have one component
vanishing.

To compute the second derivative with respect to $\a$ we
differentiate \eqref{eq:systuvadiff} once more at $\a = 0$,
obtaining
\begin{equation*}
  \left\{
    \begin{array}{ll}
      - \D_{y} \frac{\pa^2 u_\a}{\pa \a^2} + \frac{\pa^2 u_\a}{\pa \a^2} -
      p U(y)^{p-1} \frac{\pa^2 u_\a}{\pa \a^2} + 2 \mu \frac{\pa v_\a}{\pa \a}
      + 2 u_0 = \frac{\pa^2 \eta_\a}{\pa \a^2} u_0 & \hbox{ in } \R^{n-1}; \\
      - \D_{y} \frac{\pa^2 v_\a}{\pa \a^2} + \frac{\pa^2 v_\a}{\pa \a^2} -
      U(y)^{p-1} \frac{\pa^2 v_\a}{\pa \a^2} + 2 \mu \frac{\pa u_\a}{\pa \a} + 2 v_0
       = \frac{\pa^2 \eta_\a}{\pa \a^2} v_0 & \hbox{ in }
    \R^{n-1}.
    \end{array}
  \right.
\end{equation*}
As for the previous case we get
$$
  \frac{\pa^2 \eta_\a}{\pa \a^2}|_{\a=0}  =
  2 + 2 \mu \int_{\R^{n-1}} \left( u_0 \frac{\pa v_\a}{\pa \a}|_{\a = 0}
  + v_0 \frac{\pa u_\a}{\pa \a}|_{\a = 0} \right)
$$
so, using the smallness of $\mu$, the claim follows.

For the second derivative of $\s_\a$ the procedure is similar, but
notice that in this case we might obtain a multivalued function, due
to the multiplicity ($n$) of $\s_0$, see {\bf b}. However, if in the
last formula we plug in the corresponding eigenfunctions, see {\bf
d}, we still obtain a sign condition for each of the two branches of
$\s_\a$ (one of them will have multiplicity $n-1$ by the rotation
invariance of the equations).
\end{pf}

\begin{rem}\label{r:Asmall0}
Using the same argument in the previous proof one can show that
$$
\frac{\pa \eta}{\pa\a}|_{\a=\ov{\a}}=2 \ov{\a}+2\mu\int_{\R^{n-1}}Z_{\ov\a}W_{\ov\a}.
$$
\end{rem}

\begin{rem}\label{r:Asmall2}
Proposition \ref{p:alphamu} is the only result where the smallness
of the constant $\mathcal{A}$ is used, see Theorem \ref{t:main}.
Remark that $V \equiv \hat{V}$ implies $\hat{f} = \mathcal{A}
\hat{h}^\s$, and that $\mu = 2 \frac{\hat{f}}{\hat{k}}$, so the
smallness of $\mathcal{A}$ is equivalent to that of $\mu$. Notice
that by \eqref{eq:VV}, \eqref{eq:ovhovk}, when $\mathcal{A} \to 0$,
$\hat{h}$ and $\hat{k}$ stay uniformly bounded and bounded away from
zero.

We believe that dropping this smallness condition might lead to
further resonance phenomena in addition to these encountered here
(see the introduction and the last section).
\end{rem}

\begin{rem}\label{r:profeigenf} Considering \eqref{eq:systuvadiff} with $\s_\a$
replacing $\eta_\a$ and for $\a = 0$ one finds that $\mathcal{L}_r^0
\frac{\pa u_\a}{\pa \a}|_{\a = 0} = - \mu v_0$ and $\mathcal{L}_i^0
\frac{\pa v_\a}{\pa \a}|_{\a = 0} = - \mu u_0$, where
$\mathcal{L}_r^0 v = - \D_y v + v - p U(y)^{p-1} v$,
$\mathcal{L}_i^0 v = - \D_y v + v -
 U(y)^{p-1} v$. Since for $\a = 0$ we have $(u_0, v_0) = (\pa_j U, 0)$ or $(u_0,
v_0) = (0,U)$ (see {\bf d}) and
\begin{equation}\label{eq:Lri}
    \mathcal{L}_r^0 \left( - \frac{1}{p-1} U - \frac 12 \n U(y) \cdot y
    \right) = U; \qquad \qquad \mathcal{L}_i^0 (y_j U(y)) = - 2 \pa_j U,
\end{equation}
see Subsection 3.2 in \cite{mmm1}, one finds respectively that
$$
  \frac{\pa v_\a}{\pa \a}|_{\a = 0} = \frac{\mu}{2} y_j U(y) ; \qquad \qquad
  \frac{\pa u_\a}{\pa \a}|_{\a = 0} = \mu \left( \frac{1}{p-1} U +\frac 12 \n U(y)
  \cdot y \right).
$$
These expressions, together with $(59)$ in \cite{mm1} and some
integration by parts allow us to compute the explicitly $\frac{\pa^2
\s_\a}{\pa \a^2}$, whose values along the two branches are
\begin{equation*}
    \frac{\pa^2 \s_\a}{\pa \a^2} = \frac{2}{(p-1)} \left( (p-1) - 2 \mathcal{A}^2
    \th \hat{h}^{2\s-p+1} \right); \qquad \frac{\pa^2 \s_\a}{\pa \a^2} =
    \frac{2}{(p-1)} \left( (p-1) - 2 \mathcal{A}^2 \s \hat{h}^{2\s-p+1} \right).
\end{equation*}
Therefore, we find that the second derivatives of the eigenfunctions
satisfy respectively the equations
\begin{equation}\label{eq:lrj}
    \mathcal{L}_r^0 \frac{\pa^2 u_\a}{\pa \a^2} = \frac{2}{p-1} \left( (p-1) -
   2 \mathcal{A}^2 \th \hat{h}^{2\s-p+1} \right) \n_j U - 2 \n_j U -
   4 \mathcal{A}^2 \hat{h}^{2 \s - p + 1} y_j U;
\end{equation}
\begin{equation}\label{eq:lrr}
    \mathcal{L}_i^0 \frac{\pa^2 v_\a}{\pa \a^2} = \frac{2}{p-1} \left( (p-1)
  - 2 \mathcal{A}^2  \s \hat{h}^{2\s-p+1} \right) U - 2 U - 8 \mathcal{A}^2
  \hat{h}^{2 \s - p + 1} \tilde{U}.
\end{equation}

These formulas will be used crucially later on. Below, we will
denote for brevity
\begin{equation}\label{eq:frakvW}
    \hat{\mathfrak{V}}_j := \frac 12 \frac{\pa^2 u_\a}{\pa \a^2}, \quad j = 1, \dots,
    n-1; \qquad \qquad \hat{\mathfrak{W}} := \frac 12 \frac{\pa^2 v_\a}{\pa \a^2}.
\end{equation}
The factor $\frac 12$ arises in the Taylor expansion of the
eigenfunctions in $\a$, and $j$ is the index in \eqref{eq:lrj}.
\end{rem}

\no We next consider the case of variable coefficients, which can be
reduced to the previous one through a localization argument in  $s$.
To have a more accurate model for $L_\e$ the constants $\hat{k}$ and
$\hat{f}$ in \eqref{eq:coupledsys} have to be substituted with the
functions $k(\e s)$ and $f(\e s)$ satisfying \eqref{eq:ovhovks}.
Precisely, in $N \g_\e$ we define
\begin{equation}\label{eq:tildeLep}
    L^1_\e u = - \D_{\hat{g}_\e} u + V(\e s) u - h(\e s)^{p-1}
      U(k(\e s) y)^{p-1} u - (p-1) h(\e s)^{p-1} U(k(\e s) y)^{p-1} e^{- i
      \frac{f(\e s)}{\e}} \Re(e^{- i \frac{f(\e s)}{\e}} \ov{u})
\end{equation}
(recall the definition of $\hat{g}_\e$ in Subsection \ref{ss:model}:
in particular, working with the coordinates $(s,y)$ integrals will
be computed using the co-area formula \eqref{eq:coarea2}). Before
proving rigorous results, we first discuss heuristically what the
approximate kernel of $L^1_\e$ should look like. Using Fourier
expansions as above (freezing the coefficients at some $\ov{s}$),
the profile of the functions which lie in an approximate kernel of
$L^1_\e$ will be given by the solution of (recall \eqref{eq:systuv})
\begin{equation}\label{eq:systuvss}
  \left\{
    \begin{array}{ll}
      - \D_{y} u + (1 + \a^2) u - p U(y)^{p-1} u + 2
\frac{f'(\ov{s})}{k(\ov{s})} \a v = \tilde{\l} u & \hbox{ in }
    \R^{n-1}; \\
      - \D_{y} v + (1 + \a^2) v - U(y)^{p-1} v + 2
\frac{f'(\ov{s})}{k(\ov{s})} \a u = \tilde{\l} v & \hbox{ in }
    \R^{n-1},
    \end{array}
  \right.
\end{equation}
where $\tilde{\l}$ is close to zero. For $\a$ small (low Fourier
modes), Proposition \ref{p:alphamu} {\bf d} gives the profile $\n_y
U(k(\ov{s}) y)$ or $i U(k(\ov{s}) y)$ (recall also the scaling in
$y$ before \eqref{eq:coupledsys}). The remaining part of the
approximate kernel is the counterpart of that given in Proposition
\ref{p:alphamu} {\bf e}: for variable coefficients it is uniquely
defined a function $\a(\ov{s})$ such that
\begin{equation}\label{eq:defaovs}
    \eta_{\a(\ov{s})} = 0,
\end{equation}
where $\eta_\a$ here stands for the first eigenvalue of
\eqref{eq:systuvss}. We denote by $\left(Z_{\a(\ov{s})}(k(\ov{s})
y), W_{\a(\ov{s})}(k(\ov{s}) y)\right)$ the components of the
relative eigenfunction.

We next consider two bases of eigenfunctions for the weighted
eigenvalue problems (the operators $\mathfrak{J}$ and $T$ are
defined in \eqref{eq:2ndvarfin4}, \eqref{eq:defT} and are
self-adjoint)
\begin{equation}\label{eq:eigen}
  \mathfrak{J} \, \varphi_j(\ov{s}) = h(\ov{s})^\th \l_j \varphi_j(\ov{s});
  \qquad \qquad T \o_j = h(\ov{s})^{-\s} \rho_j \o_j.
\end{equation}
Because of the weights on the right-hand sides, we can choose these
eigenfunctions to be normalized so that $\int_0^L h^\th \varphi_j
\varphi_l = \d_{jl}$ and $\int_0^L h^{-\s} \o_j \o_l = \d_{jl}$:
this choice will be useful in Subsections \ref{sss:pr1},
\ref{sss:pr2}.

These heuristic arguments suggest that the following subspaces
$K_{1,\d}, K_{2,\d}$, where $\d$ is a small positive constant, once
multiplied by $e^{- \frac{i f(\e s)}{\e}}$ consist of approximate
eigenfunctions for $L^1_\e$ with eigenvalues close to zero (this
will be verified below, in the proof of Proposition \ref{p:inv}, see
also Remark \ref{r:proof})
\begin{equation}\label{eq:K1d}
    K_{1,\d} = span \left\{ h(\e s)^{\frac{p+1}{4}} \left( \langle \var_j(\e s), \n_y U(k y)
  \rangle + i \e \langle \var'_j(\e s), y \rangle \frac{f'}{k} U(k y) - \frac{\e^2}{k^2}
  \langle \var''_j(\e s), \mathfrak{V}(k y) \rangle \right) \right\};
\end{equation}
\begin{equation}\label{eq:K2d}
    K_{2,\d} = span \left\{ h(\e s)^{\frac 12} \left( \o_j(\e s) i
   U(k y) + \frac{2 \e f'(\e s)}{k} \o'_j(\e s) \tilde{U}(y)
  - \frac{\e^2}{k^2} \o''_j(\e s) \mathfrak{W}(k y) \right) \right\},
\end{equation}

$j = 0, \dots, \frac{\d}{\e}$, where
$$
  \tilde{U} = \left( \frac{1}{h^{p-1}(p-1)} U(k y) + \frac{1}{2k} \n U(k y) \cdot y \right).
$$
Here $\mathfrak{V} = (\mathfrak{V}_j)_{j=1, \dots, n-1}$ is the
counterpart of $\hat{\mathfrak{V}}$ in \eqref{eq:frakvW}
substituting $\hat{h}$ with $h(\ov{s})$ (the same holds for
$\mathfrak{W}$). The choice of the weights (as powers of $h$) in
\eqref{eq:eigen} and in $K_{1,\d}$, $K_{2,\d}$ are again done for
technical reasons, and will be convenient below, see in particular
Subsections \ref{sss:pr1} and \ref{sss:pr2}.

We also need to construct an approximate kernel with  the profile
$(Z,W)$, see Proposition \ref{p:alphamu} {\bf e}. To this aim we
introduce the functions (recall \eqref{eq:defaovs})
\begin{equation}\label{eq:Q1Q2Q3}
    Q_{1,\a}(\ov{s}) = \int_{\R^{n-1}}Z_{\a(\ov{s})}^2; \qquad Q_{2,\a}(\ov{s}) =
\int_{\R^{n-1}}W_{\a(\ov{s})}^2;
 \qquad Q_{3,\a}(\ov{s}) = \int_{\R^{n-1}} Z_{\a(\ov{s})} W_{\a(\ov{s})},
\end{equation}
 and consider the following eigenvalue
problem (with periodic boundary conditions)
$$
-\e^2 \xi''-k^2 \a^2\xi=\frac{\tilde{\nu}}{1+2f'\frac{Q_{3,\a}}{k\a}}\xi \qquad
\hbox{ on } [0,L].
$$
By the Weyl's asymptotic formula we have that the eigenvalues
$\tilde{\nu}_j$ (counted with multiplicity) have the qualitative
behavior $\tilde{\nu}_j \simeq -1 + \e^2 j^2$ as $j \to + \infty$.
Hence, there is a first index $j_\e$ (of order $\frac 1\e$) for
which $\tilde{\nu}_{j_\e} \geq 0$. Setting $\nu_j =
\tilde{\nu}_{j_\e + j}$ and denoting by $\xi_j$ the eigenfunctions
corresponding  to $\nu_j$, then we have
\begin{equation}\label{eq:asynul}
   -\e^2 \xi_j''- k^2 \a^2 \xi_j =\frac{\nu_j}{1+2f'\frac{Q_{3,\a}}{k\a}} \xi_j
  ; \qquad \qquad \quad \nu_j = \hat{C}_0 \e j + O(\e^2 j^2) + O(\e) \quad  \hbox{ if } |j| \leq
\frac{\d^2}{\e},
\end{equation}
where $\d > 0$ is any given positive (small) constant. Notice that
the family $(\xi_j)_j$ can be chosen normalized in $L^2$ with the
weight $\frac{1}{1+2f'\frac{Q_{3,\a}}{k\a}}$ (this follows from
\eqref{eq:asynul} and the Courant-Fischer formula). Next we  set
\begin{equation}\label{eq:bbj} \b_j=-\frac{1}{k\a}\left(
1-\frac{Q_{1,\a}}{k^2\a^2+2f'k\a Q_{3,\a}}  \nu_j \right)\e\xi_j'.
\end{equation}
By our choices, the functions $\b_j$ and $\xi_j$ satisfy (this
system will be useful in Subsection \ref{ss:contr})
\begin{equation}\label{eq:systbjxj}
  \left\{
    \begin{array}{ll}
      - \e^2 \b''_j - k^2 \a^2 \b_j - 2 f'
    \frac{Q_{3,\a}}{Q_{1,\a}} (\e \xi'_j + k \a \b_j) = \nu_j \b_j+ (O(\nu_j^2)
     + O(\e))\b_j; &  \\ - \e^2\xi''_j - k^2 \a^2 \xi_j + 2 f'
    \frac{Q_{3,\a}}{Q_{2,\a}} (\e \b'_j - k \a \xi_j) = \nu_j \xi_j+(O(\nu_j^2)
    + O(\e))\xi_j &
    \end{array}
  \right. \qquad  \hbox{ for } |j| \leq
\frac{\d^2}{\e}.
\end{equation}
Our next goal is to introduce a family of approximate eigenfunctions
of $L^1_\e$ which have the profile (from now on, we might denote
$(Z_{\a(\ov{s})}, W_{\a(\ov{s})})$, see \eqref{eq:defaovs}, simply
by $(Z_\a, W_\a)$)
\begin{equation}\label{eq:vv33dd}
    v_{3,j}:=(\b_j+q_j) Z_\a + \g_j \frac{\pa Z_\a}{\pa \a} + i
    \xi_j W_\a + i \kappa_j \frac{\pa W_\a}{\pa \a}:
\end{equation}
the functions $\b_j, \xi_j$ are as in
\eqref{eq:asynul}-\eqref{eq:bbj}, while $q_j,\g_j$ and $\kappa_j$
are small corrections to be chosen properly so that $L^1_\e
(e^{-i\frac{f(\e s)}{\e}}v_{3,j}) = \nu_j e^{-i\frac{f(\e
s)}{\e}}v_{3,j}$, up to an error $o(\nu_j) + O(\e)$.

With simple computations, using \eqref{eq:systuva},
\eqref{eq:systuvadiff}, Remark \ref{r:Asmall0} and
\eqref{eq:systbjxj}, one finds that

\begin{eqnarray*}
&&e^{i\frac{f(\e s)}{\e}} \left(L^1_\e (e^{-i\frac{f(\e s)}{\e}}
v_{3,j})-\nu_je^{-i\frac{f(\e s)}{\e}}v_{3,j}\right)\\
&=& Z_\a 2 f' \frac{Q_{3,\a}}{Q_{1,\a}} \left[ 2Q_{1,\a} \g_j k + (\e \xi'_j + k \a \b_j)
\right] + W_\a 2 f' \left( - \g_j k - (\e \xi'_j + \b_j k \a)-k\a q_j \right) \\ & + &
\frac{\pa Z_\a}{\pa \a} (- \e^2 \g''_j - \a^2 k^2 \g_j) + \frac{\pa W_\a}{\pa \a}
2 f' (- \e \kappa'_j - k \a \g_j) \\& + & i Z_\a 2 f' (- \kappa_j k - k \a \xi_j + \e \b'_j+\e q_j')
+ i W_\a 2 f' \frac{Q_{3,\a}}{Q_{2,\a}} (2Q_{2,\a} k \kappa_j + k \a \xi_j - \e \b'_j) \\
& + & i \frac{\pa Z_\a}{\pa \a} 2 f' (- k \a \kappa_j + \e \g'_j) + i  \frac{\pa W_\a}{\pa \a}
(- \e^2 \kappa''_j - \a^2 k^2 \kappa_j)+Z_{\a}(-\e^2q''_j-\a^2 k^2 q_j) \\
& + & (O(\nu_j^2) + O(\e))\xi_j  +(O(\nu_j^2)  + O(\e))\b_j.
\end{eqnarray*}
To make the coefficients of the terms $Z_{\a}$, $W_{\a}$ and
$iW_{\a}$ in the second and fourth lines vanish we choose
$$
\g_j= - \frac{1}{2kQ_{1,\a}}(\e \xi'_j +  k \a\b_j);  \qquad  \kappa_j=
-\frac{1}{2kQ_{2,\a}} (k \a \xi_j - \e \b'_j); \qquad
q_j=-\frac{1}{k\a}(\e \xi'_j + \b_j k \a+k\g_j).
$$
Using \eqref{eq:bbj} we get
\begin{equation}\label{eq:gammajkappaj}
    \g_j=-\e \xi'_j\frac{1}{2k}\frac{1}{k^2\a^2+2f'k\a Q_{3,\a}}  \nu_j; \qquad \qquad
\kappa_j = \frac{1}{2k} \frac{\nu_j}{1 + 2 f' \frac{Q_{3,\a}}{k \a}} \xi_j + O(\nu_j^2) \xi_j:
\end{equation}
these equations and \eqref{eq:asynul} imply the relations $- \e^2
\g''_j - \a^2 k^2 \g_j=O(\nu_j^2)\b_j, - \e \kappa'_j - k \a
\g_j=O(\nu_j^2)\b_j, - k \a \kappa_j + \e \g'_j=O(\nu_j^2)\xi_j$ and
$- \e^2 \kappa''_j - \a^2 k^2 \kappa_j=O(\nu_j^2)\xi_j$. Similarly
one finds
\begin{equation}\label{eq:qqj}
    q_j=\xi_j\frac{1}{2k}\frac{1}{k\a+2f'Q_{3,\a}}\nu_j.
\end{equation}
This also yields $- \g_j k - (\e \xi'_j + \b_j k \a)-k\a q_j =
O(\nu_j^2)\b_j$, $- \kappa_j k - k \a \xi_j + \e \b'_j+\e q_j' =
O(\nu_j^2)\xi_j$ and  $-\e^2q''_j-\a^2 k^2 q_j=O(\nu_j^2)\b_j$, so
we obtain
\begin{equation}\label{eq:apped}
    L^1_\e
(e^{-i\frac{f(\e s)}{\e}}v_{3,j}) = \nu_j e^{-i\frac{f(\e
s)}{\e}}v_{3,j} + (O(\nu_j^2) + O(\e))\xi_j  +(O(\nu_j^2)
 + O(\e))\b_j \qquad  \hbox{ for } |j| \leq \frac{\d^2}{\e},
\end{equation}
which was our claim. We next define
\begin{equation}\label{eq:K3d}
    K_{3,\d} = span \left\{ (\b_j+q_j) Z_\a + \g_j
\frac{\pa Z_\a}{\pa \a} + i \xi_j W_j+ i \kappa_j \frac{\pa
W_\a}{\pa \a} \; : \; j = - \frac{\d^2}{\e}, \dots, \frac{\d^2}{\e} \right\}.
\end{equation}
In the $K_{l,\d}$'s we added some corrections to the approximate
eigenfunctions which take into account the variation of the profile
with the frequency, see the derivation of \eqref{eq:systuv} and
Remark \ref{r:profeigenf}. In $K_{1,\d}$ and $K_{2,\d}$ the
corrections are up to the second order (in $\e j$), while in
$K_{3,\d}$ only up to the first: the reason is that the
corresponding eigenvalues have a {\em quadratic} dependence in $\e
j$ for $K_{1,\d}$ and $K_{2,\d}$  (they correspond to $\eta_\a$ in
Proposition \ref{p:alphamu}), and an {\em affine} dependence in $\e
j$ for $K_{3,\d}$ (corresponding to $\mu_\a$ in Proposition
\ref{p:alphamu}). Since the former dependence is more delicate in
the indices, we need a more accurate expansion of the
eigenfunctions. We finally set
\begin{equation}\label{eq:Kd}
    K_\d = span \left\{ K_{1,\d}, K_{2,\d}, K_{3,\d} \right\}.
\end{equation}

\subsection{Invertibility of $L_\e$ in the orthogonal complement of
$K_\d$}\label{ss:invcompl}

Since $K_\d$ (multiplied by $e^{- i \frac{f(\e s)}{\e}}$) is a good
candidate for the span of the eigenfunctions of $L^1_\e$ with small
eigenvalues, it seems plausible to invert $L^1_\e$ on the orthogonal
complement to $e^{- i \frac{f(\e s)}{\e}} K_\d$: this is the content
of the next result. We recall the definition of the constant
$\mathcal{A}$ in the introduction.

\begin{pro}\label{p:inv} There exists $\mathcal{A}_0$ sufficiently
small such that for any $\mathcal{A} \in [0,\mathcal{A}_0]$ the
following property holds. For $\d > 0$ small enough there exists $C
> 0$ (independent of $\d$) such that if
\begin{equation}\label{eq:ovHe}
    \Re \int_{N \g_\e} e^{- i \frac{f(\e s)}{\e}} v \, \ov{\phi}
  \, dV_{\hat{g}_\e} = 0 \qquad \quad \hbox{ for all } v \in K_\d,
\end{equation}
one has $\|\Pi_\e L^1_\e(\phi)\|_{L^2(N \g_\e)} \geq C^{-1} \d^2
\|\phi\|_{L^2(N \g_\e)}$. Here $\Pi_\e$ denotes the
 projection in $L^2(N \g_\e)$ onto the orthogonal complement of the
 set $\left\{e^{- i \frac{f(\e s)}{\e}} v \; : \; v \in K_\d \right\}$.
\end{pro}

\noindent Before starting with the proof, which relies on a
localization argument and the spectral analysis of Proposition
\ref{p:alphamu}, we introduce some notation and a preliminary Lemma.
We fix $\hat{s} \in [0,L]$ and we denote by $\hat{f}$, $\hat{h}$,
$\hat{k}$, $\hat{\a}$ the values of $f'(\hat{s}), h(\hat{s}),
k(\hat{s}), \a(\hat{s})$ respectively, so that the counterpart of
\eqref{eq:ovhovks} holds true. For a large constant $\tilde{C}_0$ to
be fixed later, we also define
$$
  \hat{K}_{1,\d} = span \left\{ \langle \hat{\var}_j(\e s), \n_y U(\hat{k} y)
  \rangle + i \e \langle \hat{\var}'_j(\e s), y \rangle \frac{\hat{f}}{\hat{k}}
  U(\hat{k} y) - \frac{\e^2}{k^2} \langle \var''_j(\e s), \hat{\mathfrak{V}}(\hat{k} y)
  \rangle  \right\};
$$
$$
  \hat{K}_{2,\d} = span
\left\{ \left( \hat{\o}_j(\e s) i
   U(\hat{k} y) + \frac{2 \e f'(\e s)}{\hat{k}} \hat{\o}'_j(\e s) \hat{\tilde{U}}(y)
  - \frac{\e^2}{\hat{k}^2} \tilde{\o}''_j(\e s) \hat{\mathfrak{W}}(\hat{k} y) \right) \right\},
$$
$j = 0, \dots, \frac{\d}{\tilde{C}_0 \e}$, and
$$
  \hat{K}_{3,1,\d} = span \left\{  \cos (\hat{\a}_j s) Z_{\hat{\a}_j}(\hat{k} y)
  - i \sin (\hat{\a}_j s) W_{\hat{\a}_j}(\hat{k} y) \; : \; j = -
  \frac{\d^2}{\tilde{C}_0 \e}, \dots, \frac{\d^2}{\tilde{C}_0 \e} \right\};
$$
$$
  \hat{K}_{3,2,\d} = span \left\{  \sin (\hat{\a}_j s) Z_{\hat{\a}_j}(\hat{k} y)
  + i \cos (\hat{\a}_j s) W_{\hat{\a}_j}(\hat{k} y) \; : \; j = -
  \frac{\d^2}{\tilde{C}_0 \e}, \dots, \frac{\d^2}{\tilde{C}_0 \e} \right\},
$$
where
\begin{equation}\label{eq:bbbello}
    \hat{\tilde{U}} = \left( \frac{1}{\hat{h}^{p-1}(p-1)} U(\hat{k} y) +
  \frac{1}{2\hat{k}} \n U(\hat{k} y) \cdot y \right); \qquad \qquad \hat{\a}_j =
  \left[ \frac{\hat{\a} \hat{k} L}{2 \pi \e} + j \right] \frac{2 \pi \e}{L},
\end{equation}
(again, the latter square bracket stands for the integer part, and
this choice makes the functions $L/\e$-periodic). In the above
formulas $(\hat{\var}_j)_j$ are the eigenfunctions of the normal
Laplacian with the flat metric on $\g$, and  $\hat{\o}_j$ those of
$\pa^2_{\ov{s} \ov{s}}$ on $[0,L]$: the symbols $Z_{\hat{\a}_j},
W_{\hat{\a}_j}$ stand for the components of the eigenfunctions of
\eqref{eq:systuv} corresponding to $\eta_{\hat{\a}_j}$. In analogy
with \eqref{eq:Kd} we also define
\begin{equation}\label{eq:hatKd}
    \hat{K}_\d = span \left\{ \hat{K}_{1,\d}, \hat{K}_{2,\d},
    \hat{K}_{3,1,\d}, \hat{K}_{3,2,\d} \right\}.
\end{equation}
Given a small constant $\eta > 0$ to be chosen later (of order
$\sqrt{\e}$), we consider then a smooth cutoff function $\chi_\eta$
(depending only on  $s$) with support near $\frac{\hat{s}}{\e}$ and
with length of order $\frac{\eta}{\e}$. For example, one can take
$\chi_\eta(s) = \chi( \tfrac{\e}{\eta} (s - \hat{s}/\e) )$ for a
fixed compactly supported cutoff $\chi$ which is $1$ in a
neighborhood of $0$. The next result uses Fourier cancelation, and
is related to Lemma 2.7 in \cite{mal}.

\begin{lem}\label{l:cancfour}
There exists $\tilde{C}_0$ sufficiently large (depending only on
$V$, $L$ and $\mathcal{A}_0$) with the following property. For any
given integer number $m$ there exists $C_m > 0$ depending on $m$ and
$\chi_\eta$ such that for $|j| \leq \frac{\d^2}{\tilde{C}_0 \e}$ and
for $|l| \geq \frac{\d^2}{\e}$ one has
$$
\left| \int \chi_\eta(s) \xi_l(s) \cos(\hat{\a}_j s) ds
\right| + \left| \int \chi_\eta(s) \xi_l(s) \sin(\hat{\a}_j s) ds
\right| \leq \frac 1 \e \frac{C_m}{|\nu_l|^m}
    \left[ \eta (1+ |\nu_l|) + \frac{\e}{\eta} \right]^m.
$$
\end{lem}

\begin{pf} We clearly have that $(\cos(\hat{\a}_j s))'' = - \hat{\a}_j^2
  \cos(\hat{\a}_j s)$: therefore, integrating by parts, after some manipulation we obtain
that
\begin{eqnarray}\label{eq:trucco} \nonumber
  \int \chi_\eta \xi_l(s)
    \cos(\hat{\a}_j s) ds & = & \frac{1}{\hat{\a}_j^2  -
   \hat{k}^2 \hat{\a}^2 + \frac{\nu_l}{1 + 2 \hat{f} \frac{Q_{3,\hat{\a}}}{\hat{k} \hat{\a}}}}
    \times \\ & \times & \left\{ \int \chi_\eta \xi_l(s)
    \cos(\hat{\a}_j s) \left( k^2 \a^2 - \hat{k}^2 \hat{\a}^2 +
 \frac{\nu_l}{1 + 2 \hat{f} \frac{Q_{3,\hat{\a}}}{\hat{k} \hat{\a}}}
- \frac{\nu_l}{1 + 2 f' \frac{Q_{3,\a}}{k \a}}\right) \right. \\
   & - & \left. \int \cos(\hat{\a}_j s) \left[ \chi''_\eta
   \xi_l(s) + 2 \chi'_\eta \xi'_l(s) \right] \right\}. \nonumber
\end{eqnarray}
By \eqref{eq:bbbello} the numbers $\hat{\a}_j$ satisfy $\hat{\a}_j
\simeq \hat{k} \hat{\a} + \frac{2 \pi \e j}{L}$ for $|j| \leq
\frac{\d^2}{\tilde{C}_0 \e}$, while by \eqref{eq:asynul} $|\nu_l|
\geq \frac{1}{2} \hat{C}_0 \d^2$ for $|l| \geq \frac{\d^2}{\e}$.
Notice also that $1 + 2 f' \frac{Q_{3,\a}}{k \a}$ is uniformly
bounded above and below by positive constants when $\mathcal{A}_0$
tends to zero (see for example the comments in Remark
\ref{r:Asmall2}). By these facts and the properties of $\chi_\eta$
we find
$$
\left| \int \chi_\eta \xi_l(s) \cos(\hat{\a}_j s) ds \right| \leq \frac 1 \e
   \frac{C}{|\nu_l|}  \left[ \eta (1+ |\nu_l|) + \frac{\e}{\eta}  \right].
$$
which yields the statement for $m = 1$ (similar computations can be
performed to deal with the $\sin$ function). The factor $\frac 1 \e$
inside the brackets arises from the fact that we are integrating
over the interval $[0,L/\e]$, and by the normalization of $\xi_j$
(see the comments before \eqref{eq:asynul}). To obtain it for
general $m$, it is sufficient to iterate the procedure for
\eqref{eq:trucco} $m$ times and integrate by parts.
\end{pf}

\

\begin{pfn} {\sc of Proposition \ref{p:inv}}
The proof mainly relies  on a localization argument and Lemma
\ref{l:cancfour}. If $\eta = \sqrt{\e}$ and $\chi_\eta$ is as in
Lemma \ref{l:cancfour}, we show next that the function $\chi_\eta
\phi$ is {\em almost} orthogonal to $e^{- i \hat{f} s} \hat{K}_\d$
if $\e$ is sufficiently small. We consider for example a function
$\hat{v} \in \hat{K}_{3,1,\d}$ of the form
$$
  \hat{v} = \sum_{j = -\frac{\d^2}{\tilde{C}_0 \e}}^{\frac{\d^2}{\tilde{C}_0 \e}} \hat{b}_j
  \left[ \cos (\hat{\a}_j s) Z_{\hat{\a}_j}(\hat{k} y)
  - i \sin (\hat{\a}_j s) W_{\hat{\a}_j}(\hat{k} y) \right],
$$
for some arbitrary coefficients $(\hat{b}_j)_j$, and we also set
$$
  \tilde{v} = \sum_{j = -\frac{\d^2}{\tilde{C}_0 \e}}^{\frac{\d^2}{\tilde{C}_0 \e}} \hat{b}_j
  \left[ \cos (\hat{\a}_j s) Z_{\a(\e s)}(k(\ov{s}) y)
  - i \sin (\hat{\a}_j s) W_{\a(\e s)}(k(\ov{s}) y) \right].
$$
We are going to evaluate the real part of the integral $\int_{N
\g_\e} e^{- i \hat{f} s} \hat{v} \chi_\eta \ov{\phi}$: first of all,
since $|k(\ov{s}) - \hat{k}| \leq C \eta$ and $|\hat{\a}_j -
\a(\ov{s})| \leq C (\eta + \d^2)$ on the support of $\chi_\eta$ we
have that
$$
\|Z_{\a(\e s)}(k(\ov{s}) y) - Z_{\hat{\a}_j}(\hat{k} y)\|_{L^2(\R^{n-1})}
= O(\eta + \d^2) \qquad \qquad \hbox{ in } supp(\chi_\eta),
$$
so, clearly
\begin{equation}\label{eq:dec0}
   \Re \int_{N \g_\e} e^{- i \hat{f} s} \hat{v} \chi_\eta \ov{\phi} = \Re \int_{N \g_\e}
   e^{- i \hat{f} s} \tilde{v} \chi_\eta \ov{\phi} + O(\eta + \d^2) \|\chi_\eta \phi\|_{L^2(N \g_\e)}
   \|\hat{v}\|_{L^2(N \g_\e)}.
\end{equation}
We next write $\hat{w}(s) = \chi_\eta(s) \sum_{j =
-\frac{\d^2}{\tilde{C}_0 \e}}^{\frac{\d^2}{\tilde{C}_0 \e}}
\hat{b}_j \sin(\hat{\a}_j s)$, and notice that
\begin{eqnarray}\label{eq:dec1}
  & & \hat{w}'(s) = \chi'_\eta(s) \sum_{j = -\frac{\d^2}{\tilde{C}_0 \e}}^{\frac{\d^2}{\tilde{C}_0 \e}}
  \hat{b}_j \sin(\hat{\a}_j s) + \chi_\eta(s) \sum_{j = -\frac{\d^2}{\tilde{C}_0 \e}}^{\frac{\d^2}{\tilde{C}_0 \e}}
  \hat{b}_j \hat{\a}_j \cos(\hat{\a}_j s) \\
   &=& k \a \chi_\eta(s) \sum_{j = -\frac{\d^2}{\tilde{C}_0 \e}}^{\frac{\d^2}{\tilde{C}_0 \e}}
  \hat{b}_j  \cos(\hat{\a}_j s) + \chi'_\eta(s) \sum_{j = -\frac{\d^2}{\tilde{C}_0 \e}}^{\frac{\d^2}{\tilde{C}_0 \e}}
  \hat{b}_j \sin(\hat{\a}_j s) - \chi_\eta(s) \sum_{j = -\frac{\d^2}{\tilde{C}_0 \e}}^{\frac{\d^2}{\tilde{C}_0 \e}}
  \hat{b}_j (k \a - \hat{\a}_j) \sin(\hat{\a}_j s). \nonumber
\end{eqnarray}
Using this formula and the same argument as for \eqref{eq:dec0}  we
get (recall our notation before \eqref{eq:vv33dd})
\begin{equation}\label{eq:dec2}
   \chi_\eta \tilde{v} = \frac{1}{k \a} \hat{w}'(s) Z_\a - i \hat{w}(s) W_\a
  + O(\eta + \d^2) \|\chi_\eta \hat{v}\|_{L^2(N \g_\e)}.
\end{equation}
Notice that, by the explicit form of $\hat{w}$ and $\hat{\a}_j$, for
any integer $m$ one has $\|\hat{w}\|_{H^m([0,L/\e])}^2 \leq C_m
\|\hat{w}\|_{L^2([0,L/\e])}^2$: therefore, if we write $\hat{w}$
with respect to the basis $\xi_l$ as (notice the shift of index
before \eqref{eq:asynul})
\begin{equation}\label{eq:hatws}
    \hat{w}(s) = \sum_{l=-j_\e}^\infty \check{b}_l \xi_l(\e s),
\end{equation}
we also find that
\begin{equation}\label{eq:powers}
    \sum_{l=-j_\e}^\infty (1 + |\nu_l|)^m \check{b}_l^2 \leq C_m
\|\hat{w}\|_{H^m([0,L/\e])}^2 \leq C_m
\|\hat{w}\|_{L^2([0,L/\e])}^2 \leq C_m \|\hat{v}\|_{L^2(N \g_\e)}^2.
\end{equation}
Differentiating \eqref{eq:hatws} with respect to $s$ and using the
definition of $\xi_j$ together with \eqref{eq:bbj} we find that
\begin{eqnarray*}
  \hat{w}'(s) = \sum_{l=-j_\e}^\infty \check{b}_l \e \xi'_l(\e s) =
  \sum_{l=-j_\e}^\infty \check{b}_l (-k \a + O(\nu_l)) \b_l(\e s).
\end{eqnarray*}
The last formula and \eqref{eq:dec2} imply
\begin{eqnarray}\label{eq:dec3} \nonumber
  \chi_\eta \tilde{v} & = & - \sum_{l=-j_\e}^\infty \check{b}_l (\b_l Z_\a
 + i \xi_l W_\a) + \sum_{l=-j_\e}^\infty \check{b}_l O(\nu_l) \b_l Z_\a
+ O(\eta + \d^2) \|\chi_\eta \hat{v}\|_{L^2(N \g_\e)}
 \\ & = & - \sum_{l=-j_\e}^\infty \check{b}_l v_{3,l}
+ \sum_{l=-j_\e}^\infty \check{b}_l (v_{3,l} - \b_l Z_\a - i \xi_l W_\a)
+ \sum_{l=-j_\e}^\infty \check{b}_l O(\nu_l) \b_l Z_\a
+ O(\eta + \d^2) \|\chi_\eta \hat{v}\|_{L^2(N \g_\e)}.
\end{eqnarray}
In the support of $\chi_\eta$ there exists $\hat{\th} \in \R$ such
that $\frac{f(\e s)}{\e} = \hat{f} s + \hat{\th} + O(\eta)$, which
yields
$$
  \int_{N \g_\e} \tilde{v} e^{- i \hat{f} s} \chi_\eta \ov{\phi}  =
\int_{N \g_\e} \tilde{v} e^{- i \frac{f(\e s)}{\e}} \chi_\eta \ov{\phi}
+ O(\eta)  \|\phi\|_{L^2(supp(\chi_\eta))} \|\hat{v}\|_{L^2(N \g_\e)}.
$$
Now, recalling that $\eta = \sqrt{\e}$ and that we have
orthogonality between $\phi$ and $e^{- i \frac{f(\e s)}{\e}} K_\d$,
from the last two formulas we obtain that
\begin{equation}\label{eq:miaomiao}
  \int_{N \g_\e} \tilde{v} e^{- i \hat{f} s} \chi_\eta \ov{\phi} = A_1 + A_2 + A_3 +
 O(\eta + \d^2)  \|\phi\|_{L^2(supp(\chi_\eta))} \|\hat{v}\|_{L^2(N \g_\e)};
\end{equation}
$$
  A_1 = - \int_{supp(\chi_\eta)} e^{- i \hat{f} s} \ov{\phi}
  \sum_{|l| \geq \frac{\d^2}{\e}} \check{b}_l v_{3,l}; \qquad
  A_2 = \int_{supp(\chi_\eta)} e^{- i \hat{f} s} \ov{\phi}
   \sum_{l=-j_\e}^\infty \check{b}_l (v_{3,l} - \b_l Z_\a - i \xi_l W_\a);
$$
$$
  A_3 = \int_{supp(\chi_\eta)} e^{- i \hat{f} s} \ov{\phi}
 \sum_{l=-j_\e}^\infty \check{b}_l O(\nu_l) \b_l Z_\a.
$$
To estimate these terms we notice first that, by the normalization
of $\xi_j$ before \eqref{eq:bbj}, the coefficients $\check{b}_l$ in
\eqref{eq:hatws} can be computed as
$$
  \check{b}_l = \e \int_{0}^{L/\e} \frac{\hat{w}(s)}{1+2f'\frac{Q_{3,\a}}{k\a}}
\xi_l(\e s) ds = \e \sum_{j =
-\frac{\d^2}{\tilde{C}_0 \e}}^{\frac{\d^2}{\tilde{C}_0 \e}}
\hat{b}_j  \int_{0}^{L/\e} \chi_\eta(s)
\frac{\sin(\hat{\a}_j s)}{1+2f'\frac{Q_{3,\a}}{k\a}} \xi_l(\e s) ds.
$$
Using this formula, Lemma \ref{l:cancfour} and the H\"older
inequality we find that for any integer $m$
\begin{equation}\label{eq:intm}
    \check{b}_l^2 \leq C_m \e^{2m+2} \bigg( \sum_{j =
-\frac{\d^2}{\tilde{C}_0 \e}}^{\frac{\d^2}{\tilde{C}_0 \e}}
|\hat{b}_j|  \bigg)^2 \leq C_m \e^{2m+1} \sum_{j =
-\frac{\d^2}{\tilde{C}_0 \e}}^{\frac{\d^2}{\tilde{C}_0 \e}}
\hat{b}_j^2 \leq C_m \e^{2m+2} \|\hat{v}\|_{L^2(N \g_\e)}^2; \qquad
\quad |l| \geq \frac{\d^2}{\e}.
\end{equation}
From the explicit expression for the functions $v_{3,\d}$ the above
term $A_1$ can be estimated as
$$
  A_1 \leq C \big( \frac 1 \e \sum_{|l| \geq \frac{\d^2}{\e}} (1 + |\nu_l|^2)^2
 \check{b}_l^2 \big)^{\frac 12}  \|\phi\|_{L^2(supp(\chi_\eta))}.
$$
As before, the factor $\frac 1 \e$ inside the brackets arises from
the fact that we are integrating over  $[0,L/\e]$.

Using the fact that $C^{-1} |\e l| \leq |\nu_l| \leq C (|\e l| +
\e^2 l^2)$ for $|l| \geq \frac{\d^2}{\e}$ (which follows from the
Weyl's asymptotic formula for the eigenvalue problem in
\eqref{eq:asynul}), \eqref{eq:powers} and \eqref{eq:intm}, one finds
that for any large integer $m$ and any $d \in (1, m/8)$
$$
  \sum_{\frac{\d^2}{\e} \leq |l| \leq \e^{-d}} (1 + |\nu_l|^2)^2 \check{b}_l^2 \leq
  C_m \e^{11 + 2m - 8 d} \|\hat{v}\|_{L^2(N \g_\e)}^2; \quad
  \sum_{|l| \geq \e^{-d}} (1 + |\nu_l|^2)^2 \check{b}_l^2 \leq C_m \e^{(d-1)(m-4)}
  \|\hat{v}\|_{L^2(N \g_\e)}^2.
$$
By the arbitrarity of $m$ it follows that for any $m' \in \N$
$$
  |A_1| \leq C_{m'} \e^{m'} \|\hat{v}\|_{L^2(N \g_\e)} \|\phi\|_{L^2(supp(\chi_\eta))}.
$$
Dividing the set of indices $l$ into $\{ |l| \leq \frac{\d^2}{\e}
\}$ and $\{ |l| \geq \frac{\d^2}{\e} \}$ and using similar arguments
(taking also into account \eqref{eq:gammajkappaj} and
\eqref{eq:qqj}) we get
$$
  |A_2| + |A_3| \leq C  \|\phi\|_{L^2(supp(\chi_\eta))} \bigg( \frac 1 \e
  \sum_{l=-j_\e}^\infty (\nu_l^2 + \nu_l^4) \check{b}_l^2 \bigg)^{\frac 12} \leq
  C \d^2 \|\phi\|_{L^2(supp(\chi_\eta))} \|\hat{v}\|_{L^2(N \g_\e)}.
$$
Therefore, using \eqref{eq:dec0} and \eqref{eq:miaomiao} one finds
$$
\Re \int_{N \g_\e} e^{- i \hat{f} s} \hat{v} \chi_\eta \ov{\phi} =
 O(\eta + \d^2) \|\phi\|_{L^2(supp(\chi_\eta))} \|\hat{v}\|_{L^2(N \g_\e)};
\qquad \quad \hat{v} \in \hat{K}_{3,1,\d}.
$$
Similar estimates hold for $\hat{v} \in span \{ \hat{K}_{1,\d},
\hat{K}_{2,\d}, \hat{K}_{3,2,\d} \}$, so we obtain
\begin{equation}\label{eq:projchiphionhatK}
  \int_{N \g_\e} e^{- i \hat{f} s} \hat{v} \chi_\eta \ov{\phi} =
  O(\d^2 + \eta) \|\phi\|_{L^2(supp(\chi_\eta))}
  \|\hat{v}\|_{L^2(N \g_\e)}
  \quad \hbox{ for every } \hat{v} \in \hat{K}_\d.
\end{equation}
Next we let $\hat{L}_\e$ denote the operator in \eqref{eq:hatLe}
with coefficients freezed at $\hat{s}$. Since $e^{- i \hat{f} s}
\hat{K}_\d$ consist of all the eigenfunctions of $\hat{L}_\e$ (up to
an error $o(\d^2)$) with eigenvalues smaller in absolute value than
$\d^2$, see Proposition \ref{p:alphamu} and Remark
\ref{r:profeigenf}, from \eqref{eq:projchiphionhatK} we then deduce
\begin{equation}\label{eq:bau}
    \|\hat{L}_\e (e^{-i \hat{f} s} \chi_\eta \phi)\|_{L^2(N \g_\e)}
    \geq \frac{\d^2}{C} \|\chi_\eta \phi\|_{L^2(N \g_\e)} + O(\d^4 +
    \d^2 \eta) \|\phi\|_{L^2(supp(\chi_\eta))}
\end{equation}
for some fixed constant $C$ independent of $\d$.

It is now possible to choose the cutoff function $\chi$ (see the
comments before Lemma \ref{l:cancfour}) so that it is even,
compactly supported in $[-2,2]$, $\chi \equiv 1$ in $[-1,1]$, and so
that $\chi(2-t) + \chi(t) \equiv 1$ for $t \in [1,2]$. With this
choice, we can find a partition of unity $(\chi_{\eta,j})_j$ of
$[0,L/\e]$ consisting of translates of $\chi_\eta$ (plus a
negligible scaling), with $j$ running between $1$ and a number of
order $\frac{L}{\sqrt{\e}}$. For each index $j$ we choose a point
$\hat{s}_j$ in the support of $\chi_{\eta,j}$ and we denote by
$\hat{L}_j$ the  operator corresponding to \eqref{eq:hatLe} with
coefficients freezed at $\hat{s}_j$. Then, using \eqref{eq:bau},
with easy computations one finds
\begin{eqnarray*}
  \|L^1_\e \phi\|^2_{L^2(N \g_\e)} &=&  \big\|
L^1_\e \sum_j \chi_{\eta,j}
  \phi \big\|^2_{L^2(N \g_\e)} = \big\| \sum_j \hat{L}_j
  (\chi_{\eta,j} \phi) \big\|^2_{L^2(N \g_\e)} + O(\sqrt{\e}) \|\phi\|^2_{L^2(N \g_\e)}
  \\  & \geq & \frac{\d^4}{C} \sum_j \big\| \chi_{\eta,j}
  \phi \big\|^2_{L^2(N \g_\e)} + O(\sqrt{\e}) \|\phi\|^2_{L^2(N \g_\e)}
  = \frac{\d^4}{C} \| \phi \|^2_{L^2(N \g_\e)}
  + O(\sqrt{\e}) \|\phi\|^2_{L^2(N \g_\e)}
\end{eqnarray*}
for some $C$ independent of $\d$. To complete the proof we need to
bound from below the norm of $\Pi_\e L^1_\e \phi$, showing that
$\|\Pi_\e L^1_\e \phi\|^2_{L^2(N \g_\e)} \geq \frac{\d^2}{C}
\|\phi\|^2_{L^2(N \g_\e)}$. To see this, by the last formula it is
sufficient to have
\begin{equation}\label{eq:show}
    (L^1_\e  \phi, e^{- i \frac{f(\e s)}{\e}}
    v)_{L^2(N \g_\e)} = o(\d^2) \|\phi\|_{L^2(N \g_\e)}
\|v\|_{L^2(N \g_\e)} \qquad \hbox{ for any } v \in K_\d.
\end{equation}
We prove this claim for $v \in K_{1,\d}$ only: for the other
$K_{j,\d}$'s the arguments are similar, see Remark \ref{r:proof}
below for more details. Setting $v = v_r + i v_i$ one finds (see
\eqref{eq:LrLi})
$$
  L^1_\e (e^{-\frac{i f(\e s)}{\e}} v) = e^{-\frac{i f(\e s)}{\e}} (\mathcal{L}_r v_r
  + i \mathcal{L}_i v_i ) - e^{-\frac{i f(\e s)}{\e}} \left( \frac{\pa^2 v_r}{\pa s^2}
  + i \frac{\pa^2 v_i}{\pa s^2} \right) + 2 i f' e^{-\frac{i f(\e s)}{\e}}
  \left( \frac{\pa v_r}{\pa s} + i \frac{\pa v_i}{\pa s} \right).
$$
When differentiating $v$ with respect to $s$, we either hit the
functions $\varphi_j$'s (and their derivatives) or other functions
like $k(\e s)$ or $f'(\e s)$ (see the definition of $K_{1,\d}$
above). The latter ones have a {\em slow} dependence in $s$ and
therefore these terms can be collected within an error of the form
$O(\e) \|v\|_{L^2(N \g_\e)}$.

However, by our choices of the second and the third parts of the
elements in $K_{1,\d}$ (see Remark \ref{r:profeigenf}, in particular
formulas \eqref{eq:Lri}, \eqref{eq:lrj} and \eqref{eq:lrr}), terms
containing zero-th or first order derivatives of $\varphi_j$ will
have coefficients bounded  by $\e$, while the only term containing
second derivatives of $\var_j$ will be a linear combination (in $j$)
of the expressions
$$
  - \e^2 h(\e s)^{\frac{p+1}{4}} \left( 1 - \frac{2 \mathcal{A}^2 \th}{p-1} h(\ov{s})^{\s-\th}
   \right) \langle \var_j''(\e s), \n_y U(k(\ov{s}) y) \rangle, \qquad \qquad
  j = 0, \dots, \frac{\d}{\e}.
$$
The remaining terms will contain third and the fourth derivatives of
$\var_j$ only (multiplied respectively by $\e^3$ and $\e^4$).
Therefore, if we set
$$
  v_{1,j} = h(\e s)^{\frac{p+1}{4}} \left( \langle \var_j(\e s), \n_y U(k y) \rangle
  + i \e \langle \var'_j(\e s), y \rangle \frac{f'}{k} U(k y) - \frac{\e^2}{k^2}
  \langle \var''_j(\e s), \mathfrak{V}(k y) \rangle \right),
$$
by the above comments and the fact that $\mathfrak{J} \var_j = h(\e
s)^\th \l_j \var_j$ (see \eqref{eq:eigen}) we have
\begin{equation}\label{eq:expa}
    L^1_\e \big(e^{- i \frac{f(\e s)}{\e}} \sum_{j=0}^{\frac \d \e} a_j v_{1,j} \big)
= e^{- i \frac{f(\e s)}{\e}} \sum_{j=0}^{\frac \d \e} \l_j a_j v_{1,j} + R(v);
 \qquad \quad v = \sum_{j=0}^{\frac \d \e} a_j v_{1,j},
\end{equation}
where $R(v)$ contains terms of order $\e$ or linear combinations of
third and fourth derivatives of $\varphi_j(\e s)$, so using Fourier
analysis one can derive the estimate
\begin{equation}\label{eq:estRv}
    \|R(v)\|_{L^2(N \g_\e)} \leq C \big( \frac 1 \e \sum_{j=0}^{\frac \d \e} a_j^2
   (\e + \e^3 j^3)^2 \big)^{\frac 12} \leq C (\e + \d^3) \|v\|_{L^2(N \g_\e)},
\end{equation}
for some constant $C > 0$. Therefore, using \eqref{eq:expa} and
\eqref{eq:estRv}  we obtain
$$
(L^1_\e
\phi, e^{- i \frac{f(\e s)}{\e}} v)_{L^2(N \g_\e)} = O(\e +
\d^3) \|\phi\|_{L^2(N \g_\e)} \|v\|_{L^2(N \g_\e)},
$$
which yields \eqref{eq:show} and concludes the proof.
\end{pfn}

\begin{rem}\label{r:proof}
The last step in the proof of Proposition \ref{p:inv} is nearly
identical for $v \in K_{2,\d}$ except that, still by the
computations in Remark \ref{r:profeigenf}, in the counterpart of
\eqref{eq:expa} we will obtain $\rho_j$ instead of $\l_j$ (see
\eqref{eq:eigen}). When considering  $K_{3,\d}$, setting $v =
\sum_{j=- \frac{\d^2}{\e}}^{\frac{\d^2}{\e}}
  a_j v_{3,j}$  (see \eqref{eq:apped}) one finds
\begin{equation}\label{eq:expa22}
    L^1_\e \big(e^{- i \frac{f(\e s)}{\e}} v \big) = e^{- i \frac{f(\e s)}{\e}} \sum_{j=-
\frac{\d^2}{\e}}^{\frac{\d^2}{\e}} \nu_j a_j v_{3,j} + \tilde{R}(v);
\end{equation}
\begin{equation*}
    \|\tilde{R}(v)\|_{L^2} \leq C \bigg( \frac 1 \e \sum_{j=-
   \frac{\d^2}{\e}}^{\frac{\d^2}{\e}} a_j^2 (\e + \e^2 j^2)^2
   \bigg)^{\frac 12} \leq C (\e + \d^4) \|v\|_{L^2(N \g_\e)}.
\end{equation*}
\end{rem}

\subsection{Invertibility of $L_\e$ in weighed
spaces}\label{ss:invweight}

Our goal is to show that the linearized operator $L_\e$ (see
\eqref{eq:defLe}) at approximate solutions is invertible on spaces
of functions satisfying suitable constraints. We begin with some
preliminary notation and lemmas: we first collect a decay properties
of Green's kernels in Euclidean space. Let us consider the equation
\begin{equation}\label{eq:equfinrn-1}
    - \D u + u = f \qquad \quad \hbox{ in } \R^{n-1},
\end{equation}
where $f$ decays to zero at infinity. The solution of the above
equation can be represented as
$$
  u(x) = \int_{\R^{n-1}} G_0(|x-y|) f(y) dy,
$$
where $G_0 : \R_+ \to \R_+$ is a function singular at $0$ which
decays exponentially to zero at infinity. Using the notation of
Subsection \ref{subsec:loc} and standard elliptic regularity theory,
one can prove the following result (the choice $\a \geq \frac 12$
for the H\"older exponent is technical, and is used in the proof of
Lemma \ref{l:solvortflat}).

\begin{lem}\label{l:decgreenflat}
Let $\ov{Q} > 0$, $\a \geq \frac 12$, let $0 < \t < 1$, $0
<\varsigma < 1$ and let $f \in C^\a_\tau$. Then equation
\eqref{eq:equfinrn-1} has a (unique) solution $u$ of class
$C^{2,\a}_{\varsigma \tau}$ which vanishes on $\partial
B_{\ov{Q}}(0)$. Moreover, there exist $\varsigma_0 > 0$ sufficiently
close to $1$ and $\check{C}_0$ sufficiently large (depending only on
$n$, $\a$, $\min\{\ov{Q},1\}$ and $\varsigma$) such that for
$\varsigma_0 \leq \varsigma < 1$
$$
  \|u\|_{C^{2,\a}_{\varsigma \tau}} \leq \check{C}_0 \|f\|_{C^\a_\varsigma}.
$$
\end{lem}

\

\noindent Let now $\t, \varsigma \in (0,1)$ (to be fixed later). For
any integer $m$ we let $\ov{C}^{m,\t}_{\varsigma,}$ denote the
weighted H\"older space
\begin{equation}\label{eq:cmas}
  \ov{C}^{m,\t}_{\varsigma} = \left\{ u : \R^{n-1} \to \C \; : \; \sup_{y \in \R^{n-1}}
  e^{\varsigma  |y|} \|u\|_{C^{m,\t}(B_1(y))} < + \infty \right\}.
\end{equation}
We also consider the following set of  functions $L/\e$-periodic in
$s$
\begin{equation}\label{eq:l2cmasfl}
    \ov{L}^2(\ov{C}^{m,\t}_{\varsigma}) = \left\{ u : [0,L/\e] \times \R^{n-1} \to \C
    \; : \; s \mapsto
    u(s,\cdot) \in L^2([0,L/\e]; \ov{C}^{m,\t}_{\varsigma})  \right\},
\end{equation}
and for $l \in \N$, we define similarly the functional space
\begin{equation}\label{eq:hkcmasfl}
    \ov{H}^l(\ov{C}^{m,\t}_{\varsigma}) = \left\{ u : [0,L/\e] \times \R^{n-1} \to \C
    \; : \; s \mapsto
    u(s,\cdot) \in H^l([0,L/\e]; C^{m,\t}_{\varsigma})  \right\}.
\end{equation}
The weights here are suited for studying functions which decay in
$y$ like $e^{- |y|}$, as the fundamental solution of $-
\D_{\R^{n-1}} u +  u = 0$. The parameter $\varsigma < 1$ has been
introduced to allow some flexibility in the decay rate. When dealing
with functions belonging to the above three spaces, the symbols $\|
\cdot \|_{\ov{C}^{m,\t}_{\varsigma}}$, $\| \cdot
\|_{\ov{L}^2(\ov{C}^{m,\t}_{\varsigma})}$, $\| \cdot
\|_{\ov{H}^l(\ov{C}^{m,\t}_{\varsigma})}$ will denote norms induced
by formulas \eqref{eq:cmas}, \eqref{eq:l2cmasfl} and
\eqref{eq:hkcmasfl}. Also, we keep the same notation for the norms
when considering functions defined on subsets of $[0,L/\e] \times
\R^{n-1}$.

We next consider some positive constants $\hat{V}, \hat{f}, \hat{h},
\hat{k}$ which satisfy the relations in \eqref{eq:ovhovk}. If $\d$
 and $\hat{K}_\d$ are as in the previous subsection and $\ov{\d}$ as
in Section \ref{s:LS}, letting
$$
  D_{L,\e} = [0,L/\e] \times B_{\e^{-\ov{\d}}+1}(0) \subseteq
 [0,L/\e] \times \R^{n-1},
$$
we define the space of functions
\begin{equation*}
    \hat{H}_\e : = \left\{ \phi \; : \; \Re \int_{D_{L,\e}} \ov{\phi}(s,y)
    e^{- i \hat{f} s} v\big( s, y/\sqrt{\hat{V}} \big)
     = 0 \qquad \hbox{ for all } v \in \hat{K}_\d \right\}.
\end{equation*}
This conditions represents, basically, orthogonality with respect to
$\hat{K}_\d$ (multiplied by the phase factor), when the function
$\phi$ is scaled in $y$ by $\sqrt{V}$. This is a choice made for
technical reasons, which will be helpful in Proposition
\ref{p:truefndec}. We next have the following result, related to
Proposition \ref{p:inv} once we scale $y$.

\begin{lem}\label{l:solvortflat}
Let $\frac 12 \leq \t < 1$ and $\varsigma \in (0,1)$. Then, for $\d$
 small there exists a
positive constant $C$, depending only on $p$, $\t$, $\varsigma$,
$L$, $\hat{V}$ and $\hat{f}$, such that the following property
holds: for $\hat{f}$ small, for $\e \to 0$ and for any function
$\mathfrak{b} \in \ov{L}^2(\ov{C}^{\t}_{\varsigma})$ there exist
  $u \in \hat{H}_\e$, and $\underline{v} \in
\hat{K}_\d$ such that, in $D_{L,\e}$
\begin{equation}\label{eq:solvort}
  \left\{
    \begin{array}{ll}
      - \frac{1}{\hat{V}} \pa^2_{s s} u - \D_y u + u - \frac{\hat{h}^{p-1}}{\hat{V}}
      U(y \hat{k}/\sqrt{\hat{V}})^{p-1} u - (p-1)
\frac{\hat{h}^{p-1}}{\hat{V}} U(y \hat{k}/\sqrt{\hat{V}})^{p-1} e^{- i \hat{f} s}
      \Re(e^{- i \hat{f} s} \ov{u})
      &  \\ = \mathfrak{b} + e^{- i \hat{f} s} \underline{v};
 \\ & \\ u = 0 \qquad \hbox{ on } \; \partial D_{L,\e}, &
    \end{array}
  \right.
\end{equation}
(notice that $\underline{v}$ above is intended scaled in $y$) and
such that we have the estimates
\begin{equation}\label{eq:solvortest}
  \|u\|_{\ov{L}^2(\ov{C}^{2,\t}_{\varsigma})} + \|u\|_{\ov{H}^1(\ov{C}^{1,\t}_{\varsigma})}
  + \|u\|_{\ov{H}^2(\ov{C}^{\t}_{\varsigma})} \leq \frac{C}{\d^2}
  \inf_{v \in \hat{K}_\d}  \|\mathfrak{b} + e^{- i \hat{f} s} v
  \|_{\ov{L}^2(\ov{C}^{\t}_{\varsigma})};
\end{equation}
\begin{equation}\label{eq:solvortest'}
    \|\underline{v}\|_{\ov{L}^2(\ov{C}^{\t}_{\varsigma})} \leq C
    \|\mathfrak{b}\|_{\ov{L}^2(\ov{C}^{\t}_{\varsigma})}.
\end{equation}
\end{lem}

\

\begin{pf} First of all we observe that a solution to
\eqref{eq:solvort} of class $L^2$ exists. In fact, replacing
$D_{L,\e}$ with $[0,L/\e] \times \R^{n-1}$, this would simply follow
from Proposition \ref{p:inv} with $V \equiv \hat{V}$. However, since
the functions in $\hat{K}_\d$ decay exponentially to zero as $|y|
\to + \infty$ the Dirichlet boundary conditions do not affect the
solvability property: for more details see for example \cite{malm},
Lemma 5.5. Notice that indeed, by \eqref{eq:decoU} and Proposition
\ref{p:alphamu} {\bf e}, the elements of $K_\d$ decay at the rate
$e^{- \hat{k} |y|}$, and by \eqref{eq:ovhovk} $\hat{k} >
\sqrt{\hat{V}}$. In particular
$\|\underline{v}\|_{\ov{L}^2(\ov{C}^{\t}_{\varsigma})}$ is finite
and \eqref{eq:solvortest'} holds. We also have \eqref{eq:solvortest}
replacing the left-hand side by the $L^2$ norm of $u$. We divide the
rest of the proof into two steps.

\

\noindent {\bf Step 1: $u \in \ov{L}^2(\ov{C}^{\t}_{\varsigma})$ and
$\|u\|_{\ov{L}^2(\ov{C}^{\t}_{\varsigma})} \leq \frac{C}{\d^2}
\|\mathfrak{b}\|_{\ov{L}^2(\ov{C}^{\t}_{\varsigma})}$.} \quad We set
$u = e^{- i \hat{f} s} v$ and $\mathfrak{c} = e^{- i \hat{f} s}
  ({\mathfrak b} + \underline{v})$,
so $v$ satisfies
\begin{equation*}
   \left\{
     \begin{array}{ll}
      - \D v + (1 + \hat{f}^2/\hat{V}) v + 2 i \hat{f}/\hat{V} \pa_s v - \frac{\hat{h}^{p-1}}{\hat{V}}
      U(y \hat{k}/\sqrt{\hat{V}})^{p-1} v \\ - (p-1) \frac{\hat{h}^{p-1}}{\hat{V}} U(y \hat{k}/\sqrt{\hat{V}})^{p-1}
      \Re(\ov{v})  =
{\mathfrak c}  & \hbox{ in } B_{\e^{-\ov{\d}}+1}(0), \\
       v = 0 & \hbox{ on } \partial B_{\e^{-\ov{\d}}+1}(0).
     \end{array}
   \right.
\end{equation*}
We now use a Fourier decomposition in the variable $s$: setting
$$
  {\mathfrak c}(s,y) = \sum_j {\mathfrak c}_j(y) e^{i j \e s}; \qquad \quad
  v(s,y) = \sum_j v_j(y) e^{i j \e s},
$$
(here we are assuming for simplicity that $L = 2 \pi$) we see that
each ${\mathfrak c}_j$ belongs to $C^\t_{\varsigma,\hat{V}}$, that
\begin{equation}\label{eq:decffour}
    \sum_j \|{\mathfrak c}_j\|_{\ov{C}^{\t}_{\varsigma}}^2 = \frac1 \e
\|{\mathfrak c}\|_{\ov{L}^2(\ov{C}^{\t}_{\varsigma})}^2; \qquad \quad
\sum_j \|v_j\|_{\ov{C}^{\t}_{\varsigma}}^2 = \frac1 \e
\|v\|_{\ov{L}^2(\ov{C}^{\t}_{\varsigma})}^2 \leq \frac{C}{\e \d^2}
\|\mathfrak{b}\|_{\ov{L}^2(\ov{C}^{\t}_{\varsigma})},
\end{equation}
and that each $v_j$ solves
\begin{equation}\label{eq:eqforuj}
   \left\{
     \begin{array}{ll}
       - \D_{y} v_j + \left( 1 + \frac{\hat{f}^2 + \e^2 j^2 - 2 \hat{f} \e j}{\hat{V}}
     \right) v_j - \frac{\hat{h}^{p-1}}{\hat{V}} U(y \hat{k}/\sqrt{\hat{V}})^{p-1} v_j
     \\ - (p-1) \frac{\hat{h}^{p-1}}{\hat{V}} U(y \hat{k}/\sqrt{\hat{V}})^{p-1}
      \Re(\ov{v}_j)  = {\mathfrak c}_j  & \hbox{ in } B_{\e^{-\ov{\d}}+1}(0), \\
       v_j = 0 & \hbox{ on } \partial B_{\e^{-\ov{\d}}+1}(0).
     \end{array}
   \right.
\end{equation}
From elliptic regularity theory, we find that for any $R
> 0$ there exists a constant $C$ depending only on $R$, $p$ and $\t$
such that $\|v_j\|_{C^\t(B_R)} \leq C \|{\mathfrak
c}_j\|_{\ov{C}^\t_{\varsigma}} + C \|v_j\|_{L^2}$. Now we choose $R$
(depending on $p$ and $\varsigma$) so large that $p
\frac{\hat{h}^{p-1}}{\hat{V}} U^{p-1}(y \hat{k}/\sqrt{\hat{V}}) <
\frac{1}{4} (1-\varsigma)$ for $|y| \geq \frac R 2$, and a smooth
radial cutoff function $\hat{\chi}$ such that $\hat{\chi}(y) = 1$
for $|y| \leq \frac R 2$, and $\hat{\chi}(y) = 0$ for $|y| \geq R$.
Next, we write equation \eqref{eq:eqforuj} as
\begin{equation*}
    \left\{
      \begin{array}{ll}
        - \D_{y} v_j + \left( 1 + \frac{\hat{f}^2 + \e^2 j^2 - 2 \hat{f}
    \e j}{\hat{V}} \right)
        v_j - (1 - \hat{\chi}) p \frac{\hat{h}^{p-1}}{\hat{V}}
   U(y \hat{k}/\sqrt{\hat{V}})^{p-1} v_j = \tilde{{\mathfrak c}}_j
        + \hat{\chi} p \frac{\hat{h}^{p-1}}{\hat{V}}
   U(y \hat{k}/\sqrt{\hat{V}})^{p-1} v_j, & \\
       v_j = 0 \qquad \hbox{ on } \partial B_{\e^{-\ov{\d}}+1}(0). &
     \end{array}
   \right.
\end{equation*}
We notice that the first linear coefficient of $v_j$ is bounded
below (uniformly in $j$) by $1$. Therefore, using the Green's
representation formula, the maximum principle and our choice of $R$
(see Lemma \ref{l:decgreenflat}) for any $\varsigma' < \varsigma$ we
then have the estimate
$$
 \|v_j\|_{\ov{C}^\t_{\varsigma'}} \leq C (\|{\mathfrak
 c}_j\|_{\ov{C}^\t_{\varsigma}} + \|v_j\|_{L^2})
$$
for some fixed constant $C$ depending only on $p$, $\varsigma$ and
$\t$. Taking the square and summing over $j$ we get
$$
 \|u\|_{\ov{L}^2(\ov{C}^{\t}_{\varsigma'})}^2 =
 \|v\|_{\ov{L}^2(\ov{C}^{\t}_{\varsigma'})}^2 \leq C
 \|{\mathfrak b}\|_{\ov{L}^2(\ov{C}^{\t}_{\varsigma})}^2
 + C \|v\|_{L^2}^2 \leq \|{\mathfrak b}\|_{\ov{L}^2(\ov{C}^{\t}_{\varsigma})}^2.
$$
We next want to replace in the last formula $\varsigma'$ with
$\varsigma$.  Rewrite \eqref{eq:solvort} as
\begin{equation*}
  \left\{
    \begin{array}{ll}
      - \frac{1}{\hat{V}} \pa^2_{s s} u - \D_y u + u = \hat{\mathfrak{c}} : = \frac{\hat{h}^{p-1}}{\hat{V}}
      U(y \hat{k}/\sqrt{\hat{V}})^{p-1} u & \\ + (p-1) \frac{\hat{h}^{p-1}}{\hat{V}}
U(y \hat{k}/\sqrt{\hat{V}})^{p-1} e^{- i \hat{f} s}
      \Re(e^{- i \hat{f} s} \ov{u}) + \mathfrak{b} + e^{- i \hat{f} s} \underline{v};
      &  \\ \\  u = 0 \qquad \hbox{ on } \partial D_{L,\e}.   &
    \end{array}
  \right.
\end{equation*}
Using the same procedure as above, write $\hat{{\mathfrak c}}(s,y) =
\sum_j \hat{{\mathfrak c}}_j(y) e^{i j \e s}$ and $u(s,y) = \sum_j
u_j(y) e^{i j \e s}$.

We consider now the function $U^{p-1} u_j$: by \eqref{eq:decoU}, if
we choose $\varsigma' + \frac{(p-1) \hat{k}}{\sqrt{\hat{V}}} >
\varsigma$, it follows from the above estimates that
$\|\hat{\mathfrak{c}}\|_{\ov{L}^2(\ov{C}^{\t}_{\varsigma})}$ is
finite and that
\begin{equation*}
    \sum_j \|\hat{{\mathfrak c}}_j\|_{\ov{C}^{\t}_{\varsigma}}^2 \leq  \frac C \e
\|{\mathfrak b}\|_{\ov{L}^2(\ov{C}^{\t}_{\varsigma})}^2.
\end{equation*}
Moreover $u_j$ satisfies
$$
 \left\{
   \begin{array}{ll}
     - \D u_j + \left(1 + \frac{\e^2 j^2}{\hat{V}} \right) u_j =
   \hat{{\mathfrak c}}_j  & \hbox{ in } B_{\e^{-\ov{\d}}+1}(0), \\
       u_j = 0 & \hbox{ on } \partial B_{\e^{-\ov{\d}}+1}(0).
 \end{array} \right.
$$
Also, it is easy to show
\begin{eqnarray}\label{eq:decnormsu} \nonumber
  & & \|u\|_{\ov{L}^2(\ov{C}^{2,\t}_{\varsigma})}^2 +
  \|u\|_{\ov{H}^1(\ov{C}^{1,\t}_{\varsigma})}^2
  + \|u\|_{\ov{H}^2(\ov{C}^{\t}_{\varsigma})}^2 \\
  & = & \frac 1 \e \sum_j \left[ \|u_j\|^2_{\ov{C}^{2,\t}_{\varsigma}} +
  (1 + \e^2 j^2) \|u_j\|^2_{\ov{C}^{1,\t}_{\varsigma}}  +
  (1 + \e^2 j^2 + \e^4 j^4) \|u_j\|^2_{\ov{C}^\t_{\varsigma}} \right],
\end{eqnarray}
and therefore we are reduced to find estimate
$\|u_j\|^2_{\ov{C}^{2,\t}_{\varsigma}}$,
$\|u_j\|^2_{\ov{C}^{1,\t}_{\varsigma}}$ and
$\|u_j\|^2_{\ov{C}^\t_{\varsigma}}$, done in the next step.

\

\noindent {\bf Step 2: proof concluded.} \quad We now set $a_j = 1 +
\frac{\e^2 j^2}{\hat{V}}$, and  $v_j(y) = u_j\left(
\frac{y}{\sqrt{a}_j} \right)$. Then, from a change of variables we
have the equation
$$
  \left\{
    \begin{array}{ll}
      - \D v_j(y) + v_j(y) = \hat{F}_j(y) := \frac{1}{a_j}
  \hat{{\mathfrak c}}_j \left( \frac{y}{\sqrt{a}_j} \right) & \hbox{ in }
   B_{\sqrt{a_j} (\e^{-\ov{\d}}+1)}(0), \\
       v_j = 0 & \hbox{ on } \partial B_{\sqrt{a_j} (\e^{-\ov{\d}}+1)}(0).
 \end{array} \right.
$$
Notice that $a_j > 0$ stays bounded from below independently of $j$,
and therefore by a scaling argument (and some elementary
inequalities) one finds
\begin{eqnarray}\label{eq:bello} \nonumber
  \sup_{y,z \in B_1(x)} |\hat{F}_j(y) - \hat{F}_j(z)| & = & \frac{1}{a_j}
  \sup_{y,z \in B_{1/\sqrt{a_j}} (x)} \left| \hat{\mathfrak{c}}_j(y/\sqrt{a_j}) -
  \hat{\mathfrak{c}}_j(z/\sqrt{a_j}) \right| \\
  & \leq & \frac{C}{a_j}
  \sup_{y,z \in B_1 (x)} \left| \hat{\mathfrak{c}}_j(y/\sqrt{a_j}) -
  \hat{\mathfrak{c}}_j(z/\sqrt{a_j}) \right| \\ & \leq & \nonumber \frac{C
  \|\hat{{\mathfrak c}}_j\|_{\ov{C}^\t_{\varsigma}}}{a_j^{1 + \frac \t 2}}
  \sup_{y,z \in B_1(x)} |y-z|^\t e^{- \frac{\varsigma |x|}{\sqrt{a_j}}}
\end{eqnarray}
where $C$ depends on $\t$ only, and hence we get
\begin{equation}\label{eq:normtildeFj}
    \|\hat{F}_j\|_{\ov{C}^\t_{\frac{\varsigma}{\sqrt{a_j}}}} \leq
  \frac{C}{a_j^{1 + \frac \t 2}} \|\hat{{\mathfrak c}}_j\|_{\ov{C}^\t_{\varsigma}}.
\end{equation}
Now Lemma \ref{l:decgreenflat} implies that
$\|v_j\|_{\ov{C}^{2,\t}_{{\varsigma/\sqrt{a_j}}}} \leq
\frac{C}{a_j^{1 + \frac \t 2}} \|\hat{{\mathfrak
c}}_j\|_{\ov{C}^\t_{\varsigma}}$. From this estimate, we will obtain
next some control on $u_j$ by scaling back the variables.

We consider an arbitrary $x \in \R^{n-1}$: similarly as before we
have
$$
  \sup_{y,z \in B_1(x)} \frac{|u_j(y) - u_j(z)|}{|y-z|^\t} = \sup_{y,z \in
  B_{1}(x)} \frac{|v_j(\sqrt{a_j} y) - v_j(\sqrt{a_j} z)|}{|y-z|^\t}.
$$
Since $a_j$ can be arbitrarily large, we cannot evaluate the
difference $v_j(\sqrt{a_j} y) - v_j(\sqrt{a_j} z)$ directly using
the weighted norm in the definition \eqref{eq:cmas} (as we did for
the first inequality in \eqref{eq:bello}), since the two points
$\sqrt{a_j} y$ and $\sqrt{a_j} z$ might not belong to the same unit
ball. We avoid this problem choosing $[\sqrt{a_j}]$ (the integer
part of $\sqrt{a_j}$) points $(y^l)_l$ lying on the segment
$[\sqrt{a_j} y, \sqrt{a_j} z]$ at equal distance one from each
other, and using the triangular inequality. Now the distance of two
consecutive points $y^l$ and $y^{l+1}$ will stay uniformly bounded
from above, and the minimal norm of the $y^l$'s is bounded from
below by $C^{-1} \sqrt{a_j} (|x|-1)$. Therefore, adding
$[\sqrt{a_j}]$ times the inequality and using \eqref{eq:normtildeFj}
we obtain
\begin{eqnarray*}
  \sup_{y,z \in B_1(x)} \frac{|u_j(y) - u_j(z)|}{|y-z|^\t} & \leq & \frac{C
  \sqrt{a_j}}{|y-z|^\t} \frac{C}{a_j^{1 + \frac \t 2}} \left| \frac{y-z}{\sqrt{a_j}}
  \right|^\t e^{- \sqrt{a_j} \varsigma  \frac{(|x|-1)}{\sqrt{a_j}}}
  \|\hat{{\mathfrak c}}_j\|_{\ov{C}^\t_{\varsigma}} \\
  & \leq & \frac{C}{a_j^{1 + \t - \frac 12}} e^{- \varsigma |x|}
  \|\hat{{\mathfrak c}}_j\|_{\ov{C}^\t_{\varsigma}} \leq \frac{C}{a_j}
  e^{-  \varsigma |x|} \|\hat{{\mathfrak c}}_j\|_{\ov{C}^\t_{\varsigma}} \leq \frac{C}{a_j}
  e^{-  \varsigma |x|} \|\tilde{\mathfrak c}_j\|_{\ov{C}^\t_{\varsigma}},
\end{eqnarray*}
since we chose $\t \geq \frac 12$ and since $a_j$ is uniformly
bounded from below. Similarly, taking first and second derivatives
we find that
$$
  \sup_{y,z \in B_1(x)} \frac{|\n u_j(y) - \n u_j(z)|}{|y-z|^\t}
  \leq \frac{C}{\sqrt{a_j}} e^{- \varsigma  |x|}
  \|\tilde{\mathfrak c}_j\|_{\ov{C}^\t_{\varsigma}};
$$
$$
 \sup_{y,z \in B_1(x)} \frac{|\n^2 u_j(y) - \n^2 u_j(z)|}{|y-z|^\t} \leq C e^{-
 \varsigma |x|} \|\tilde{\mathfrak c}_j\|_{\ov{C}^\t_{\varsigma}}
$$
where, again, $C$ depends only on $\t$. Recalling that $a_j =
\hat{V} + \e^2 j^2$, we have in this way proved that
$$
  \|u_j\|^2_{\ov{C}^{2,\t}_{\varsigma}} \leq C \|\tilde{\mathfrak c}_j\|_{\ov{C}^\t_{\varsigma}}^2;
  \qquad \|u_j\|^2_{\ov{C}^{1,\t}_{\varsigma}} \leq \frac{C}{1 + \e^2 j^2}
  \|\tilde{\mathfrak c}_j\|_{\ov{C}^\t_{\varsigma}}^2; \qquad \|u_j\|^2_{\ov{C}^\t_{\varsigma}} \leq
  \frac{C}{(1 + \e^2 j^2)^2} \|\tilde{\mathfrak c}_j\|_{\ov{C}^\t_{\varsigma}}^2.
$$
Now the conclusion follows from \eqref{eq:decffour},
\eqref{eq:decnormsu}, the last formula and the fact that
$$\|\tilde{{\mathfrak c}}\|_{\ov{L}^2(\ov{C}^{\t}_{\varsigma})} \leq C
\inf_{v \in \hat{K}_\d} \|\mathfrak{b} + v\|_{\ov{L}^2(\ov{C}^{\t}_{\varsigma})},
$$
see the beginning of Step 1.
\end{pf}

\

\noindent We next consider the operator $L_\e$  in $\tilde{D}_\e$,
see \eqref{eq:defLe}, acting on a suitable subset of
$H_{\tilde{D}_\e}$ (verifying an orthogonality condition similar to
\eqref{eq:ovHe}). We want to allow some flexibility in the choice of
approximate solutions: to do this we consider a normal section
$\Phi$ to $\g$ which verifies the following two conditions
\begin{equation}\label{eq:condPhi}
    \Phi \in span \left\{ h^{\frac{p+1}{4}} \varphi_j \; : \; j = 0, \dots, \frac{\d}{\e}
    \right\}; \qquad \qquad \|\Phi\|_{H^2(0,L)} \leq c_1 \e.
\end{equation}
Here $(\varphi_j)_j$ are as in \eqref{eq:eigen}, while $c_1$ is a
large constant to be determined later. Notice that by
\eqref{eq:condPhi} we have $\|\Phi'''\|_{L^2[0,L]} \leq C$, which
implies $\|\Phi''\|_{L^\infty} \leq C$, so also \eqref{eq:bdPhi}
holds true. This will allow us, in the next section, to apply
Proposition \ref{p:exvarphi}. Next, we define the  variables
\begin{equation}\label{eq:defzz}
    z = y - \Phi(\e s).
\end{equation}
In the above coordinates $(s,z)$, we will consider the approximate
solution
\begin{equation}\label{eq:apslll}
    \tilde{\psi}_{\e} = e^{- i \frac{\tilde{f}(\e s)}{\e}} \eta_\e
  \left( h(\e s) U(k(\e s) z) + U_1(s,z) \right)
  := \tilde{\psi}_{0,\e} + \hat{\psi}_\e,
\end{equation}
where $\e s = \ov{s}$, and where $\tilde{f}$, $U_1$ satisfy, for
some fixed $C > 0$ and $\t \in (0,1)$
\begin{equation}\label{eq:tildefU1}
    \|\tilde{f} - f\|_{H^2([0,L])} \leq C \e^2; \qquad \qquad
    |U_1|(s,z) \leq C \e (1 + |z|^C) e^{- k(\ov{s}) |z|};
\end{equation}
\begin{equation}\label{eq:tildefU2}
     \| \, |h U(k \cdot) + U_1|^{p-1} - |h U(k \cdot)|^{p-1}\|_{C^\t}
     \leq C \e \qquad \hbox{ in } \tilde{D}_\e.
\end{equation}
With this choice of $\tilde{\psi}_{\e}$, we are going to study the
analogue of Lemma \ref{l:solvortflat} for $L_\e$, see
\eqref{eq:defLe}, using a perturbation method.

To state our final result we need to introduce some more notation.
Recalling the definition in \eqref{eq:cmas}, still using the
coordinates $(s,z)$, for $\t \in (0,1)$ and $\varsigma > 0$ we
define the function space
\begin{equation}\label{eq:l2cmas}
    L^2(C^{m,\t}_{\varsigma,V}) = \left\{ u : \tilde{D}_\e
    \to \C \; : \; s \mapsto u(s, \cdot/\sqrt{V(\e s)}) \in L^2([0,L/\e];
C^{m,\t}_{\varsigma,1})  \right\}.
\end{equation}
Also, for $m \in \N$, we define similarly
\begin{equation}\label{eq:hkcmas}
    H^l(C^{m,\t}_{\varsigma,V}) = \left\{ u : \tilde{D}_\e
    \to \C  \; : \; s \mapsto u(s,\cdot/\sqrt{V(\e s)})
\in H^l([0,L/\e]; C^{m,\t}_{\varsigma,1})  \right\}.
\end{equation}
We next let $\tilde{K}_\d$ be the counterpart of $K_\d$ (see
\eqref{eq:Kd}), when we replace the coordinates $y$ by  $z$.
Finally, we denote by $\ov{H}_\e$ the following subspace of
functions
\begin{equation}\label{eq:orth}
    \ov{H}_\e : = \left\{ \phi \in H_{\tilde{D}_\e} \; : \; \Re \int_{\tilde{D}_\e}
    e^{- i \frac{\tilde{f}(\e s)}{\e}} v \, \ov{\phi} = 0
    \qquad \hbox{ for all } v \in \tilde{K}_\d \right\}.
\end{equation}
Defining
\begin{equation}\label{eq:normstar}
    \|\cdot\|_{\varsigma,V} := \|\cdot\|_{L^2(C^{2,\t}_{\varsigma,V})} +
  \|\cdot\|_{H^1(C^{1,\t}_{\varsigma,V})}
  + \|\cdot\|_{H^2(C^{\t}_{\varsigma,V})},
\end{equation}
we have then the following result (recall the definition of
$\tilde{D}_\e$ in \eqref{eq:tildeDe}).

\begin{pro}\label{p:truefndec}
Suppose $0 < \varsigma < 1$ and $\frac 12 \leq \t < 1$.
  Suppose $\tilde{\psi}_\e$ is as in \eqref{eq:apslll}, with $\tilde{f}, U_1$
satisfying  \eqref{eq:tildefU1}. Then, if $K^2(\e s) = V(\e s)$, if
$\mathcal{A}$ and $\d$ are sufficiently small, in the limit $\e \to
0$ the following property holds: for any  function $b \in
L^2(C^{\t}_{\varsigma,V})$ there exist $\tilde{u} \in \ov{H}_\e$,
and $\underline{\tilde{v}} \in \tilde{K}_\d$ such that
\begin{equation}\label{eq:solvortfin}
  \left\{
    \begin{array}{ll}
      - \D_{g_\e} \tilde{u} + V(\e x) \tilde{u} - |\tilde{\psi}_\e|^{p-1} \tilde{u}
   - (p-1) |\tilde{\psi}_\e|^{p-3} \tilde{\psi}_\e \Re
   (\tilde{\psi}_\e \ov{\tilde{u}})  = b + e^{- i \frac{\tilde{f}(\e s)}{\e}} \underline{\tilde{v}}
   & \hbox{ in } \tilde{D}_\e \\[2mm]
     \tilde{u} = 0 & \hbox{ on } \partial \tilde{D}_\e
    \end{array}
  \right.
\end{equation}
is solvable, and such that for every $\varsigma' < \varsigma$ there
exists some $C > 0$ for which we have the estimates
\begin{equation}\label{eq:solvortestfin}
\|\tilde{u}\|_{\varsigma',V} \leq \frac{C}{\d^2} \inf_{\tilde{v} \in \tilde{K}_\d}
\|b + e^{- i \frac{\tilde{f}(\e s)}{\e}} \tilde{v} \|_{L^2(C^{\t}_{\varsigma,V})};
\qquad \quad \|\underline{\tilde{v}}\|_{L^2(C^{\t}_{\varsigma,V})}  \leq C
\|b\|_{L^2(C^{\t}_{\varsigma,V})}.
\end{equation}
\end{pro}

\begin{pf} We divide the proof into two steps.

\

\noindent {\bf Step 1: solvability of \eqref{eq:solvortfin}}. First
of all we notice that, from Proposition \ref{p:inv} and from
elliptic regularity results, if $H_\e$ denotes the subspace of
function in $H^2(N \g_\e)$ satisfying \eqref{eq:ovHe}, then the
operator $L^1_\e$ is invertible from $(H_\e, \|\cdot\|_{H^2(N
\g_\e)})$ onto $(\Pi_\e L^2(N \g_\e), \|\cdot\|_{L^2(N \g_\e)})$;
moreover the norm of the inverse operator
 is bounded by $\frac{C}{\d^2}$.

By the comments at the beginning of the proof of Lemma
\ref{l:solvortflat}, we also deduce the following property. Given $b
\in L^2\left( \{ |y| \leq (\e^{-\ov{\d}}+1)/K(\e s) \} \right)$
there exist $u \in H^2\left( \{ |y| \leq (\e^{-\ov{\d}}+1)/K(\e s)
\} \right)$ and $\underline{v} \in K_\d$ such that
\begin{eqnarray}\label{eq:ciaoello} \nonumber
      \tilde{L}_\e u & := & - \D_{\hat{g}_\e} u + V(\e s) u - h(\e s)^{p-1}
      U(k(\e s) y)^{p-1} u \\ & - & (p-1) h(\e s)^{p-1} U(k(\e s) y)^{p-1} e^{- i
      \frac{\tilde{f}(\e s)}{\e}} \Re(e^{- i \frac{\tilde{f}(\e s)}{\e}} \ov{u})
      = b + e^{- i \frac{\tilde{f}(\e s)}{\e}} \underline{v} \qquad \hbox{ in }
     \{ |y| \leq (\e^{-\ov{\d}}+1)/K \};
\end{eqnarray}
\begin{equation*}
    \Re \int_{\{ |y| \leq (\e^{-\ov{\d}}+1)/K \}} \ov{u} e^{- i \frac{f(\e s)}{\e}} v
    dV_{\hat{g}_\e} = 0 \qquad \hbox{ for every } v \in K_\d.
\end{equation*}
Again, we have the estimates
\begin{equation}\label{eq:estest}
    \|u\|_{H^2(\{ |y| \leq (\e^{-\ov{\d}}+1)/K \})} \leq \frac{C}{\d^2}
   \|b\|_{L^2(\{ |y| \leq (\e^{-\ov{\d}}+1)/K \})};
    \quad  \|\underline{v}\|_{L^2(\{ |y| \leq (\e^{-\ov{\d}}+1)/K \})}
   \leq C \|b\|_{L^2(\{ |y| \leq (\e^{-\ov{\d}}+1)/K \})}.
\end{equation}
Using a perturbative argument, we show  that we can recover the same
invertibility result for \eqref{eq:solvortfin} where, compared to
\eqref{eq:ciaoello}, we need to substitute $y$ with $z$,
$\D_{\hat{g}_\e}$ with $\D_{\tilde{g}_\e}$, $f$ with $\tilde{f}$ and
$e^{- i\frac{f(\e s)}{\e}} h U$ with $\tilde{\psi}_\e$.

In fact, let us denote by $\Pi_y$ and $\Pi_z$ the orthogonal
projections in $L^2$ onto the orthogonal complements of the sets $\{
e^{- i \frac{f(\e s)}{\e}} v \; : \; v \in K_\d \}$, $\{ e^{- i
\frac{\tilde{f}(\e s)}{\e}} v \; : \; v \in \tilde{K}_\d \}$ with
respect to the scalar products induced by the metrics $\hat{g}_\e$
and $\tilde{g}_\e$ respectively. By \eqref{eq:condPhi}, Lemma  3.1
in \cite{mmm1} and \eqref{eq:tildefU1} for every $u \in
H^2(\tilde{D}_\e)$ and every $b \in L^2(\tilde{D}_\e)$ one has
\begin{equation*}
\| L_\e u - \tilde{L}_\e u \|_{L^2(\tilde{D}_\e)} \leq C(c_1) \e \|u\|_{H^2(\tilde{D}_\e)};
\qquad \quad \|\Pi_y b - \Pi_z b\|_{L^2(\tilde{D}_\e)} \leq C(c_1) \e \|b\|_{L^2(\tilde{D}_\e)},
\end{equation*}
where $C(c_1)$ is a positive constant which depends on $\g$, $V$ and
the constant $c_1$ in \eqref{eq:condPhi}.

From \eqref{eq:estest} and the last formula we deduce the
solvability of \eqref{eq:solvortfin}, together with the estimates
$\|\tilde{u}\|_{H^2(\tilde{D}_\e)} \leq \frac{C}{\d^2}
\|b\|_{L^2(\tilde{D}_\e)}$ and
$\|\tilde{\underline{v}}\|_{L^2(\tilde{D}_\e)} \leq C
\|b\|_{L^2(\tilde{D}_\e)}$.

\

\noindent {\bf Step 2: proof of \eqref{eq:solvortestfin}}. Recall
that the coordinates $y$ (see the beginning of this section) are not
global, since they are defined locally in $s$ by normal parallel
transport: the same holds of course for the coordinates $z$.
 Therefore, if we prolong the $z$'s along $\g_\e$, there will
be a discontinuity between $0$ and $L/\e$.

To reduce ourselves to the periodic case, as in Lemma
\ref{l:solvortflat}, we apply a rotation $R_\e = R_\e(\e s)$ to the
$z$ axes which makes the coordinates $\tilde{z} := R(\e s) z$
periodic in $s$. To compute the Laplace-Beltrami operator in the new
 coordinates $\tilde{z}$ one
should apply the chain rule in this way
$$
  \pa_{z_j} u = (R_\e)_{jl} \pa_{\tilde{z}_l} u; \qquad \pa^2_{s z_j}
  u = \e \pa_{\ov{s}} (R_\e)_{jl} \pa_{\tilde{z}_l} u +  (R_\e)_{jl}
  \pa^2_{s \tilde{z}_l} u; \qquad
  \pa^2_{z_j z_l} u =  R_{mj} R_{tl} \pa^2_{\tilde{z}_m \tilde{z}_t} u.
$$
In particular, since $R_\e$ is orthogonal, $\pa^2_{z_j z_j} u =
R_{mj} R_{tj} \pa^2_{\tilde{z}_m
  \tilde{z}_t} u = (R_\e)_{mj} (R_\e^{-1})_{jt} \pa^2_{\tilde{z}_m \tilde{z}_t}
u = \pa^2_{\tilde{z}_m \tilde{z}_m} u$, namely the main term in the
Laplacian stays invariant. Taking into account  Lemma 3.1 in
\cite{mmm1} and the last formulas, for $\varsigma'' \in
(\varsigma',\varsigma)$ one finds
\begin{equation}\label{eq:diffcoordztz}
    \|\D^{\tilde{z}}_{\tilde{g}_\e} u - \D^{z}_{\tilde{g}_\e}
u\|_{L^2(C^\t_{\varsigma'',k})} \leq C(c_1) \e \|u\|_{\varsigma'',V}.
\end{equation}
We use next a localization argument as in the proof of Proposition
\ref{p:inv}. If $\hat{s}_j$ and $\chi_{\eta,j}$ are as in that
proof, by \eqref{eq:tildefU1} we can find $\hat{\th}_j \in \R$ such
that $\frac{\tilde{f}(\e s)}{\e} - \hat{f}_j s - \hat{\th}_j =
O(\sqrt \e)$ in the support of $\chi_{\eta,j}$. If we set
$\D_{\R^n}^{(s,\tilde{z})} = \D_{\R^{n-1}}^{\tilde{z}} +
\pa^2_{ss}$, and if we scale the $\tilde{z}$ variables by $K(\e s) =
\sqrt{V(\e s)}$, the function $\chi_{\eta,j}(s) u(s,\tilde{z})$
(which is now periodic in $s$) satisfies the equation
$$
  \left\{
    \begin{array}{ll}
      - \frac{1}{\hat{V}_j} \pa^2_{ss} \chi_{\eta,j} u - \D_{\R^{n-1}}^{\tilde{z}}
   (\chi_{\eta,j} u) + \chi_{\eta,j} u -
  \frac{\hat{h}_j^{p-1}}{\hat{V}_j} U\big(\tilde{z} \hat{k}_j/\sqrt{\hat{V}_j}\big)^{p-1}
   \chi_{\eta,j} u  &  \\  - (p-1) \frac{\hat{h}_j^{p-1}}{\hat{V}_j}
  U\big(\tilde{z} \hat{k}_j/\sqrt{\hat{V}_j}\big)^{p-1} e^{- i (\hat{f} s
  + \hat{\th_j})} \Re(e^{- i (\hat{f} s + \hat{\th_j})} \chi_{\eta,j} \ov{u}) =
  \mathcal{F}_j & \hbox{ in } \{ |\tilde{z}| \leq (\e^{-\ov{\d}}+1)/K \};  \\
    \chi_{\eta,j} u = 0  & \hbox{ on }  \{ |\tilde{z}| = (\e^{-\ov{\d}}+1)/K \},
    \end{array}
  \right.
$$
where
\begin{eqnarray*}
  \mathcal{F}_j & = & \frac{1}{V(\ov{s})} \chi_{\eta,j} e^{- i \frac{\tilde{f}(\e s)}{\e}}
  (b + \underline{v}) + \frac{1}{V(\ov{s})} (\D^{\tilde{z}}_{\tilde{g}_\e} -
  \D_{\R^n}^{(s,\tilde{z})}) \chi_{\eta,j} u - \frac{1}{V(\ov{s})}
  (2 \n_{\tilde{g}_\e}^{\tilde{z}} u \cdot
  \n_{\tilde{g}_\e}^{\tilde{z}} \chi_{\eta,j} + u \D^{\tilde{z}}_{\tilde{g}_\e}
  \chi_{\eta,j}) \\ & + & \frac{1}{V(\ov{s})} (\hat{V}_j - V)
 \chi_{\eta,j} u + \frac{1}{V(\ov{s})} |\tilde{\psi}_\e|^{p-1}
  \chi_{\eta,j} u  + \frac{1}{V(\ov{s})} (p-1) |\tilde{\psi}_\e|^{p-2}
  \tilde{\psi}_\e \Re (\tilde{\psi}_\e \chi_{\eta,j} \ov{u}) \\ & - &
  \frac{\hat{h}_j^{p-1}}{\hat{V}_j} U\big(\tilde{z} \hat{k}_j/\sqrt{\hat{V}_j}\big)^{p-1}
   \chi_{\eta,j} u  - (p-1) \frac{\hat{h}_j^{p-1}}{\hat{V}_j}
  U\big(\tilde{z} \hat{k}_j/\sqrt{\hat{V}_j}\big)^{p-1} e^{- i (\hat{f} s
  + \hat{\th_j})} \Re(e^{- i (\hat{f} s + \hat{\th_j})} \chi_{\eta,j} \ov{u}).
\end{eqnarray*}
In the last formula, the functions $b, \underline{v}$, $V$ and
$\tilde{\psi}_\e$ are intended scaled in $\tilde{z}$ by $\sqrt{V(\e
s)}$. Reasoning as for \eqref{eq:projchiphionhatK}, from
\eqref{eq:orth} one finds that $\int_{\tilde{D}_\e} e^{- i \hat{f}_j
s - \hat{\th}_j} \hat{v} \chi_{\eta,j} \ov{\phi} = O(\d^2 +
\sqrt{\e}) \|\phi\|_{L^2(supp(\chi_{\eta,j}))}
\|\hat{v}\|_{L^2(\tilde{D}_\e)}$
 for every $\hat{v} \in \hat{K}_\d$. Moreover, as for
\eqref{eq:diffcoordztz} one can show that
$$
\|(\D^{\tilde{z}}_{\tilde{g}_\e} -  \D_{\R^n}^{(s,\tilde{z})})
\chi_{\eta,j} u\|_{\ov{L}^2(\ov{C}^\t_{\varsigma''})} \leq C(c_1) \e
\left( \|\chi_{\eta,j} u\|_{\ov{L}^2(\ov{C}^{2,\t}_{\varsigma''})} +
  \|\chi_{\eta,j} u\|_{\ov{H}^1(\ov{C}^{1,\t}_{\varsigma''})}
  + \|\chi_{\eta,j} u\|_{\ov{H}^2(\ov{C}^{\t}_{\varsigma''})} \right).
$$
Therefore, using Lemma \ref{l:solvortflat}, \eqref{eq:tildefU1} and
\eqref{eq:tildefU2} we obtain the estimate
\begin{eqnarray}\label{eq:fethi} \nonumber
  \|\chi_{\eta,j} u\|_{\ov{L}^2(\ov{C}^{2,\t}_{\varsigma''})} +
  \|\chi_{\eta,j} u\|_{\ov{H}^1(\ov{C}^{1,\t}_{\varsigma''})}
  + \|\chi_{\eta,j} u\|_{\ov{H}^2(\ov{C}^{\t}_{\varsigma''})} \leq
  \frac{C}{\d^2} \|\mathcal{F}_j\|_{\ov{L}^2(\ov{C}^{\t}_{\varsigma''})} \\ \leq \frac{C}{\d^2}
   \|\chi_{\eta,j} b \|_{\ov{L}^2(\ov{C}^{\t}_{\varsigma''})} +
   \frac{C}{\d^2} \|\chi_{\eta,j} \underline{v} \|_{\ov{L}^2(\ov{C}^{\t}_{\varsigma''})}
   \\ + C \sqrt{\e} \left( \|u\|_{\ov{L}^2(\ov{C}^{2,\t}_{\varsigma''}, supp(\chi_{\eta,j}))} +
  \|u\|_{\ov{H}^1(\ov{C}^{1,\t}_{\varsigma''}, supp(\chi_{\eta,j}))}
  + \|u\|_{\ov{H}^2(\ov{C}^{\t}_{\varsigma''}, supp(\chi_{\eta,j}))} \right), \nonumber
\end{eqnarray}
where the last symbols denote the restrictions of the weighted norms
to $supp(\chi_{\eta,j})$. Recall that the functions in the previous
formula have been scaled in $\tilde{z}$ by $\sqrt{V (\e s)}$:
therefore, from the uniform continuity of $V(\ov{s})$, for some $C >
0$ we have (recall that $\varsigma'' \in (\varsigma',\varsigma)$)
$$
  \frac 1 C \|\chi_{\eta,j} u\|_{\varsigma',V} \leq \|\chi_{\eta,j}
u (\cdot, \sqrt{V(\ov{s}) \, \cdot}) \|_{\ov{L}^2(\ov{C}^{2,\t}_{\varsigma''})} +
  \|\chi_{\eta,j} u (\cdot, \sqrt{V(\ov{s}) \, \cdot}) \|_{\ov{H}^1(\ov{C}^{1,\t}_{\varsigma''})}
  + \|\chi_{\eta,j} u (\cdot, \sqrt{V(\ov{s}) \, \cdot}) \|_{\ov{H}^2(\ov{C}^{\t}_{\varsigma''})}.
$$
A similar inequality holds for the restriction of $u$ to the support
of $\chi_{\eta,j}$, together with
$$
  \|\chi_{\eta,j} b(\cdot, \sqrt{V(\ov{s}) \, \cdot})
\|_{\ov{L}^2(\ov{C}^{\t}_{\varsigma''})} +
    \|\chi_{\eta,j} \underline{v} (\cdot, \sqrt{V(\ov{s}) \, \cdot})
\|_{\ov{L}^2(\ov{C}^{\t}_{\varsigma''})} \leq C \left(  \|\chi_{\eta,j}
 b\|_{L^2(C^{\t}_{\varsigma,V})}
  + \|\chi_{\eta,j} \underline{v} \|_{L^2(C^{\t}_{\varsigma,V})} \right).
$$
Using the last two inequalities, taking the square of
\eqref{eq:fethi} and summing over $j$, we can bring the last term in
the  right-hand side to the left, so we get
\eqref{eq:solvortestfin}.
\end{pf}

\section{Approximate solutions}\label{s:as}

In this section we construct some approximate solutions to
\eqref{eq:new} which depend on suitable parameters, and find
 rigorous estimates on the error terms.

As in the previous subsection, we let $y$ be a system of Fermi
coordinates in $N \g_\e$, and for a normal section $\Phi$ of $N
\g_\e$ of class $H^2$  we define the coordinates (see
\eqref{eq:defzz})
\begin{equation*}
    z = y - \Phi(\e s), \qquad \qquad z \in \R^{n-1}.
\end{equation*}
 By the results
in Subsection \ref{subsec:loc}, we will restrict our attention to
the set $\tilde{D}_\e$.

\begin{rem}\label{r:neglect}
In the spirit of Proposition \ref{p:sec2final}, we will work with
approximate solutions $\tilde{\psi}_\e$ supported in $\tilde{D}_\e$.
Therefore, using the above coordinates, $\tilde{\psi}_\e(s,z)$ has
to vanish for $|z|$ sufficiently large. This can be achieved by
defining formally $\tilde{\psi}_\e(s,z)$ on $N \g_\e$, and
multiplying it by a cutoff function $\eta_\e$ as in Subsection
\ref{subsec:loc}. However, since the functions we are dealing with
decay exponentially to zero as $|z| \to \infty$, the effect of this
cutoff on the expansions below is exponentially small in $\e$, and
it will turn out to be negligible for our purposes. Therefore, for
reasons of brevity and clarity, we will tacitely assume that
$\tilde{\psi}_\e(s,z)$ is multiplied by such a cutoff, without
writing it explicitly.
\end{rem}

\noindent Recall that in \eqref{eq:defse} we defined
$$
  S_\e(\psi) = - \D_g \psi + V(\e x) \psi - |\psi|^{p-1} \psi.
$$
We set $\tilde{f}_0(\ov{s}) = f(\ov{s}) + \e f_1(\ov{s})$, where $f$
is given in \eqref{eq:f'Cintr} and $f_1$ (depending on $\Phi$ and
$\mathcal{A}$) was defined at the end of Subsection 4.1 in
\cite{mmm1}), in order to satisfy the equation
\begin{equation}\label{eq:f1111}
    \pa_{\ov s}\left( \frac{h^2f_1'}{(p-1)k^{n+1}}\left[(p-1)h^{p-1}-
2 \s \mathcal{A}^2h^{2\s} \right] \right) = 2\mathcal{A} \left(
\frac{p-1}{2\th} - 1
  \right)  \partial_{\ov s} \langle {\bf H},\Phi
  \rangle.
\end{equation}
This equation is indeed solvable explicitly and the solution is
given by
\begin{equation}\label{eq:f1}
    f'_1=\frac{2\mathcal{A}(p-1)k^{n+1}}{(p-1)h^{p+1}- 2 \s
\mathcal{A}^2h^{2\s+2} } \left( \frac{p-1}{2\th} - 1
  \right)   \langle {\bf H},\Phi \rangle
  + \mathcal{A}' \frac{(p-1) k^{n+1}}{(p-1) h^{p+1} - 2 \s \mathcal{A}^2
  h^{2\s+2}},
\end{equation}
(we refer to \cite{mmm1} for the definition of $\mathcal{A}'$). If
$w_r$ and $w_i$ are smooth functions of $\ov{s}$ and $z$ we have,
formally
$$
  S_\e\left(e^{-i \frac{\tilde{f}_0(\e s)}{\e}} (h U(k z) + \e (w_r + i w_i))
  \right) = e^{-i \frac{\tilde{f}_0(\e s)}{\e}} \left( \e (\mathcal{R}_r
  + i \mathcal{R}_i) \right) + o(\e),
$$
for some quantities $\mathcal{R}_r, \mathcal{R}_i$ given in
Subsection 3.2 of \cite{mmm1} (where we refer to also for the
derivation of the last formula). $\mathcal{R}_r$ and $\mathcal{R}_i$
can be written as $\mathcal{R}_r = \mathcal{L}_r w_r -
\mathcal{F}_r$, $\mathcal{R}_i = \mathcal{L}_i w_i - \mathcal{F}_i$,
where
\begin{equation}\label{eq:Rr}
    \mathcal{F}_r =- 2f'f_1'hU- 2 f'^2 h U(kz) \langle
    {\bf H}, z + \Phi \rangle -
    h k \langle {\bf H}, \n U(kz) \rangle -  \langle \n^N V, z + \Phi \rangle h U(kz);
\end{equation}
\begin{equation}\label{eq:Ri}
    \mathcal{F}_i = - \left[ f'' h U(kz) + 2 f' h' U(kz) + 2
    f' h k' \n U(kz) \cdot z \right] + 2 \sum_j [\Phi'_j f' h k \pa_j
    U(kz)],
\end{equation}
and where the operators $\mathcal{L}_r$, $\mathcal{L}_i$ are defined
in \eqref{eq:LrLi}. Therefore, for canceling the errors of order
$\e$ we require $w_r$ and $w_i$ to be formally determined by the
equations $\mathcal{L}_r w_r = \mathcal{F}_r$, $\mathcal{L}_i w_i =
\mathcal{F}_i$.

Dividing the right-hand sides of \eqref{eq:Rr} and \eqref{eq:Ri}
into their even and odd parts (in the variables $z$), we obtained
that $w_r = w_{r,e} + w_{r,o}$, and $w_i = w_{i,e} + w_{i,o}$,
where
\begin{equation}\label{eq:we}
    w_{r,e} = \left[\frac{p-1}{\th} h^p \langle {\bf H}, \Phi \rangle+2f'f_1'h\right] \left(
  \frac{1}{(p-1) h^{p-1}} U(kz) + \frac{1}{2k} \n U(kz) \cdot z
  \right);
\end{equation}
\begin{equation}\label{eq:wiewio}
    w_{i,e} = \frac{p-1}{4} f' h' |z|^2 U(kz); \qquad \qquad
  w_{i,o} = - \sum_j \Phi'_j f' h z_j U(kz),
\end{equation}
and where $w_{r,o}$ is given implicitly by the equation
\begin{equation}\label{eq:wro}
    \mathcal{L}_r w_{r,o} = - 2 (f')^2 h U(kz) \langle {\bf H}, z \rangle
  - h k \sum_j H^j \partial_j U(kz) - \langle \n^N V, z \rangle h
  U(kz).
\end{equation}
As noticed in Subsection 3.2 in \cite{mmm1} (see also the
introduction here), the solvability of the last equation is
guaranteed by the stationarity condition \eqref{eq:eulerintr}.
Moreover, it is standard to check that $w_{r,o}$ has exponential
decay in $z$, as for the other correction terms. Defining
$$
  \tilde{\psi}_{1,\e} = e^{- i \frac{\tilde{f}_0(\e s)}{\e}} \left(
  h U(k z) + \e w_r + i \e w_i \right),
$$
from the expansions in Subsection 3.3 of \cite{mmm1} we can write
that
\begin{eqnarray}\label{eq:esSe}
e^{i \frac{\tilde{f}_0(\e s)}{\e}} S_\e(\psi_{1,\e})&=&
\e^2(\tilde{R}_{r,e}+\tilde{R}_{r,o})
+\e^2(\tilde{R}_{r,e,f_1}+\tilde{R}_{r,o,f_1})
\\&+&\e^2i(\tilde{R}_{i,e}+\tilde{R}_{i,o})
+\e^2i(\tilde{R}_{i,e,f_1}+\tilde{R}_{i,o,f_1})
+o(\e^2).\nonumber
\end{eqnarray}
In the last formula, the $\tilde{R}$'s represent the terms of order
$\e^2$ appearing in the expansion, see Subsection 3.3 of
\cite{mmm1}, while $o(\e^2)$ stands for the terms which are formally
of higher order. Here indeed we want to prove rigorous estimates, so
we want to be careful in treating the latter term.

To allow some more flexibility in the choice of approximate
solutions, we substitute the phase $\tilde{f}_0$ with the function
$\tilde{f} = f + \e f_1 + \e^2 f_2$, where $f_2$ is some function of
class $H^2$. On $\Phi$ and $f_2$ we assume the following conditions
for some constants $c_1, c_2$ to be determined later
\begin{equation}\label{eq:bdc1c2}
    \|\Phi\|_{H^2}\le c_1 \e; \qquad \qquad \|f_2\|_{H^2}\le c_2.
\end{equation}
Moreover, letting $\d$ be as in Subsections \ref{ss:appker} and
$\var_j, \o_j$ as in \eqref{eq:eigen}, we also assume that
\begin{equation}\label{eq:bdPhif2d}
    \Phi \in span \left\{ h^{\frac{p+1}{4}} \varphi_j \; : \; j = 0, \dots, \frac{\d}{\e} \right\};
 \qquad \quad f_2 \in span \left\{ h^{\frac 12} \o_j \; : \; j = 0, \dots, \frac{\d}{\e} \right\}.
\end{equation}
To deal with the resonance phenomenon mentioned in the introduction,
related to the components in $K_{3,\d}$ of the approximate kernel,
we add to the approximate solutions a function $v_\d$ like
\begin{equation}\label{eq:vD}
    v_\d = \b(\e s) Z_{\a(\e s)} + i \xi(\e s) W_{\a(\e s)}
\end{equation}
(see \eqref{eq:defaovs} and the lines after), with $\b, \xi$ given
by
\begin{equation}\label{eq:exprbeta}
    \b(\e s) = \sum_{j=-\frac{\d^2}{\e}}^{\frac{\d^2}{\e}} b_j \b_j(\e s); \qquad
   \qquad \xi = \sum_{j=-\frac{\d^2}{\e}}^{\frac{\d^2}{\e}} b_j \xi_j(\e s),
\end{equation}
where, we recall, $\xi_j$ solves \eqref{eq:asynul} and is related to
$\b_j$ by \eqref{eq:bbj}. Below, we will regard $\b$ as an
independent variable, and $\xi$ as a function of $\b$. Introducing
the norm
\begin{equation}\label{eq:normsharp}
     \|\b\|_\sharp := \bigg( \sum_{j=-\frac{\d^2}{\e}}^{\frac{\d^2}{\e}}
 b_j^2 (1 + |j|)^2 \bigg)^{\frac 12},
\end{equation}
we will assume later on that
\begin{equation}\label{eq:bdbe3}
    \|\b\|_{\sharp} \leq c_3 \e^2
\end{equation}
for some constant $c_3 > 0$ to be specified later.

We will look for approximate solutions of the form
\begin{equation}\label{eq:P2eeee}
    \tilde{\Psi}_{2,\e}(s,z):=e^{- i \frac{\widetilde{f}(\e s)}{\e}}
\bigg\{ h(\e s) U \left( k( \e
  s) z \right) + \e \left[ w_r + i w_i \right] +\e^2
  \widetilde{v}+ \e^{2} v_0  + v_\d \bigg\}.
\end{equation}
In this formula $\tilde{f}$ is as above, while $\widetilde{v}$ and
$v_0$ are corrections whose choice is given below, in order to
improve the accuracy of the approximate solutions.

Our goal is to estimate with some accuracy the quantity
$S_\e(\tilde{\Psi}_{2,\e})$: for simplicity, to treat separately
some terms in this expression, we will write $\tilde{\Psi}_{2,\e}$
as
\begin{equation}\label{eq:psi2edec}
    \tilde{\Psi}_{2,\e}(s,z)=\tilde{\Psi}_{1,\e}(s,z)+E(s,z)+F(s,z)+G(s,z),
\end{equation}
where $\tilde{\Psi}_{1,\e},\,E,\,F$ and $G$ are respectively defined
by
\begin{equation*}
    \tilde{\Psi}_{1,\e}(s,z):=e^{- i \frac{\widetilde{f}(\e s)}{\e}} \bigg\{
  h(\e s) U \left( k( \e s) z \right) + \e \left[ w_r + i w_i \right]
  \bigg\} := e^{-i \frac{\e^2 f_2(\e s))}{\e}}
  \tilde{\psi}_{1,\e};
\end{equation*}
\begin{equation*}
E(s,z):= \e^2 e^{- i \frac{\widetilde{f}(\e s)}{\e}}\widetilde{v};
 \qquad \qquad F(s,z):=\e^2e^{- i \frac{\widetilde{f}(\e
s)}{\e}}v_0;\qquad\quad G(s,z):=e^{- i \frac{\widetilde{f}(\e
s)}{\e}} v_\d.
\end{equation*}

\

\no To expand $S_\e(\Psi_{2,\e})$ conveniently, we can write
$$
  S_\e(\tilde{\Psi}_{2,\e}) = S_\e(\tilde{\Psi}_{1,\e}) + \mathfrak{A}_3 + \mathfrak{A}_4
  + \mathfrak{A}_5 + \mathfrak{A}_6,
$$
where $\mathfrak{A}_3, \dots, \mathfrak{A}_6$ are respectively the
linear terms in the equation which involve $E$, $F$ and $G$ (see
\eqref{eq:psi2edec}):
\begin{equation}\label{eq:defA3}
    \mathfrak{A}_3 = - \D_g E + V(\e x) E - |\tilde{\Psi}_{1,\e}|^{p-1} E - (p-1)
  |\tilde{\Psi}_{1,\e}|^{p-3} \tilde{\Psi}_{1,\e} \Re(\tilde{\Psi}_{1,\e}
  \ov{E});
\end{equation}
\begin{equation}\label{eq:defA4}
    \mathfrak{A}_4 = - \D_g F + V(\e x) F - |\tilde{\Psi}_{1,\e}|^{p-1} F - (p-1)
  |\tilde{\Psi}_{1,\e}|^{p-3} \tilde{\psi
  }_{1,\e} \Re(\tilde{\Psi}_{1,\e}
  \ov{F});
\end{equation}
\begin{equation}\label{eq:defA5}
    \mathfrak{A}_5 = - \D_g G + V(\e x) G - |\tilde{\Psi}_{1,\e}|^{p-1} G - (p-1)
  |\tilde{\Psi}_{1,\e}|^{p-3} \tilde{\Psi}_{1,\e} \Re(\tilde{\Psi}_{1,\e}
  \ov{G}),
\end{equation}
and where $\mathfrak{A}_6$ contains the contribution of the
nonlinear part
\begin{equation}\label{eq:defA7}
  \mathfrak{A}_6 = - |\tilde{\Psi}_{2,\e}|^{p-1} \tilde{\Psi}_{2,\e} + |\tilde{\Psi}_{1,\e}|^{p-1}
  (E+F+G) + (p-1) |\tilde{\Psi}_{1,\e}|^{p-3} \tilde{\Psi}_{1,\e}
  \Re(\tilde{\Psi}_{1,\e} (\ov{E}+ \ov{F} + \ov{G} )).
\end{equation}
Next we also write (tautologically)
\begin{equation}\label{eq:defA1}
  S_\e(\tilde{\Psi}_{1,\e}) = e^{-i \frac{\e^2 f_2(\e s)}{\e}}
  S_\e(\tilde{\psi}_{1,\e}) + \mathfrak{A}_1; \qquad \qquad
 \mathfrak{A}_1 = S_\e(\tilde{\Psi}_{1,\e}) - e^{-i \frac{\e^2 f_2(\e s)}{\e}}
  S_\e(\tilde{\psi}_{1,\e}),
\end{equation}
and set
\begin{eqnarray}\label{eq:defA2} \nonumber
  \mathfrak{A}_2 & = & e^{-i \frac{\tilde{f}(\e s)}{\e}} \bigg(
  e^{i \frac{\tilde{f}_0(\e s)}{\e}}S_\e(\tilde{\psi}_{1,\e})
   - \e^2(\tilde{R}_{r,o}+\tilde{R}_{r,e}) - \e^2(\tilde{R}_{r,o,f_1}
 +\tilde{R}_{r,e,f_1})
 \\ & - &   \e^2i(\tilde{R}_{i,e}
 +\tilde{R}_{i,o})
-\e^2i(\tilde{R}_{i,e,f_1}+\tilde{R}_{i,o,f_1})
\bigg),
\end{eqnarray}
so that $\mathfrak{A}_2$ represents the terms which are formally
of order $\e^3$ and higher in $S_\e(\tilde{\psi}_{1,\e})$
(multiplied by a phase factor). Therefore, from the definitions
\eqref{eq:defA3}-\eqref{eq:defA2} we find that
\begin{eqnarray}\label{eq:ssss} \nonumber
  S_\e(\tilde{\Psi}_{2,\e}) &=& e^{-i \frac{\tilde{f}(\e s)}{\e}} \bigg(
  \e^2(\tilde{R}_{r,o}+\tilde{R}_{r,e}) +
  \e^2(\tilde{R}_{r,o,f_1}
 +\tilde{R}_{r,e,f_1}) \\
   & + & \e^2i(\tilde{R}_{i,e}
 +\tilde{R}_{i,o})
+\e^2i(\tilde{R}_{i,e,f_1}+\tilde{R}_{i,o,f_1})
\bigg) + \mathfrak{A}_1+ \mathfrak{A}_2+ \mathfrak{A}_3+
\mathfrak{A}_ 4+ \mathfrak{A}_5+ \mathfrak{A}_6.
\end{eqnarray}

\

\no To estimate rigorously the $\mathfrak{A}_i$'s, we display the
first and second order derivatives of $\widetilde{\Psi}_{2,\e}$
\begin{eqnarray*}
\partial_s \tilde{\Psi}_{2,\e} & = & - i \widetilde{f}'(\e s)e^{- i
\frac{\widetilde{f}(\e s)}{\e}} \bigg[
  h(\e s) U (k (\e
  s) z) + \e \left[ w_r + i w_i \right] +\e^2 \widetilde{v}+\e^{2}v_0+ v_\d\bigg]
  \\ & + &  e^{- i \frac{\widetilde{f}(\e s)}{\e}}
  \bigg[ \e h'U(kz)+\e hk'\n U\cdot z+\e^2\pa_s w_r+i\e^2\pa_s w_i+\e^3 \pa_s
  \widetilde{v}+\e^3\pa_s v_0 + \pa_s v_\d \bigg]; \nonumber
\end{eqnarray*}
  \begin{equation*}
\partial_{j} \tilde{\Psi}_{2,\e} = e^{- i \frac{\widetilde{f}(\e s)}{\e}} \bigg[
  h(\e s) k(\e s)\pa_jU (k (\e
  s) z) + \e \left[ \pa_jw_r + i \pa_jw_i \right] +\e^2
  \pa_j\widetilde{v}+\e^{2}\pa_jv_0+ \pa_j v_\d \bigg];
\end{equation*}
\begin{eqnarray*}
    & & \partial^2_{ss} \tilde{\Psi}_{2,\e} = (-\widetilde{f}'^2 - i\e
    \widetilde{f}'')(\e s)e^{- i \frac{\widetilde{f}(\e s)}{\e}} \bigg[
  h(\e s) U (k (\e
  s) z) + \e \left[ w_r + i w_i \right] +\e^2 \widetilde{v}+\e^{2}v_0+a(\e s)Z(kz)  \bigg]
  \\ \nonumber & + &  e^{- i \frac{\widetilde{f}(\e s)}{\e}}
  \bigg[\e^2 h''U(kz)+2\e^2 h'k'\n U\cdot z+\e^2 hk''\n U\cdot z+\e^2 hk'^2
  \n^2U[z,z]+\e^3\pa^2_{ss} w_r
  +i\e^3\pa^2_{ss} w_i\\
  \nonumber&+&\e^4\pa^2_{ss} \widetilde{v}+\e^4\pa^2_{ss} v_0+ \pa^2_{ss} v_\d \bigg]\\
 \nonumber &-& 2 i\widetilde{f}'(\e s)e^{- i \frac{\widetilde{f}(\e s)}{\e}} \bigg[
 \e h'U(kz)+\e hk'\n U\cdot z+\e^2\pa_s w_r
  +i\e^2\pa_s w_i+\e^3 \pa_s \widetilde{v}+\e^3\pa_s v_0 + \pa_s v_\d \bigg]; \nonumber
\end{eqnarray*}
  \begin{equation*}
\partial^2_{jl} \tilde{\Psi}_{2,\e} = e^{- i \frac{\widetilde{f}(\e s)}{\e}} \bigg[
  h(\e s) k^2\pa^2_{jl}U (k (\e
  s) z) + \e \left[ \pa^2_{jl}w_r + i \pa^2_{jl}w_i \right] +\e^2
  \pa^2_{jl}\widetilde{v} +\e^{2}\pa^2_{jl}v_0+ \pa^2_{jl} v_\d \bigg];
\end{equation*}
\begin{eqnarray*}
  \partial^2_{sj} \tilde{\Psi}_{2,\e} & = & - i \widetilde{f}'(\e s)e^{- i
  \frac{\widetilde{f}(\e s)}{\e}} \bigg[
  h(\e s)k\pa_j U (k (\e
  s) z) + \e \left[ \pa_jw_r + i\pa_j w_i \right] +\e^2
  \pa_j\widetilde{v}+\e^{2}\pa_jv_0 \nonumber \\ & + & a(\e s)k\pa_jZ(kz)
    \bigg] +  e^{- i \frac{\widetilde{f}(\e s)}{\e}}
  \bigg[\e h'k\pa_jU(kz)+\e hk'\pa_j U+\e hkk'z_l\pa^2_{jl}U+\e^2\pa^2_{sj} w_r
  \\ & + & i\e^2\pa^2_{sj} w_i+\e^3 \pa^2_{sj} \widetilde{v}+
  \e^3\pa^2_{sj} v_0\nonumber
   + \pa^2_{sj} v_\d \bigg].\nonumber
\end{eqnarray*}

To simplify the expressions of the error terms, we introduce some
convenient notation. For any positive integer $q$, the two symbols
$\mathfrak{R}_q(\Phi,\Phi')$ and $\mathfrak{R}_q(\Phi,\Phi',\Phi'')$
will denote error terms satisfying the following bounds, for some
fixed constants $C, d$ (which depend on $q$, $c_1, c_2, c_3$ but not
on $\e$, $s$ and $\d$)
\begin{equation*}
    \left\{
      \begin{array}{ll}
        |\mathfrak{R}_q(\Phi,\Phi')| \leq C \e^q (1+|z|^d) e^{- k |z|}; &  \\
        |\mathfrak{R}_q(\Phi,\Phi') - \mathfrak{R}_q(\tPhi,\tPhi')| \leq C
    \e^q (1+|z|^d) [|\Phi-\tPhi| + |\Phi'-\tPhi'|] e^{- k |z|}, &
      \end{array}
    \right.
\end{equation*}
while the  term $\mathfrak{R}_q(\Phi,\Phi',\Phi'')$ (which involves
also second derivatives of $\Phi$) stands for a quantity for which
\begin{equation*}
    |\mathfrak{R}_q(\Phi,\Phi',\Phi'')| \leq C \e^q (1+|z|^d) e^{- k |z|}
    + C \e^{q+1} (1+|z|^d) e^{- k |z|} |\Phi''|;
\end{equation*}
\begin{eqnarray*}
 \nonumber & & |\mathfrak{R}_q(\Phi,\Phi',\Phi'') -
  \mathfrak{R}_q(\tPhi,\tPhi',\tPhi'')| \leq C
    \e^q (1+|z|^d) [|\Phi-\tPhi| + |\Phi'-\tPhi'|] e^{- k |z|} \\
  & + & C \e^{q+1} (1+|z|^d) \left( |\Phi'' + \tPhi''| (|\Phi-\tPhi|+|\Phi'-\tPhi'|) + |\Phi'' -
\tPhi''| \right) e^{- k |z|}.
\end{eqnarray*}
Similarly, we will let $\mathfrak{R}_q(\ov{s})$ denote a quantity
(depending only on $\ov{s}$ and $z$) such that
\begin{equation*}
  |\mathfrak{R}_q(\ov{s})| \leq C \e^q (1+|z|^d) e^{- k |z|},
\end{equation*}
and which depends smoothly on $\ov{s}$. In the estimates below, the
assumptions \eqref{eq:bdc1c2}-\eqref{eq:bdPhif2d} will be used: one
hand by \eqref{eq:bdc1c2} we have $L^\infty$ estimates on $\Phi,
f_2$ and their first derivatives; one the other by
\eqref{eq:bdPhif2d} we have $L^2$ estimates on the higher order
derivatives, of the type $\|\Phi^{(l)}\|_{L^2} \leq C_l
\frac{\d^l}{\e^l} \|\Phi\|_{L^2}$, for $l \in \N$.

We will also use notations like $\Phi \mathfrak{R}_q(\Phi,\Phi')$,
$f_2'' \mathfrak{R}_q(\Phi,\Phi')$, etc., to denote error terms
which are products of functions of $\ov{s}$, like $\Phi$ or $f'_2$,
and the above $\mathfrak{R}_q$'s.

\

\no Having defined this notation, we can compute (and estimate)
$S_\e(\Tilde{\Psi}_{2,\e})$ term by term.

\bigskip

\noindent {$\bullet$ \bf Estimate of $\mathfrak{A}_1$}

\bigskip

\no From the expression of the Laplace-Beltrami operator (see
Subsection 3.1 in \cite{mmm1}) it follows that
\begin{eqnarray*}
  e^{i \frac{\e^2 f_2(\e s)}{\e}}
  S_\e(\tilde{\Psi}_{1,\e}) - S_\e(\tilde{\psi}_{1,\e}) & = &
  \tilde{g}^{11} \left[ \e^4 (f'_2)^2 \tilde{\psi}_{1,\e} +
  i \e^3 f''_2 \tilde{\psi}_{1,\e} + 2 i \e^2 f'_2
  \pa_s \tilde{\psi}_{1,\e} \right] \\
  & + & 2 i \sum_l \tilde{g}^{1l} \e^2 f'_2 \pa_l
  \tilde{\psi}_{1,\e} + \frac{i}{\sqrt{\det \tilde{g}}} \pa_A \left(
  g^{A1} \sqrt{\det \tilde{g}} \right) \e^2 f'_2
  \tilde{\psi}_{1,\e}.
\end{eqnarray*}
Using the expressions of $w_r$, $w_i$ and the expansions of the
metric coefficients in Subsection 3.1 of \cite{mmm1}, multiplying
the last equation by $e^{i \frac{\tilde{f}_0(\e s)}{\e}}$ one
obtains
\begin{eqnarray}\label{eq:estA1} \nonumber
e^{i \frac{\tilde{f}(\e s)}{\e}} \mathfrak{A}_1 & = & e^{i
\frac{\tilde{f}_0(\e s)}{\e}} \left( e^{i \frac{\e^2 f_2(\e
s)}{\e}}  S_\e(\tilde{\Psi}_{1,\e}) - S_\e(\tilde{\psi}_{1,\e})
\right) = e^{i \frac{\tilde{f}(\e s)}{\e}} S_\e(\tilde{\Psi}_{1,\e})
- e^{i \frac{\tilde{f}_0 (\e s)}{\e}} S_\e(\tilde{\psi}_{1,\e}) \\
& = & \mathfrak{A}_{1,0} + \tilde{\mathfrak{A}}_1 :=
\mathfrak{A}_{1,0}+\mathfrak{A}_{1,r,e}+\mathfrak{A}_{1,r,o}+
\mathfrak{A}_{1,i,e}+\mathfrak{A}_{1,i,o}+\mathfrak{A}_{1,1},
\end{eqnarray}
where
\begin{eqnarray*}
{\mathfrak{A}}_{1,0} =  2 \e^2 f' f'_2 h
U; \qquad
{\mathfrak{A}}_{1,r,e} =  \e^3 f'_2 \left[ 2 f'_1 h U + 4 \langle
{\bf H},\Phi \rangle f' h U + 2 f' w_{r,e} \right],
\end{eqnarray*}
\begin{eqnarray*}
{\mathfrak{A}}_{1,r,o}&=& \e^3 f'_2 \left[ 4 \langle {\bf H}, z \rangle f' h U + 2 f' w_{r,o}
\right];
\end{eqnarray*}
\begin{eqnarray}\label{eq:mfA1ie}
{\mathfrak{A}}_{1,i,e} =  i\e^3
f_2''hU + 2i\e^3f_2'\left[f' w_{i,e}+h'U+hk'\n U \cdot z \right]; \qquad
{\mathfrak{A}}_{1,i,o} = 2i\e^3f'f_2'w_{i,o}
\end{eqnarray}
\begin{eqnarray*}
{\mathfrak{A}}_{1,1}&=&  (f_2')^2\mathfrak{R}_4(\Phi,\Phi')
+ f_2''\mathfrak{R}_5(\Phi,\Phi') + f_2' \Phi'' \mathfrak{R}_4(\Phi,\Phi') +f_2'
\mathfrak{R}_4(\Phi,\Phi').
\end{eqnarray*}

\

\no {$\bullet$ \bf  Estimate of $\mathfrak{A}_2$}

\bigskip

\no Reasoning as for the previous estimate, collecting the terms of
order $\e^3$ and higher in $S_\e(\tilde{\psi}_{1,\e})$, we obtain
\begin{equation}\label{eq:estA2}
 e^{i \frac{\tilde{f}(\e s)}{\e}} \mathfrak{A}_2 =
\mathfrak{A}_{2,0} + \tilde{\mathfrak{A}}_2 :=
\mathfrak{A}_{2,r,e}+\mathfrak{A}_{2,r,o}+
\mathfrak{A}_{2,i,e}+\mathfrak{A}_{2,i,o}+\mathfrak{A}_{2,1},
\end{equation}
where $\mathfrak{A}_{2,0} = 0$ and where the remaining terms are
given by
\begin{equation*}
\mathfrak{A}_{2,r,e} =
\e^3\Phi''F_e(\ov s);
\end{equation*}
\begin{equation*}
\mathfrak{A}_{2,r,o} = 2\e^3hk\langle{\bf H},\Phi\rangle
\sum_j\Phi_j''\pa_jU+\e^3f'hf_1'\sum_j z_j\Phi''_jU
+2\e^3f'^2h\langle{\bf H},\Phi\rangle \sum_j\Phi''_j\,z_jU;
\end{equation*}
\begin{equation*}
\mathfrak{A}_{2,i,e} = -2i\e^3f'h\sum_{j,l}\Phi'_l\Phi''_j\,\pa_{l}(z_jU)
+ 2i\e^4f'h\langle{\bf H},z\rangle \sum_j\Phi'''_j\,z_jU;
\end{equation*}
\begin{equation*}
\mathfrak{A}_{2,i,o} = i\e^3\sum_j\Phi_j''z_j\left( f''hU
+f'h'U+f'hk'\n U\cdot z \right) +i\e^3\Phi'''F_o(\ov s) + 2i\e^4f'h\langle{\bf H},\Phi\rangle
\sum_j\Phi'''_j\,z_jU;
\end{equation*}
\begin{equation*}
\mathfrak{A}_{2,1} = \mathfrak{R}_3(\ov{s}) + (\Phi + \Phi')
\mathfrak{R}_3(\Phi,\Phi',\Phi'') + \mathfrak{R}_4(\Phi,\Phi',\Phi''),
\end{equation*}
where $F_e(\ov s)$ and $F_o(\ov s)$ are respectively an even real
function and an odd real function in the variables $z$, with smooth
coefficients in $\ov s = \e s$, and satisfying the decay property
$|F_e(\ov s)| + |F_o(\ov s)| \leq C (1+|z|^d) e^{- k |z|}$.

\

\no {$\bullet$ \bf  Choice of $\tilde{v}$ and estimate of
$\mathfrak{A}_3$}

\bigskip

\no We choose the function $\tilde{v}$ in such a way to annihilate
(roughly) one of the main terms in \eqref{eq:estA1}, namely
$2\e^2f'f'_2hU(kz)$. Hence we define $\tilde{v}$ so that it solves
\begin{equation}\label{eq:Lrtildev}
    \mathcal{L}_r\widetilde{v}=-2f'f'_2hU(kz).
\end{equation}
Reasoning as for the definition of $w_r$ (see \cite{mmm1},
Subsection 3.2), $\widetilde{v}$ can be  explicitly determined as
\begin{equation*}
    \widetilde{v}=2f'f'_2h\widetilde{U}(kz).
\end{equation*}
With this definition, using the above estimates on the metric
coefficients and the expressions of error terms, the linear terms
involving $E$ in $S_\e(\Tilde{\Psi}_{2,\e})$ can be written as
\begin{eqnarray}\label{eq:estA3}
& & e^{i \frac{\tilde{f}(\e s)}{\e}} \mathfrak{A}_3 =
\mathfrak{A}_{3,0} + \tilde{\mathfrak{A}}_3 :=  \mathfrak{A}_{3,0} +
\mathfrak{A}_{3,r,e}+\mathfrak{A}_{3,r,o}+
\mathfrak{A}_{3,i,e}+\mathfrak{A}_{3,i,o}+\mathfrak{A}_{3,1},
\end{eqnarray}
where
\begin{equation*}
    \mathfrak{A}_{3,0} = \e^2 \mathcal{L}_r \tilde{v};
\end{equation*}
\begin{eqnarray*}
 \mathfrak{A}_{3,r,e}&=& 2\e^3 f' h f'_2 \left( 2 (f')^2 \langle {\bf H}, \Phi \rangle
   + \langle \n^N V, \Phi \rangle - p (p-1)
   h^{p-2} U^{p-2} w_{r,e} \right) \tilde{U}\\
  &+&4\e^3 (f')^2 f'_1 h f'_2 \tilde{U}(kz)
  -2\e^4f'f_2'''h\widetilde{U}
\end{eqnarray*}
\begin{eqnarray*}
 \mathfrak{A}_{3,r,o}&=&2 \e^3f' h f'_2 \left( 2 (f')^2 \langle {\bf H}, z \rangle
   + \langle \n^N V,z \rangle - p (p-1)
   h^{p-2} U^{p-2} w_{r,o} \right) \tilde{U}\\
     &+& 2 \e^3f' f'_2 h k\sum_j H^j \pa_j
\tilde{U}(kz)-2\e^4f'hkf_2'\sum_j\Phi_j''\pa_j \widetilde{U}
\end{eqnarray*}
\begin{eqnarray}\label{eq:estA3ie}
 \mathfrak{A}_{3,i,e}&=&4 i f'\e^3 \pa_{\ov{s}} \left( h f' f'_2 \tilde{U}
    \right) + 2 i \e^3f' f'_2 \left( f'' h  - (p-1)
   h^{p-1} U^{p-1} w_{i,e}  \right) \tilde{U}
\end{eqnarray}
\begin{eqnarray*}
 \mathfrak{A}_{3,i,o}&=&- 2(p-1) i \e^3f' f'_2
   h^{p-1} U^{p-1} w_{i,o} \tilde{U}- 4 i \e^3(f')^2
   f'_2 hk \sum_j \Phi'_j \pa_j \tilde{U}
\end{eqnarray*}
\begin{eqnarray*}
 \mathfrak{A}_{3,1}&=&f_2'\left[\mathfrak{R}_4(\Phi,\Phi',\Phi'')\right]
+f_2''\left[\mathfrak{R}_4(\Phi,\Phi')+\mathfrak{R}_6(\Phi,\Phi',\Phi'')\right]
+\e^4f_2'''\left[\mathfrak{R}_1(\Phi,\Phi')\right] \nonumber
\\ & + & \mathfrak{R}_4(\Phi,\Phi')  f'_2 \left[ f'_2
(1 + \e^2 f'_2) + \e f''_2 \right] + \mathfrak{R}_5(\Phi,\Phi') f''_2
f'_2 + \mathfrak{R}_6(\Phi,\Phi',\Phi'') (f'_2)^2.
\end{eqnarray*}

\

\no {$\bullet$ \bf Choice of $v_0$ and estimate of $\mathfrak{A}_4$}

\bigskip

\no In order to make the approximate solution as accurate as
possible,  we add a correction $\e^2 v_0$ in such a way to
compensate (most of) the terms $\e^2 (\tilde{R}_{r,e} + i
\tilde{R}_{i,o})$, see Subsection 3.3 in \cite{mmm1}. We notice that
these terms contain parts which are independent of $\Phi$, which  we
denote them by $\tilde{R}_{r,e}^0$ and $\tilde{R}_{i,o}^0$
respectively, and parts which are quadratic in $\Phi$ or its
derivatives, $\tilde{R}_{r,e}^\Phi$ and $\tilde{R}_{i,o}^\Phi$
respectively. Since we will take $\Phi$ of order $\e$, we regard the
latter ones as higher order terms, and we add corrections to cancel
$\tilde{R}_{r,e}^0$ and $\tilde{R}_{i,o}^0$. Precisely we define
$v_{r,e}^0$ and $v_{i,o}^0$ by
\begin{eqnarray}\label{eq:Rre} \nonumber
    - \mathcal{L}_{r} v_{r,e}^0
    & = & - \tfrac 12 (f')^2 h U(kz)  \sum_{l,m} \partial^2_{lm}
    g_{11} z_m z_l   +2(f')^2 \langle {\bf H}, + w_{r,o}z\rangle
    \nonumber + 4(f')^2 h U(kz)
    \langle {\bf H},z\rangle^2 + 2 f' \partial_s w_{i,e}  \\ \nonumber & + & f'' w_{i,e}
    - \left[ h'' U(kz) + 2 h' k'
    \n U(kz) \cdot z + h k'' \n U(kz) \cdot z + h (k')^2 \n^2 U (kz) [z,z]
    \right] \\ & + & \tfrac 12
    \sum_{l,m} \partial^2_{lm}
    g_{tj} z_m z_l h k^2
    \partial^2_{tj} U(kz) +  k h [\langle {\bf H},z\rangle H^m - \tfrac 12
    \sum_l \partial^2_{ml} g_{11} z_l ] \partial_{m} U(kz)\\ & + &h k  \sum_l \partial^2_{lm} g_{mj} z_l
     \partial_{j} U(kz) + \sum_l H^l \partial_{l} w_{r,o}
    + h k \langle {\bf H},z\rangle H^l  \partial_{l} U(kz) \nonumber + \langle
    \n^N V,w_{r,o}z \rangle \\ & - & \tfrac 12 (p-1) h^{p-2}
    U(kz)^{p-2} w_{i,e}^2 -
    \nonumber \tfrac 12 p (p-1) h^{p-2} U(kz)^{p-2} w_{r,o}^2
      + \tfrac{1}{2} \sum_{m,j} \partial^2_{mj} V z_m z_j h U(kz); \nonumber
\end{eqnarray}
\begin{eqnarray}\label{eq:Rio} \nonumber
    - \mathcal{L}_{i} v_{i,o}^0 & = & 2 \left[ f'' h U(kz) + 2 f' h' U(kz) + 2 f' h k'
    \n U(kz) \cdot z \right] \langle {\bf H},z \rangle
     + \sum_i H^j\partial_{j} w_{i,e} \\
    & + & 2(f')^2\langle {\bf H},w_{i,e}z  \rangle
     + 2 f' \partial_s w_{r,o} + f'' w_{r,o} -  f' h k \sum_j\partial_{j} U(kz)
      \sum_{l,m} \partial^2_{lm} g_{1j} z_m z_l
     \\ & - & f' h U(kz) \left( \sum_m \partial^2_{1m} g_{11} z_m
    \right) -f'h \left( \sum_{j,l} \partial^2_{lj} g_{1j} z_l \right)  U(kz) + \frac 12f' h
\left( \sum_l \partial^2_{1l} g_{11} z_l
    \right)  U(kz) \nonumber \\ & - & (p-1) h^{p-2} U(kz)^{p-2}  w_{r,o} w_{i,e}
     + \langle \n^N V,w_{i,e}z \rangle. \nonumber
\end{eqnarray}
We notice that the right-hand side of \eqref{eq:Rre} is even in $z$,
and hence orthogonal to the kernel of $\mathcal{L}_r$. As a
consequence the equation is indeed solvable in $v_{r,e}^0$, see the
comments after \eqref{eq:LrLi}. The same comment applies to
\eqref{eq:Rio}, where the right-hand side which is odd in $z$.
Furthermore the right-hand sides decay at infinity at most like
$(1+|z|^d) e^{-k|z|}$ for some integer $d$, so the same holds true
for $v_{r,e}^0$ and $v_{i,o}^0$. In conclusion, after some
computations one finds
$$
  e^{i \frac{\tilde{f}(\e s)}{\e}} \mathfrak{A}_4 = \mathfrak{A}_{4,0} +
  \tilde{\mathfrak{A}}_4 :=  \mathfrak{A}_{4,0} +
\mathfrak{A}_{4,r,e}+\mathfrak{A}_{4,r,o}+
\mathfrak{A}_{4,i,e}+\mathfrak{A}_{4,i,o}+\mathfrak{A}_{4,1},
$$
where
\begin{equation*}
  \mathfrak{A}_{4,0} = \e^2\mathcal{L}_r v_{r,e}^0 + i \e^2 \mathcal{L}_i v_{i,o}^0;
\end{equation*}
\begin{equation}\label{eq:estA401}
\mathfrak{A}_{4,r,e} = \e^3 F_{4,r,e}(\ov{s}); \qquad \mathfrak{A}_{4,r,o} = \e^3
F_{4,r,o}(\ov{s}); \qquad \mathfrak{A}_{4,i,e} = \e^3 F_{4,i,e}(\ov{s}); \qquad
\mathfrak{A}_{4,i,o} = \e^3 F_{4,i,o}(\ov{s});
\end{equation}
\begin{equation}\label{eq:estA4}
\mathfrak{A}_{4,1} = \mathfrak{R}_4(\Phi,\Phi') + (\Phi + \Phi')  (1 + f'_2)
\mathfrak{R}_3(\Phi,\Phi') + f''_2 \mathfrak{R}_5(\Phi,\Phi')
+ (f'_2)^2 \mathfrak{R}_6(\ov{s}) + \Phi'' \mathfrak{R}_4(\Phi,\Phi').
\end{equation}
As for  $F_e(\ov s)$ and $F_o(\ov s)$ in $\mathfrak{A}_2$, the
$F_4$'s depend only on $V$, $\g$, $M$, and are bounded above by $C
(1+|z|^d) e^{- k |z|}$.

\

\no {$\bullet$ \bf Estimate of $\mathfrak{A}_5$}

\bigskip

\no The term involving $v_\d$ in $S_\e(\Tilde{\Psi}_{2,\e})$ is
given by

\begin{eqnarray}\label{eq:estA5}
\mathfrak{A}_5=\mathfrak{A}_{5,0}+\tilde{\mathfrak{A}}_5 :=
\mathfrak{A}_{5,0}+\mathfrak{A}_{5,r,e}+\mathfrak{A}_{5,r,o}
+\mathfrak{A}_{5,i,e}+\mathfrak{A}_{5,i,o}+\mathfrak{A}_{5,1}
\end{eqnarray}
where
\begin{equation}\label{eq:estA50}
\mathfrak{A}_{5,0} = \b \mathcal{L}_r Z_{\a}(k z) - \e^2 \b'' Z_{\a}(k z)
- 2 \e \xi' W_{\a} f' + i
\xi \mathcal{L}_i W_{\a} - i \e^2
\xi'' W_{\a} + 2 i \e \b' f' Z_{\a};
\end{equation}
\begin{eqnarray*}
\mathfrak{A}_{5,r,e} & = & - \e f'' \xi W_{\a}
- 2 \e^2 \b' \left( \frac{\pa Z_{\a}}{\pa {\a}} \a' + k'
\n Z_{\a}(k z) \cdot z \right) \\ & - & 2 \e f' \xi
\left( \frac{\pa W_{\a}}{\pa {\a}} \a' + k'
\n W_{\a}(k z) \cdot z \right) - (p-1) \e h^{p-2} U^{p-2}
\xi w_{i,e} W_{\a};
\end{eqnarray*}
\begin{eqnarray*}
\mathfrak{A}_{5,r,o}&=& \e \b \sum_j H^j\pa_j Z_{\a} +2 \e \langle{\bf
H},z\rangle \left[ (f')^2 \b Z_{\a} + \e^2 \b'' Z_{\a} - 2 \e
\xi' f' W_{\a} \right] \\ & + &
 \e \langle \n^N V,z \rangle \b Z_{\a} -p (p-1) \e
h^{p-2}U^{p-2} w_{r,o} \b Z_{\a};
\end{eqnarray*}
\begin{eqnarray*}
\mathfrak{A}_{5,i,e} & = & \e f'' \b Z_{\a}
- 2 \e^2 \xi' \left( \frac{\pa
W_{\a}}{\pa \a} \a' + k'
\n W_{\a}(k z) \cdot z \right) \\ & + & 2 \e f' \b
\left( \frac{\pa Z_{\a}}{\pa {\a}} \a' + k'
\n Z_{\a}(k z) \cdot z \right) - (p-1) \e h^{p-2} U^{p-2}
\b w_{i,e} Z_{\a};
\end{eqnarray*}
\begin{eqnarray*}
\mathfrak{A}_{5,i,o} & = & \e \xi \sum_j H^j\pa_j W_{\a}
+2 \e \langle{\bf H},z\rangle \left[ (f')^2 \xi W_{\a}
+ \e^2  \xi'' W_{\a} + 2 \e
\b' f' Z_{\a} \right] \\ & + &
 \e \langle \n^N V,z \rangle \xi W_{\a} - (p-1) \e
h^{p-2}U^{p-2} w_{r,o} \xi W_{\a}.
\end{eqnarray*}
The error term $\mathfrak{A}_{5,1}=\mathfrak{A}_{5,1}(\b,\Phi,f_2)$
satisfies the following estimates
$$
|\mathfrak{A}_{5,1}(\b,\Phi,f_2)|\le C(\e^2+\e^2|f_2'|+\e^3|f_2''|+\e^3|\Phi''|)
(1+|z|^d) e^{- k |z|} (|\b| + \e |\b'| + \e^2 |\b''| + \e^3 |\b'''|);
$$
\begin{eqnarray*}
\left|  \mathfrak{A}_{5,1}(\b,\Phi,f_2)-
\mathfrak{A}_{5,1}(\tilde{\b},\tilde{\Phi},\tilde{f}_2) \right|
&\le& C(\e |\Phi-\tilde{\Phi}|+\e|\Phi'-\tilde{\Phi}'|+\e^3
|\Phi''-\tilde{\Phi}''|+\e^2|f_2'-\tilde{f}_2'|+\e^3|f_2''-\tilde{f}_2''|) \\
& \times & (1+|z|^d) e^{- k |z|} (|\b| + \e |\b'| + \e^2 |\b''| + \e^3 |\b'''|)
\\ & + & C(\e^2+\e^2|f_2'|+\e^3|f_2''|+\e^3|\Phi''|) \\ & \times &
(1+|z|^d) e^{- k |z|} (|\b - \tilde{\b}| + \e |\b' - \tilde{\b}'| + \e^2
|\b'' - \tilde{\b}''| + \e^3 |\b''' - \tilde{\b}'''|).
\end{eqnarray*}
By the form of the function $\b$, see \eqref{eq:asynul},
\eqref{eq:bbj} and \eqref{eq:exprbeta}, its Fourier modes are
naively concentrated around indices of order $\frac{1}{\e}$. As a
consequence, $L^2$ norms of functions like $\e \b, \e^2 \b'', \e^3
\b'''$, etc. can be controlled with the $L^2$ norm of $\b$, see also
the comments before \eqref{eq:hatws}.

\

\no $\bullet$ {\bf Estimate of $\mathfrak{A}_6$}

\bigskip

\noindent First of all we notice that we are taking $\Phi'$ and
$f'_2$ in $H^1([0,L])$, and hence they belong to $L^\infty([0,L])$.
As a consequence, since we have the bound  $\|\b\|_{L^\infty([0,L])}
+ \e \|\b'\|_{L^\infty([0,L])} + \e^2 \|\b''\|_{L^\infty([0,L])} +
\e^3 \|\b'''\|_{L^\infty([0,L])} \leq C \e^2$ (which follows from
\eqref{eq:bdbe3} and the above comments), one has the estimate
\begin{equation}\label{eq:EFG}
    |E| + |F| + |G| \leq C \e^2 (1+|z|^d) e^{-k |z|}.
\end{equation}
If then we choose $\d$ sufficiently small (recall also the
expressions of $w_r$, $w_i$ and \eqref{eq:psi2edec}), we deduce that
$$
  |\tilde{\Psi}_{2,\e} - \tilde{\Psi}_{1,\e}| \leq |\tilde{\Psi}_{1,\e}| \qquad
  \quad \hbox{ in } \tilde{D}_\e.
$$
This estimate implies that $\mathfrak{A}_6$ admits a uniform
quadratic Taylor expansion in $|\tilde{\Psi}_{2,1} -
\tilde{\Psi}_{1,\e}|$ and is bounded by $|\tilde{\Psi}_{1,\e}|^{p-2}
|\tilde{\Psi}_{2,\e} - \tilde{\Psi}_{1,\e}|^2$. Precisely, we can
write
\begin{equation}\label{eq:estA6}
    \mathfrak{A}_6 = \mathfrak{A}_{6,0} + \tilde{\mathfrak{A}}_6 :=
\mathfrak{A}_{6,0}+\mathfrak{A}_{6,r,e}+\mathfrak{A}_{6,r,o}
+\mathfrak{A}_{6,i,e}+\mathfrak{A}_{6,i,o}+\mathfrak{A}_{6,1},
\end{equation}
where
\begin{equation}\label{eq:estA62}
  \mathfrak{A}_{6,0} =
\mathfrak{A}_{6,r,e} = \mathfrak{A}_{6,r,o} = \mathfrak{A}_{6,i,e} =
\mathfrak{A}_{6,i,o} = 0; \qquad \quad \mathfrak{A}_{6,1} =
R_4(f'_2,\Phi,\Phi',\b),
\end{equation}
where $R_4(f'_2,\Phi,\Phi',\b)$ is a quantity satisfying the
estimates
$$
  |R_4(f'_2,\Phi,\Phi',\b)| \leq C \left[ \e^4 + (\e^2 + \|\b\|_{L^\infty} +
  \e \|\b'\|_{L^\infty}) (|\b| + \e |\b'|) \right] (1+|z|^d) e^{-k|z|};
$$
\begin{eqnarray*}
  |R_4(f'_2,\Phi,\Phi',\b) - R_4(\tilde{f}'_2,\tilde{\Phi},
  \tilde{\Phi'},\tilde{\b})| & \leq & C  \left( \e^2+ |\b| + \e |\b'|
   + |\tilde{\b}| + \e |\tilde{\b}'|   \right) (1+|z|^d) e^{-k|z|} \\
     & \times & \left(\e |\Phi-\tilde{\Phi}|+ \e |\Phi'-\tilde{\Phi}'|
    +\e^2 |f'_2-\tilde{f}'_2|+|\b-\tilde{\b}| + \e |\b'-\tilde{\b}'| \right).
\end{eqnarray*}

\no $\bullet$ {\bf Final estimate of $S_\e(\Tilde{\Psi}_{2,\e})$}

\

\noindent By \eqref{eq:ssss}, in the above notation we have
\begin{eqnarray*}
  e^{i \frac{\tilde{f}(\e s)}{\e}} S_\e(\tilde{\Psi}_{2,\e}) &=&
  \e^2(\tilde{R}_{r,o}+\tilde{R}_{r,e}) +
  \e^2(\tilde{R}_{r,o,f_1}
 +\tilde{R}_{r,e,f_1}) \\
   & + & \e^2i(\tilde{R}_{i,e}
 +\tilde{R}_{i,o})
+\e^2i(\tilde{R}_{i,e,f_1}+\tilde{R}_{i,o,f_1})
 + \sum_{i=1}^6 \mathfrak{A}_{i,0} + \sum_{i=1}^6
 \tilde{\mathfrak{A}}_i.
\end{eqnarray*}
Recalling the choices of $\tilde{v}$, $v_{r,e}^0$ and $v_{i,o}^0$ in
\eqref{eq:Lrtildev} and \eqref{eq:Rre}, \eqref{eq:Rio} (and
recalling the notation for the $R$'s after \eqref{eq:esSe}) we
finally obtain the following result.

\begin{pro}\label{p:errest}
Suppose $\Phi$, $f_2$ and $\b$ satisfy \eqref{eq:bdc1c2},
\eqref{eq:bdPhif2d} and \eqref{eq:bdbe3} for some $c_1, c_2, c_3 >
0$.  Let $\tilde{f} = f + \e f_1  + \e^2 f_2$, where $f$ is given in
\eqref{eq:f'Cintr} and $f_1$ in \eqref{eq:f1}. Let also $w_{r} =
w_{r,e} + w_{r,o}$, with $w_{r,e}, w_{r,o}$ given respectively in
\eqref{eq:we}, \eqref{eq:wro}, and $w_i = w_{i,e} + w_{i,o}$, where
$w_{i,e}$ and $w_{i,o}$ are given respectively in \eqref{eq:wiewio}.
Let  $\tilde{\Psi}_{2,\e}$ be defined in \eqref{eq:P2eeee}. Then, as
$\e$ tends to zero, we have that
\begin{eqnarray}\label{eq:finalfinal} \nonumber
 e^{i \frac{\tilde{f}}{\e}} S_\e(\Tilde{\Psi}_{2,\e}) & = &
 \e^2 (\tilde{R}_{r,e}^\Phi +\tilde{R}_{r,o}+\tilde{R}_{r,e,f_1}+
 \tilde{R}_{r,o,f_1}) + \e^2i(\tilde{R}_{i,e}+\tilde{R}_{i,o}^\Phi
 +\tilde{R}_{i,e,f_1}+\tilde{R}_{i,o,f_1})
 \\
 &+& \b \mathcal{L}_r Z_{\a}(k z) + \e^2 \b'' Z_{\a}(k z)
- 2 \e \xi' W_{\a} f' +i
\xi \mathcal{L}_i W_{\a}  \\ & + & i \e^2
\xi'' W_{\a} + 2 i \e \b' f'
Z_{\a} + e^{i \frac{\tilde{f}}{\e}} \sum_{j=1}^6 \tilde{\mathfrak{A}}_j,
\nonumber
\end{eqnarray}
where the $R$'s are as in \eqref{eq:esSe}, where
$\tilde{R}_{r,e}^\Phi$, $\tilde{R}_{i,o}^\Phi$ are the terms
quadratic in $\Phi, \Phi'$ within $\tilde{R}_{r,e}$,
$\tilde{R}_{i,o}$, and where the latter error terms are given in
\eqref{eq:estA1}, \eqref{eq:estA2}, \eqref{eq:estA3},
\eqref{eq:estA4}, \eqref{eq:estA5}, \eqref{eq:estA6} respectively.
\end{pro}

\section{Proof of Theorem \ref{t:main}}

In this section we prove our main theorem. First we solve the
equation in the $\ov{H}_\e$ components, see \eqref{eq:orth}, using a
Lyapunov-Schmidt reduction. Then we turn to the components in
$\tilde{K}_\d$ and solve the bifurcation equation as well: in this
last step we use crucially the non-degeneracy assumption on $\g$ and
an accurate choice for the values of the parameter $\e$.

\subsection{Solvability in the component of $\ov{H}_\e$}\label{ss:hatpsi}

In Proposition \ref{p:sec2final} we showed that  problem
\eqref{eq:new} is reduced to finding a solution of $L_\e(\phi) =
\tilde{S}_\e(\phi)$ in $\tilde{D}_\e$, see \eqref{eq:defLe},
\eqref{eq:LetildeSe} and \eqref{eq:tildeSe}, if we take $K^2(\e s) =
V(\e s)$. Choosing in Proposition \ref{p:truefndec} as approximate
solution $\tilde{\psi}_\e = \tilde{\Psi}_{2,\e}$ (the function
constructed in the previous subsection), we have the following
result where, as usual, $\d$ is sufficiently small. We recall
Proposition \ref{p:sec2final}, formulas
\eqref{eq:l2cmas}-\eqref{eq:normstar} and the definition of
$\tilde{K}_\d$ after \eqref{eq:hkcmas}: also, we denote by
$\tilde{\Pi}_\e$ the orthogonal projection onto the set $\{ e^{- i
\frac{\tilde{f}(\e s)}{\e}} \tilde{v} \; : \; \tilde{v} \in
\tilde{K}_\d \}$.

\begin{pro}\label{p:exhatphi}
Let $\tilde{\Psi}_{2,\e}$ be as in Proposition \ref{p:errest}. Then
there exists $\check{v}_\d \in \tilde{K}_\d$, depending on the
parameters $\Phi, f_2, \b$, such that the following problem admits a
solution
\begin{equation}\label{eq:solvortfinnn}
  \left\{
    \begin{array}{ll}
      - \D_{g_\e} \hat{\phi} + V(\e x) \hat{\phi} - |\tilde{\Psi}_{2,\e}|^{p-1}
      \hat{\phi} - (p-1) |\tilde{\Psi}_{2,\e}|^{p-3} \tilde{\Psi}_{2,\e} \Re
   (\tilde{\Psi}_{2,\e} \ov{\hat{\phi}}) =
   \tilde{S}_\e(\hat{\phi}) + e^{- i \frac{\tilde{f}(\e s)}{\e}}
   \check{v}; \\[2mm]  \hat{\phi} \in \ov{H}_\e, \quad \check{v} \in \tilde{K}_\d. & \end{array}
  \right.
\end{equation}
Furthermore, if $m \in N$, if $\tilde{\tilde{\Psi}}_{2,\e}$ is an
approximate solution corresponding to different $\Phi, f_2, \b$, for
a fixed constant $C$ independent of $\e$ and $\d$,  for $\t = \frac
12$ and $0 < \varsigma' < \varsigma < 1$ sufficiently small, we have
\begin{equation}\label{eq:solvortestfin3}
\|\hat{\phi}\|_{\varsigma',V} \leq \frac{C}{\d^2}
  \|\tilde{\Pi}_\e S_\e(\tilde{\Psi}_{2,\e}) \|_{L^2(C^{\t}_{\varsigma,V})} + C \e^m;
   \qquad \|\check{v}\|_{L^2(C^{\t}_{\varsigma,V})} \leq C
    \|S_\e(\tilde{\Psi}_{2,\e})\|_{L^2(C^{\t}_{\varsigma,V})}.
\end{equation}
\begin{equation}\label{eq:solvortestfin3'}
\|\hat{\phi} - \hat{\tilde{\phi}}\|_{\varsigma',V} \leq \frac{C}{\d^2}
\|\tilde{\Pi}_\e (S_\e(\tilde{\Psi}_{2,\e}) -
S_\e(\tilde{\tilde{\Psi}}_{2,\e})) \|_{L^2(C^{\t}_{\varsigma,V})}.
\end{equation}
\end{pro}

\begin{pf} The proof relies on Proposition \ref{p:exvarphi},
Proposition \ref{p:truefndec} and the contraction mapping theorem.
By Proposition \ref{p:truefndec}, the operator $L_\e$ (see
\eqref{eq:defLe}) is invertible from $(\ov{H}_\e,
\|\cdot\|_{\varsigma',V})$ into $L^2(C^{\t}_{\varsigma,V})$, and the
norm of the inverse is  uniformly bounded by $C/\d^2$. By this
invertibility, \eqref{eq:solvortfinnn} is satisfied if and only if
$\hat{\phi}$ is a fixed point of the operator $\check{F}_\e :
(\ov{H}_\e, \|\cdot\|_{\varsigma',V}) \to (\ov{H}_\e,
\|\cdot\|_{\varsigma',V})$ defined by
\begin{eqnarray*}
  \check{F}_\e(\hat{\phi}) & = & L_\e^{-1} \left[ \tilde{\Pi}_\e \left(
  \tilde{S}_\e(\hat{\phi}) \right) \right] := L_\e^{-1} \left[
  \tilde{\Pi}_\e \left( S_\e(\tilde{\Psi}_{2,\e})
    +  N_\e(\eta_\e \hat{\phi} + \varphi(\hat{\phi}))  \right. \right.
  \\ & + & \left. \left. |\tilde{\Psi}_{2,\e}|^{p-1} \varphi(\hat{\phi})
  + (p-1) |\tilde{\Psi}_{2,\e}|^{p-3} \tilde{\Psi}_{2,\e} \Re
  (\tilde{\Psi}_{2,\e} \ov{\varphi(\hat{\phi})}) \right) \right].
\end{eqnarray*}
We recall that, in the last formula, $\varphi(\hat{\phi})$ is given
by Proposition \ref{p:exvarphi}, while $N_\e$ is defined in
\eqref{eq:defNe}.

Our next goal is to show that $\check{F}_\e$ is a contraction on a
metric ball (in the $\|\cdot\|_{\varsigma',V}$ norm) of radius
$\frac{C}{\d^2} \|\tilde{\Pi}_\e S_\e(\tilde{\Psi}_{2,\e})
\|_{L^2(C^{\t}_{\varsigma,V})} + C \e^m$ for $C$ large enough and
$m$ arbitrary integer. Setting for simplicity
$$
  \check{G}_\e(\hat{\phi}) = N_\e(\eta_\e \hat{\phi} + \varphi(\hat{\phi})) +
  |\tilde{\Psi}_{2,\e}|^{p-1} \varphi(\hat{\phi})
  + (p-1) |\tilde{\Psi}_{2,\e}|^{p-3} \tilde{\Psi}_{2,\e} \Re
  (\tilde{\Psi}_{2,\e} \ov{\varphi(\hat{\phi})}),
$$
one clearly finds
\begin{equation}\label{eq:estcheckF}
    \left\{
    \begin{array}{ll}
      \|\check{F}_\e(\hat{\phi})\|_{L^2(C^{\t}_{\varsigma',V})} \leq \frac{C}{\d^2}
  \left( \| \tilde{\Pi}_\e S_\e(\tilde{\Psi}_{2,\e})\|_{L^2(C^{\t}_{\varsigma,V})}
  + \|\check{G}_\e(\hat{\phi})\|_{L^2(C^{\t}_{\varsigma,V})} \right); &  \\
      \|\check{F}_\e(\hat{\phi}_1) - \check{F}_\e(\hat{\phi}_2)
    \|_{L^2(C^{\t}_{\varsigma',V})} \leq \frac{C}{\d^2}
   \|\check{G}_\e(\hat{\phi}_1) -
\check{G}_\e(\hat{\phi}_2) \|_{L^2(C^{\t}_{\varsigma,V})}. &
    \end{array}
  \right.
\end{equation}
We next evaluate
$\|\check{G}_\e(\hat{\phi})\|_{L^2(C^{\t}_{\varsigma,V})}$, and show
that it is {\em superlinear} in
$\|\hat{\phi}\|_{L^2(C^{\t}_{\varsigma',V})}$ up to negligible
terms: we make first the following claim.

\

\noindent {\bf Claim}: in the notation \eqref{eq:wnmak1}, letting
$k_1(\ov{s}) = (\varsigma')^2 \sqrt{V(\ov{s})}$ we have
$\|\hat{\phi}\|_{C^{1,\frac 12}_{k_1}} \leq C
\|\hat{\phi}\|_{\varsigma',V}$ for some $C > 0$.

\

\noindent Assuming the claim true and choosing $\varsigma'' <
(\varsigma')^2$, we can apply Proposition \ref{p:exvarphi} with $\t
= \frac 12$, $k_0(\ov{s}) = \varsigma \sqrt{V(\ov{s})}$,
$k_1(\ov{s}) = (\varsigma')^2 \sqrt{V(\ov{s})}$ and $k_2(\ov{s}) =
\varsigma'' \sqrt{V(\ov{s})}$, to find
\begin{equation}\label{eq:qqqqq11}
\|\varphi(\hat{\phi})\|_{C^{\frac 12}_{-k_2}} \leq C \left( e^{-
\inf \frac{k_2 + k_0}{K} \e^{-\ov{\d}}} \|S_\e(\tilde{\Psi}_{2,\e})\|_{C^{\frac 12}_{k_0}}
    + e^{- \inf \frac{k_2 + k_1}{K} \e^{-\ov{\d}}} \|\hat{\phi}\|_{C^{1,\frac
    12}_{k_1}} \right).
\end{equation}
From the expression of $w_r, w_i, \tilde{v}, v_0$ and formula
\eqref{eq:EFG}, one can deduce that $|\tilde{\Psi}_{2,\e}| \leq C
e^{- k_0 |z|}$: moreover, from the estimates in the proof of
Proposition \ref{p:errest} one also finds that $\|
S_\e(\tilde{\Psi}_{2,\e})\|_{L^2(C^{\t}_{\varsigma,V})} \to 0$ as
$\e \to 0$. By \eqref{eq:Ne1} (recall that $\zeta > 0$), the latter
bounds on $\tilde{\Psi}_{2,\e}$, the previous claim and
\eqref{eq:qqqqq11}, if $m$ is an arbitrary integer and if
$\varsigma''$ is sufficiently close to $1$ after some elementary
computations we deduce
$$
  \|\check{G}_\e(\hat{\phi})\|_{L^2(C^{\t}_{\varsigma,V})}
     \leq C \left( \|\hat{\phi}\|_{L^2(C^{\t}_{\varsigma,V})}^{1 + \zeta}
   + \|\hat{\phi}\|_{L^2(C^{\t}_{\varsigma,V})}^{p} + \e^m \left( 1 +
  \|\hat{\phi}\|_{L^2(C^{\t}_{\varsigma,V})} \right) \right).
$$
Similarly, if $\|\hat{\phi}_1\|_{L^2(C^{\t}_{\varsigma,V})}$,
$\|\hat{\phi}_2\|_{L^2(C^{\t}_{\varsigma,V})}$ are finite one also
finds
$$
   \|\check{G}_\e(\hat{\phi}_1) - \check{G}_\e(\hat{\phi}_2)
    \|_{L^2(C^{\t}_{\varsigma,V})} \leq C \left[ \max_{l=1,2} \left\{
     \|\hat{\phi}_l\|_{L^2(C^{\t}_{\varsigma,V})}^{\zeta \wedge (p-1)} \right\} +
   \e^m \right]
   \|\hat{\phi}_1 - \hat{\phi}_2\|_{L^2(C^{\t}_{\varsigma,V})},
$$
where the symbol $\wedge$ stands for the minimum. Formula
\eqref{eq:estcheckF} and the latter one show that $\check{F}_\e$ is
a contraction, and we obtain \eqref{eq:solvortestfin3};
\eqref{eq:solvortestfin3'} follows similarly.

\

\noindent {\bf Proof of the claim.} According to our previous
notation, the norm $\| \cdot \|_{\varsigma',V}$ is evaluated using
the variables $(s,z)$, where the $z$'s are defined in
\eqref{eq:defzz}. If we want to estimate the $\| \cdot
\|_{C^{1,\frac 12}_{k_1}}$ norm instead, we should use
lipschitzianity with respect to $s$ and $y$.

Given $s_1, s_2 \in \R$ and $y_1, y_2 \in \R^{n-1}$ we want to
consider the difference $\n \hat{\phi}(s_1, y_1) - \n
\hat{\phi}(s_2, y_2)$. Recalling \eqref{eq:defzz} we can write that
\begin{eqnarray*}
  \pa_s \hat{\phi}(s_1, y_1) - \pa_s \hat{\phi}(s_2, y_2) & = & \pa_s
  \hat{\phi}(s_1, z_1 + \Phi(\e s_1)) - \pa_s \hat{\phi}(s_2, z_1 + \Phi(\e s_1)) \\
   & + & \pa_s \hat{\phi}(s_2, z_1 + \Phi(\e s_1)) - \pa_s \hat{\phi}(s_2, z_2 + \Phi(\e s_2)).
\end{eqnarray*}
By the definition of $\| \cdot \|_{\varsigma',V}$,  $\pa_s
\hat{\phi} \in H^1(C^{\t}_{\varsigma',V}) \subseteq C^{\frac
12}(C^{\t}_{\varsigma',V})$. This fact, the smoothness of
$V(\ov{s})$ and $\|\Phi\|_\infty + \|\Phi'\|_\infty \leq C(c_1) \e$
(which follows from \eqref{eq:bdc1c2}) imply that if $(s_1,y_1),
(s_2,y_2) \in B_1(s,y)$ then
$$
  e^{(\varsigma')^2 \sqrt{V(\ov{s})} |z|} |\pa_s \hat{\phi}(s_1, y_1) - \pa_s
  \hat{\phi}(s_2, y_2)| \leq C(c_1) \|\hat{\phi}\|_{\varsigma',V} \left(
  |s_1 - s_2|^{\frac 12} + |z_1 - z_2|^{\frac 12} + \e |s_1 - s_2| \right).
$$
A similar estimate holds for the derivatives of $\hat{\phi}$ with
respect to $y$, so from \eqref{eq:wnmak1} we get the conclusion.
\end{pf}

\

\noindent To apply Proposition \ref{p:exhatphi} we establish
explicit estimates on $\tilde{\Pi}_\e S_\e(\tilde{\Psi}_{2,\e})$ and
$\tilde{\Pi}_\e (S_\e(\tilde{\Psi}_{2,\e}) -
S_\e(\tilde{\tilde{\Psi}}_{2,\e}))$. Precisely, assuming from now on
$\t = \frac 12$, we have the following result.

\begin{pro}\label{p:esthatphi}
Assume $\Phi, f_2, \b, \tilde{\Phi}, \tilde{f}_2, \tilde{\b}$
satisfy conditions \eqref{eq:bdc1c2}, \eqref{eq:bdPhif2d} and
\eqref{eq:bdbe3}. Then, if $\hat{\phi}$ is defined as in Proposition
\ref{p:exhatphi} we have the estimates
\begin{equation}\label{eq:esthatphi}
 \sqrt{\e} \d^2 \|\hat{\phi}(\b,\Phi,f_2)\|_{*} \leq  C(c_1, c_2, c_3) \e^3;
\end{equation}
\begin{equation}\label{eq:esthpsilip}
  \sqrt{\e} \d^2 \|\hat{\phi}(\b,\Phi,f_2)-\hat{\phi}(\tilde{\b},\tilde{\Phi},\tilde{f}_2)\|_{*}
  \le  C(c_1, c_2, c_3) \left[  \e^2 \|\Phi-\tilde{\Phi}\|_{H^2} + \e^3
  \|f_2-\tilde{f}_2\|_{H^2}+ \e \|\b - \tilde{\b}\|_\sharp \right],
\end{equation}
where $C(c_1, c_2, c_3)$ is a positive constant depending on $c_1,
c_2, c_3$ but independent of $\e$ and $\d$.
\end{pro}

\begin{pf}
We prove \eqref{eq:esthatphi} only: \eqref{eq:esthpsilip} will
follow from similar considerations. To show \eqref{eq:esthatphi} we
use Proposition \ref{p:exhatphi}, so we are reduced to estimate
$\|\tilde{\Pi}_\e S_\e(\tilde{\Psi}_{2,\e})
\|_{L^2(C^{\t}_{\varsigma,V})} $, for which we can employ
\eqref{eq:finalfinal}.

By our assumptions on $\Phi, f_2, \b$ and by the estimates of the
previous subsection, it is easy to see that
$$
  \left\|
 \e^2 (\tilde{R}_{r,e}^\Phi +\tilde{R}_{r,o}+\tilde{R}_{r,e,f_1}+
 \tilde{R}_{r,o,f_1}) + \e^2i(\tilde{R}_{i,e}+\tilde{R}_{i,o}^\Phi
 +\tilde{R}_{i,e,f_1}+\tilde{R}_{i,o,f_1})
\right\|_{L^2(C^{\t}_{\varsigma,V})} \leq \frac{C(c_1, c_2, c_3) \e^3}{\sqrt{\e}};
$$
$$
\big\| e^{i \frac{\tilde{f}}{\e}} \sum_{i=1}^6 \tilde{\mathfrak{A}}_i
\big\|_{L^2(C^{\t}_{\varsigma,V})} \leq \frac{C(c_1, c_2, c_3) \e^3}{\sqrt{\e}}.
$$
Recall that in the choice of approximate solutions we have formally
corrected all the terms of order up to $\e^2$, so we are left with
terms of order $\e^3$ and higher. The factor $\sqrt{\e}$ in the
denominator arises from the fact that the length of $\g_\e$ is
$L/\e$: this gives a factor $\frac 1 \e$ when computing the $L^2$
norm squared, and we need then to take the square root. For the
estimates in $\tilde{\mathfrak{A}}_6$, which also require the
$L^\infty$ norm of $\b$, we can use the interpolation inequalities
$$
\|\b\|_{L^\infty([0,L])}\le C\|\b\|_{L^2([0,L])}^\frac12\|\b'\|_{L^2([0,L])}^\frac12\le
C\e^{\frac32}; \quad \quad \|\b'\|_{L^\infty([0,L])}\le
C\|\b'\|_{L^2([0,L])}^\frac12\|\b''\|_{L^2([0,L])}^\frac12 \le C\e^{\frac12}.
$$
It remains to consider now the other terms in the right-hand side of
\eqref{eq:finalfinal}, involving the functions $Z_\a$ and $W_\a$.
Let us call $\tilde{L}_\e^1$ the operator obtained from $L_\e^1$
(see \eqref{eq:tildeLep}) by replacing the variables $y$ with $z$
and $f$ with $\tilde{f}$. Let us first notice that the terms under
interest, with this notation, are nothing but $\tilde{\Pi}_\e
\tilde{L}^1_\e v_\d$.

Let us now recall the expression of $\b$ in \eqref{eq:exprbeta} and
$v_\d$ in \eqref{eq:vD}: if $\tilde{v}_{3,j}$ stand for the
functions in $K_{3,\d}$ (see \eqref{eq:K3d}) replacing $y$ with $z$,
we define the function
$$
\tilde{v}_\d = \sum_{j = - \frac{\d^2}{\e}}^{\frac{\d^2}{\e}}
b_j \tilde{v}_{3,j}.
$$
From the expression of $\tilde{v}_{3,j}$, see \eqref{eq:vv33dd}, one
finds that
\begin{equation}\label{eq:vdtvd}
    \|v_\d - \tilde{v}_\d\|_{\varsigma,V} \leq \frac{C}{\sqrt{\e}} \big( \sum_{j = -
\frac{\d^2}{\e}}^{\frac{\d^2}{\e}} b_j^2 \e^2 (1+j^2) \big)^{\frac 12};
\end{equation}
\begin{equation}\label{eq:oooo}
\tilde{L}_\e^1 (e^{- i \frac{\tilde{f}(\e s)}{\e}} v_\d) =
\tilde{L}_\e^1 (e^{- i \frac{\tilde{f}(\e s)}{\e}} \tilde{v}_\d) + O
\left( \frac{1}{\sqrt{\e}} \right) \big( \sum_{j = -
\frac{\d^2}{\e}}^{\frac{\d^2}{\e}} b_j^2 \e^2 (1+j^2)
\big)^{\frac 12} \quad (\hbox{in the } \|\cdot\|_{L^2(C^{\t}_{\varsigma,V})} \hbox{ norm}).
\end{equation}
Similarly to \eqref{eq:expa22}, recalling the asymptotic of $\nu_j$
(see \eqref{eq:asynul} and the lines before) one finds that
\begin{equation}\label{eq:aaaa}
    \tilde{L}_\e^1 (e^{- i \frac{\tilde{f}(\e s)}{\e}} \tilde{v}_\d) = e^{- i
\frac{\tilde{f}(\e s)}{\e}} \sum_{j=-\frac{\d^2}{\e}}^{\frac{\d^2}{\e}} \nu_j
b_j \tilde{v}_{3,j}+ R_1,
\end{equation}
where $\|R_1\|_{L^2(C^{\t}_{\varsigma,V})} \leq \frac{C}{\sqrt{\e}}
\big( \sum_{j = - \frac{\d^2}{\e}}^{\frac{\d^2}{\e}} b_j^2 \e^2
(1+j^2) \big)^{\frac 12} \leq C \sqrt{\e} \|\b\|_\sharp$. This
implies the conclusion, by \eqref{eq:bdbe3}.
\end{pf}

\subsection{Projections onto $\tilde{K}_\d$}\label{ss:proj}

In this section we estimate the projections of the equation onto the
components of $\tilde{K}_\d$. We estimate first  their size and
their Lipschitz dependence in the data $\Phi, f_2$ and $\b$. Then we
use the contraction mapping theorem to annihilate the function
$\check{v}$ in Proposition \ref{p:exhatphi}, which implies the
solvability of \eqref{eq:new}.

\subsubsection{Projection onto $\tilde{K}_{1,\d}$}\label{sss:pr1}

We want to evaluate the $\tilde{K}_{1,\d}$ component of the function
$\check{v}_\d$ in \eqref{eq:solvortfinnn}. To do this we consider a
normal section $\underline{\Phi}$ to $\g$ which satisfies the first
relation in \eqref{eq:bdPhif2d}, and the function
$$
  v_{\underline{\Phi}} := h(\e s)^{\frac{p+1}{4}} \left(
 \langle \underline{\Phi}(\e s), \n_z U(k z) \rangle + i \e
  \langle \underline{\Phi}'(\e s), z \rangle \frac{f'}{k} U(k z) - \frac{\e^2}{k^2}
  \langle \underline{\Phi}''(\e s), \mathfrak{V}(k z) \rangle \right).
$$
We then multiply both the left-hand side of \eqref{eq:solvortfinnn}
and $\tilde{S}_\e(\hat{\phi})$ (see \eqref{eq:tildeSe}) by the
conjugate of $e^{- i \frac{\tilde{f}(\e s)}{\e}}
v_{\underline{\Phi}}$, integrate over $\tilde{D}_\e$ and take the
real part. When multiplying the left-hand side, we can integrate by
parts and let the operator $L_\e$ act on $e^{- i \frac{\tilde{f}(\e
s)}{\e}} v_{\underline{\Phi}}$: using the arguments in the proofs of
Proposition \ref{p:inv} (see in particular \eqref{eq:expa} and
\eqref{eq:estRv}) and of Proposition \ref{p:esthatphi} one finds
that
$$
  L_\e (e^{- i \frac{\tilde{f}(\e s)}{\e}} v_{\underline{\Phi}}) =
  e^{- i \frac{\tilde{f}(\e s)}{\e}} \tilde{v}_{\underline{\Phi}} +
  R(v_{\underline{\Phi}}),
$$
where $\tilde{v}_{\underline{\Phi}} \in \tilde{K}_{1,\d}$, and where
$\|R(v_{\underline{\Phi}})\|_{L^2(C^{\t}_{\varsigma,V})} \leq C(\e +
\d^3) \|v_{\underline{\Phi}}\|_{L^2(C^{\t}_{\varsigma,V})} \leq
\frac{C}{\sqrt{\e}} (\e + \d^3) \|\underline{\Phi}\|_{L^2([0,L])}$.
Therefore, since $\hat{\phi}$ is orthogonal to $\tilde{K}_\d$, from
\eqref{eq:esthatphi} we deduce that
\begin{equation}\label{eq:small}
    \left| \Re \int_{\tilde{D}_\e} e^{i \frac{\tilde{f}(\e s)}{\e}}
    \overline{v_{\underline{\Phi}}} L_\e \hat{\phi} \right| dV_{\tilde{g}_\e}
 \leq \frac{C}{\sqrt{\e}} (\e + \d^3) \|\underline{\Phi}\|_{L^2([0,L])}
 \|\hat{\phi}\|_{L^2(C^{\t}_{\varsigma,V})} \leq C(c_1,c_2,c_3) \d \e^2
  \|\underline{\Phi}\|_{L^2([0,L])}.
\end{equation}
We next have to consider $\tilde{S}_\e(\hat{\phi})$, whose main term
is $S_\e(\Tilde{\Psi}_{2,\e})$: for this we use formula
\eqref{eq:finalfinal}. Here we have three kinds of terms: the
$\tilde{R}$'s, those involving $Z_\a$, $W_\a$ (which coincide with
$\mathfrak{A}_{5,0}$, with our notation in \eqref{eq:estA50}) and
the $\tilde{\mathfrak{A}}$'s.

For the $\tilde{R}$'s, since $v_{\underline{\Phi}}$ is odd in $z$,
the products with the even terms will vanish. The products of the
odd terms (notice that the two phases cancel and we use the change
of variables $s \mapsto \e s$) instead give us
$$
  \e^2 \Re \int_{\tilde{D}_\e}  (\tilde{R}_{r,o}+ \tilde{R}_{r,o,f_1})
  \overline{v_{\underline{\Phi}}} dV_{\tilde{g}_\e} + \e^2 \Re \int_{\tilde{D}_\e} i
(\tilde{R}_{i,o}^\Phi +\tilde{R}_{i,o,f_1}) \overline{v_{\underline{\Phi}}}
 dV_{\tilde{g}_\e} = - \e \frac{p-1}{2\th} C_0
\int_0^L \langle \mathfrak{J} (\Phi),  \underline{\Phi} \rangle \; d \ov{s} + \tilde{R}_0,
$$
where $C_0 = \int_{\R^{n-1}} U(y)^2 dy$ and $|\tilde{R}_0| \leq C \d
\e \|\underline{\Phi}\|_{L^2([0,L])}$. To explain why this estimate
holds, we notice first that $- \frac{p-1}{2\th} C_0 \langle
\mathfrak{J} (\Phi), \underline{\Phi} \rangle$ is exactly the first
term of $v_{\underline{\Phi}}$ multiplied by $\tilde{R}_{r,o} +
\tilde{R}_{r,o,f_1}$, as shown in Subsections 4.1 and 4.2 in
\cite{mmm1} (the factor $h^{\frac{p+1}{4}}$ in \eqref{eq:K1d} is
needed precisely to cancel the factor $\frac{1}{hk}$ in the last
formula of Subsection 4.2 in \cite{mmm1}). The remaining terms in
the last equation are given either by products of the imaginary part
of $v_{\underline{\Phi}}$ and the imaginary $\tilde{R}$'s or that of
$\tilde{R}_{r,o} + \tilde{R}_{r,o,f_1}$ and the last term in
$v_{\underline{\Phi}}$. In the latter case for example, we obtain a
quantity bounded by
$$
  C \e^2 \int_0^{L/\e} (|\Phi| + |\Phi'| + |\Phi''|) \e^2
 \underline{\Phi}'' d s \leq C \e^2 \d^2 \|\underline{\Phi}\|_{L^2([0,L])}.
$$
The last inequality follows from \eqref{eq:bdc1c2} and the fact that
$\underline{\Phi}$ satisfies the first condition in
\eqref{eq:bdPhif2d}. On the other hand, the terms involving
$\underline{\Phi}'$ once integrated will be bounded by $C \e^2 \d
\|\underline{\Phi}\|_{L^2([0,L])}$, still by \eqref{eq:bdPhif2d}.

Concerning $\mathfrak{A}_{5,0}$, we next claim that for any $m \in
\N$ one has
\begin{equation}\label{eq:claim22}
  \left| \Re \int_{\tilde{D}_\e} \overline{v_{\underline{\Phi}}} \;
  \mathfrak{A}_{5,0} dV_{\tilde{g}_\e} \right| \leq C \e^m
  \|\underline{\Phi}\|_{L^2([0,L])} \qquad \qquad \hbox { as } \e \to 0.
 \end{equation}
To see this, notice that $\Phi$ satisfies \eqref{eq:bdPhif2d} while
$\mathfrak{A}_{5,0}$ arises from functions involving $v_\d$ (in
particular $\b$, see \eqref{eq:exprbeta}): since $j$ ranges between
$- \frac{\d^2}{\e}$ and $\frac{\d^2}{\e}$, the main modes of $\b$
are much higher than the ones of $\Phi$. Hence, using Fourier
cancelation as in Lemma \ref{l:cancfour}, one can deduce
\eqref{eq:claim22}. It is also easy to see that
\begin{equation}\label{eq:claim33}
   \left| \Re \int_{\tilde{D}_\e} \overline{v_{\underline{\Phi}}} \;
  \sum_{j=1}^6 \tilde{\mathfrak{A}}_j dV_{\tilde{g}_\e} \right| \leq C(c_1,c_2,c_3) \e^2
 \|\underline{\Phi}\|_{L^2([0,L])}.
 \end{equation}
It remains finally to consider the product of $v_{\underline{\Phi}}$
and the last three terms in \eqref{eq:tildeSe}. Indeed, since these
are either superlinear in $\hat{\phi}$ (see \eqref{eq:Ne1}) or
contain $\var(\hat{\phi})$ (see \eqref{eq:qqqqq1}), they are of
lower order compared to \eqref{eq:small}.

Using \eqref{eq:small}-\eqref{eq:claim33} and the above arguments we
finally obtain that, if $\check{v}$ is as in Proposition
\ref{p:exhatphi}, then
\begin{equation}\label{eq:estR1}
  \int_{\tilde{D}_\e} \check{v} \; \overline{v_{\underline{\Phi}}} dV_{\tilde{g}_\e} =
  - \e \frac{p-1}{2\th} C_0 \int_0^L \langle \mathfrak{J} (\Phi), \underline{\Phi} \rangle
  \; d \ov{s} + R_1; \qquad \quad
  |R_1| \leq C(c_1, c_2, c_3) \e^2 \|\underline{\Phi}\|_{L^2([0,L])}.
\end{equation}
Similarly, using the estimates in Section \ref{s:as} one finds that
if $\tilde{\check{v}}$ corresponds to the triple $(\tilde{\Phi}$,
$\tilde{f}_2$, $\tilde{\b})$, then
\begin{equation}\label{eq:claim44}
    \int_{\tilde{D}_\e} (\check{v} - \tilde{\check{v}}) \; \overline{v_{\underline{\Phi}}}
   dV_{\tilde{g}_\e} = - \e \frac{p-1}{2\th} C_0 \int_0^L \langle \mathfrak{J}
 (\Phi - \tilde{\Phi}), \underline{\Phi} \rangle \; d \ov{s} + \tilde{R}_1,
\end{equation}
 where $\tilde{R}_1$ satisfies
\begin{equation}\label{eq:lipR1}
    |\tilde{R}_1| \leq  C(c_1, c_2, c_3) \left( \d \e \|\Phi - \tilde{\Phi}\|_{H^2([0,L])}
    + \e^2 \|f_2 - \tilde{f}_2\|_{H^2([0,L])} + \d \|\b - \tilde{\b}\|_\sharp \right)
  \|\underline{\Phi}\|_{L^2}.
\end{equation}

\subsubsection{Projection onto $\tilde{K}_{2,\d}$}\label{sss:pr2}

For this projection we will be more sketchy since most of the
arguments of the previous one can be applied. If $\underline{f}_2$
satisfies the second condition in \eqref{eq:bdPhif2d}, we consider
the function
$$
  v_{\underline{f}_2} = h(\e s)^{\frac 12} \left( i \underline{f}_2(\e s) U(k z)
  + 2 \e \frac{f' \underline{f}'_2(\e s)}{k} \tilde{U}(k z) - i \e^2
  \frac{\underline{f}''_2(\e s)}{k^2} \mathfrak{W}(k z) \right).
$$
As for the previous case, the main contribution to the projection is
given by the product of the first term in $v_{\underline{f}_2}$ and
the imaginary parts of $S_\e(\tilde{\Psi}_{2,\e})$ listed in
\eqref{eq:finalfinal} which are even in $z$.

We denote by $\hat{R}_{i,e,f_2}$ the sum of all imaginary even terms
of order $\e^3$ appearing in the equation, namely
$\mathfrak{A}_{1,i,e}$, $\mathfrak{A}_{3,i,e}$ and
$\mathfrak{A}_{4,i,e} = F_{4,i,e}(\ov{s})$, see \eqref{eq:mfA1ie},
\eqref{eq:estA3ie} and \eqref{eq:estA401}
\begin{eqnarray*}
  \hat{R}_{i,e,f_2} & = & 2h'f_2'U+2hf_2'k'\n
U\cdot z+2f'f_2'w_{i,e}+f_2''hU+4f'\pa_s(hf'f_2'\tilde{U}) \\
  & + & 2f''hf'f_2'\tilde{U}-2(p-1) h^{p-1} |U|^{p-2} f'f_2'
  \tilde{U} w_{i,e} + F_{4,i,e}(\ov{s}) := \tilde{R}_{i,e,f_2} + F_{4,i,e}(\ov{s}).
\end{eqnarray*}
Notice that $\tilde{R}_{i,e,f_2}$ coincides with the function
$\tilde{R}_{i,e,f_1}$ in \eqref{eq:esSe} (see Subsection 3.3 in
\cite{mmm1} for the precise expression) if we replace $f_1$ with
$f_2$. Therefore, from estimates similar to the previous ones (which
use especially the computations in Subsection 4.1 in \cite{mmm1}) we
find
\begin{equation}\label{eq:claim55}
    \int_{\tilde{D}_\e} \check{v} \; \overline{v_{\underline{f}_2}}
  dV_{\tilde{g}_\e} = \e^2 C_0 \int_0^L T(f_2) \underline{f}_2 \; d \ov{s} +
  \e^2 \int_0^L \left( \int_{\R^{n-1}} F_{4,i,e} U(k(\ov{s})) \right)
   \underline{f}_2 \; d \ov{s} + R_2,
\end{equation}
where $C_0 = \int_{\R^{n-1}} U(y)^2 dy$, where
\begin{equation}\label{eq:defT}
T(f_2) = \pa_{\ov s}\left( \frac{h^2f_2'}{(p-1)k^{n+1}}\left[(p-1)h^{p-1}- 2
\s \mathcal{A}^2h^{2\s} \right] \right),
\end{equation}
and where $R_2$ satisfies
\begin{equation}\label{eq:estR2}
    |R_2| \leq C(c_1, c_2, c_3) \d \e^2 \|\underline{f}_2\|_{L^2([0,L])}.
\end{equation}
Moreover,  if $\tilde{\check{v}}$ corresponds to the triple
$(\tilde{\Phi}$, $\tilde{f}_2$, $\tilde{\b})$, then
\begin{equation}\label{eq:lipr20}
    \int_{\tilde{D}_\e} (\check{v} - \tilde{\check{v}}) \;  \overline{v_{\underline{f}_2}}
  dV_{\tilde{g}_\e} =
  \e^2 C_0 \int_0^L T(f_2 - \tilde{f}_2) \underline{f}_2 \; d \ov{s} + \tilde{R}_2,
\end{equation}
 with
\begin{equation}\label{eq:lipR2}
    |\tilde{R}_2| \leq C(c_1, c_2, c_3) \left( \d \e^2 \|f_2 -
    \tilde{f}_2\|_{H^2([0,L])} + \d \e \|\Phi - \tilde{\Phi}\|_{H^2([0,L])}
    +  \d \|\b - \tilde{\b}\|_\sharp  \right) \|\underline{f}_2\|_{L^2([0,L])}.
\end{equation}

\subsubsection{Projection onto $\tilde{K}_{3,\d}$}\label{sss:pr3}

To compute the last components of the projection we recall our
notation in Subsection \ref{ss:appker}, and define
\begin{equation*}
    \underline{\b}(\e s) = \sum_{j=-\frac{\d^2}{\e}}^{\frac{\d^2}{\e}} \underline{b}_j
    \b_j(\e s); \qquad \qquad
    v_{\underline{\b}} = \sum_{j=-\frac{\d^2}{\e}}^{\frac{\d^2}{\e}}
\underline{b}_j \tilde{v}_{3,j}.
\end{equation*}
As for the previous cases, the main contribution to the projection
comes here still from $S_\e(\tilde{\Psi}_{2,\e})$. In particular,
following the arguments for $\tilde{K}_{1,\d}$, when testing on
$v_{\underline{\b}}$, by Fourier cancelation and parity the major
terms are indeed $\mathfrak{A}_{5,0}$, $\mathfrak{A}_{5,r,e}$ and
$\mathfrak{A}_{5,i,e}$. With straightforward computations one finds
that
\begin{equation}\label{eq:prrrrr}
    \int_{\tilde{D}_\e} \check{v} \; \overline{v_{\underline{\b}}} dV_{\tilde{g}_\e} =
  \frac{1}{\e} \int_0^L \L(\b,\xi,\underline{\b},\underline{\xi}) \; d \ov{s} + \hat{R}_3,
\end{equation}
where
\begin{eqnarray*}
& & \L(\b,\xi,\underline{\b},\underline{\xi}) = \b\underline{\b} Q_{4,\a} - \e^2\b''\underline{\b}
Q_{1,\a} -2\e \xi'
f'\underline{\b}Q_{3,\a} + \xi \underline{\xi} Q_{5,\a} -
\e^2\xi''\, \underline{\xi} Q_{2,\a}
- \e f''\left(\xi \underline{\b} - \underline{\xi} \b\right)
Q_{3,\a} \\ & + & 2\e\b'f' \underline{\xi} Q_{3,\a} - 2\e^2\a'\left(\b'\underline{\b} Q_{6,\a}
+ \underline{\xi} \xi'
Q_{7,\a} \right)  - 2\e^2k'\left(\b'\underline{\b} Q_{10,\a}
+ \underline{\xi} \xi' Q_{11,\a} \right)
\nonumber\\
&-&2\e f' k' \left( \xi \underline{\b} Q_{12,\a} -
\b \underline{\xi} Q_{13,\a} \right) \nonumber
 - 2\e f'\a'\left( \xi \underline{\b} Q_{8,\a}
- \underline{\xi} \b Q_{9,\a} \right)-(p-1)\e h^{p-2} \left( \underline{\xi} \b + \xi \underline{\b} \right)
Q_{14,\a}; \nonumber
\end{eqnarray*}
$$
  Q_{4,\a}(\ov{s}) =
  \int_{\R^{n-1}}Z_{\a(\ov{s})} \mathcal{L}_rZ_{\a(\ov{s})}; \qquad Q_{5,\a}(\ov{s}) =
\int_{\R^{n-1}}W_{\a(\ov{s})}\mathcal{L}_iW_{\a(\ov{s})}; \qquad Q_{6,\a}(\ov{s}) =
  \int_{\R^{n-1}}Z_{\a(\ov{s})} \frac{\pa Z_{\a(\ov{s})}}{\pa \a};
$$
$$
  Q_{7,\a}(\ov{s}) =
\int_{\R^{n-1}}W_{\a(\ov{s})} \frac{\pa W_{\a(\ov{s})}}{\pa \a} \qquad
  Q_{8,\a}(\ov{s}) = \int_{\R^{n-1}}Z_{\a(\ov{s})} \frac{\pa W_{\a(\ov{s})}}{\pa\a};
\qquad Q_{9,\a}(\ov{s}) =
\int_{\R^{n-1}}W_{\a(\ov{s})} \frac{\pa Z_{\a(\ov{s})}}{\pa\a};
$$
$$
  Q_{10,\a}(\ov{s}) = \int_{\R^{n-1}}Z_{\a(\ov{s})} \n_z Z_{\a(\ov{s})} \cdot z;
\qquad Q_{11,\a}(\ov{s}) =
  \int_{\R^{n-1}}W_{\a(\ov{s})} \n_z W_{\a(\ov{s})} \cdot z;
$$
$$
  Q_{12,\a}(\ov{s}) = \int_{\R^{n-1}}Z_{\a(\ov{s})} \n_zW_{\a(\ov{s})} \cdot z;
\qquad Q_{13,\a}(\ov{s}) =
  \int_{\R^{n-1}}W_{\a(\ov{s})} \n_z Z_{\a(\ov{s})} \cdot z;
$$
$$
  Q_{14,\a}(\ov{s}) = \int_{\R^{n-1}}U(kz)^{p-2}w_{i,e}W_{\a(\ov{s})} Z_{\a(\ov{s})}.
$$
\begin{equation}\label{eq:estR3}
  |\hat{R}_3| \leq C(c_1, c_2, c_3) \d \e^2 \|\underline{\b}\|_{L^2([0,L])}.
\end{equation}
After some manipulation using the fact that $(Z_\a, W_\a)$ solve
\eqref{eq:systuva} with $\eta_\a = 0$, the normalization
$\int_{\R^{n-1}} (Z_\a^2 + W_\a^2) = 1$ and  some integration by
parts in $z$ we find that
\begin{equation}\label{eq:P3}
    \frac{1}{\e} \int_0^L \L(\b, \xi, \underline{\b}, \underline{\xi}) \; d \ov{s} =
   \frac 1 \e \int_0^L \L_0(\b, \xi, \underline{\b}, \underline{\xi}) \; d \ov{s}  +
    \int_0^L \L_1(\b, \xi, \underline{\b}, \underline{\xi}) \; d \ov{s},
\end{equation}
where
\begin{eqnarray}\label{eq:L0} \nonumber
  \L_0(\b,\xi,\underline{\b},\underline{\xi}) & = & Q_{1,\a} \left( \e^2  \b' \underline{\b}' -
  k^2 \a^2 \b \underline{\b} \right)  + Q_{2,\a}
  \left( \e^2   \xi' \underline{\xi}' - \a^2 k^2 \xi
  \underline{\xi}  \right) \\ & + & 2 f' Q_{3,\a} \left( \e \b'
  \underline{\xi} - \e \xi' \underline{\b} - k \a \b
  \underline{\b} - k \a \xi \underline{\xi} \right);
\end{eqnarray}
$$
  \L_1 = (\b \underline{\xi} + \xi \underline{\b}) \mathfrak{g}(\ov{s}); \qquad
\quad  \mathfrak{g}(\ov{s}) = \left[ f'' Q_{3,\a} + 2 f' k' Q_{13,\a}
+ 2 f' \a' Q_{9,\a} - (p-1) h^{p-2}
  Q_{14,\a} \right].
$$

Now we notice that, by \eqref{eq:bbj}, one has
$$
  \b \underline{\xi} + \xi \underline{\b} = - \frac{\e}{k \a}
  \left[ \sum_{j,l=-\frac{\d^2}{\e}}^{\frac{\d^2}{\e}} b_j \underline{b}_l (\xi'_j \xi_l
  + \xi_j \xi'_l) - F_1 \sum_{j,l=-\frac{\d^2}{\e}}^{\frac{\d^2}{\e}} b_j \underline{b}_l
  (\nu_j \xi_l \xi'_j + \nu_l \xi_j \xi'_l) \right],
$$
where $F_1 = \frac{Q_{1,\a}}{k^2\a^2+2f'k\a Q_{3,\a}}$. Integrating
by parts in $\ov{s}$ and using \eqref{eq:asynul} we find that
\begin{equation}\label{eq:F1}
  \int_0^L \L_1(\b, \xi, \underline{\b}, \underline{\xi}) \; d \ov{s}
  = \e \int_0^L \xi \underline{\xi} \left( \frac{\mathfrak{g}}{k\a} \right)'(\ov{s}) d \ov{s}
  + O(\d^2) \|\b\|_{L^2([0,L])} \|\underline{\b}\|_{L^2([0,L])}.
\end{equation}
Finally, combining \eqref{eq:prrrrr}, \eqref{eq:estR3},
\eqref{eq:P3} and \eqref{eq:F1} we deduce
\begin{equation}\label{eq:PP3}
    \int_{\tilde{D}_\e} \check{v} \; \overline{v_{\underline{\b}}}
     dV_{\tilde{g}_\e} = \frac{1}{\e} \int_0^L \L_0(\b, \xi, \underline{\b}, \underline{\xi})
    \; d \ov{s} + R_3,
\end{equation}
where
\begin{equation}\label{eq:estR32}
  |R_3| \leq C(c_1, c_2, c_3) (\d \e^2 + (\e + \d^2) \|\b\|_{L^2([0,L])})
  \|\underline{\b}\|_{L^2([0,L])}.
\end{equation}
Analogously we obtain
\begin{equation}\label{eq:Lippr3}
    \int_{\tilde{D}_\e} (\check{v} - \tilde{\check{v}}) \;
\overline{v_{\underline{\b}}}  dV_{\tilde{g}_\e} = \frac{1}{\e}
  \int_0^L \L_0(\b - \tilde{\b}, \xi - \tilde{\xi}, \underline{\b},
\underline{\xi}) + \tilde{R}_3,
\end{equation}
where $\tilde{R}_3$ satisfies
\begin{equation}\label{eq:lipR3}
  |\tilde{R}_3| \leq C(c_1, c_2, c_3) \d \left( \e \|\Phi -\tilde{\Phi}\|_{H^2([0,L])}
  + \e^2 \|f_2 -\tilde{f}_2\|_{H^2([0,L])} + \|\b -\tilde{\b}\|_\sharp \right)
\|\underline{\b}\|_{L^2([0,L])}.
\end{equation}

\begin{rem}\label{r:diagLo} Let us consider the eigenvalue problem
in $(\b,\xi)$
$$
  \int_0^L \L_0(\b, \xi, \underline{\b}, \underline{\xi}) = \nu \int_0^L (Q_{1,\a}
  \b \underline{\b} + Q_{2,\a} \xi \underline{\xi}) \qquad \hbox{ for all } (\underline{\b},
\underline{\xi})
$$
where $Q_1, Q_2$ are defined in \eqref{eq:Q1Q2Q3}. Then the
eigenvalue equation is the following
\begin{equation}\label{eq:systbjxj2}
  \left\{
    \begin{array}{ll}
      - \e^2 \frac{(Q_{1,\a} \b')'}{Q_{1,\a}} - k^2 \a^2 \b - 2 f'
    \frac{Q_{3,\a}}{Q_{1,\a}} (\e \xi' + k \a \b) = \nu \b; &  \\ - \e^2
\frac{(Q_{2,\a} \xi')'}{Q_{2,\a}} - k^2 \a^2 \xi + 2 f'
    \frac{Q_{3,\a}}{Q_{2,\a}} (\e \b' - k \a \xi) = \nu \xi. &
    \end{array}
  \right.
\end{equation}
By \eqref{eq:systbjxj}, the couple of functions $(\b_j,\xi_j)$
constructed in Subsection \ref{ss:appker} represents a family of
approximate eigenfunctions corresponding to $\nu = \nu_j$.
\end{rem}

\subsection{The contraction argument}\label{ss:contr}

The usual procedure in performing a fixed point argument is to apply
to the equation an invertible linear operator first. From the
expansions in the last subsection, we showed that the main terms in
the projections onto $\tilde{K}_\d$ are the operators
$\mathfrak{J}$, $T$ and $\L_0$ (the latter is identified by duality
with the associated quadratic form), see \eqref{eq:estR1},
\eqref{eq:claim55} and \eqref{eq:PP3}. By our non-degeneracy
assumption on $\g$, $\mathfrak{J}$ is invertible and the same holds
also for $T$, since it is coercive (and in divergence form). It
remains then to invert $\L_0$, which is the content of the next
result: before stating it we introduce some notation. Using the
symbology of Subsection \ref{ss:appker} we define the spaces
$$
  X_{1,\d} = span \left\{ \varphi_j \; : \; j = 0, \dots, \frac{\d}{\e} \right\};
 \qquad  X_{2,\d} = span \left\{ \o_j \; : \; j = 0, \dots, \frac{\d}{\e} \right\};
\qquad
$$
$$
  X_{3,\d} = span \left\{ \b_j \; : \; j = -\frac{\d^2}{\e},
   \dots, \frac{\d^2}{\e} \right\},
$$
with $X_{1,\d}$, $X_{2,\d}$ endowed with the $H^2$ norm on $[0,L]$,
and $X_{3,\d}$ with the $\| \cdot \|_\sharp$ norm.

We also call $Y_{1,\d}, Y_{2,\d}, Y_{3,\d}$ the same spaces of
functions, but endowed with weighted $L^2$ norms: by the
normalization after \eqref{eq:eigen} it is natural to put the
weights $h^\th$ and $h^{-\s}$ on $Y_{1,\d}$ and $Y_{2,\d}$
respectively. Concerning $Y_{3,\d}$, by Remark \ref{r:diagLo}, we
will endow it with the  product $(\b,\underline{\b})_{Y_{3,\d}} =
\int_0^L (Q_{1,\a}
  \b \underline{\b} + Q_{2,\a} \xi \underline{\xi}) d \ov{s}$ where,
as above, $\xi$ is related to $\b$ by \eqref{eq:bbj} and
\eqref{eq:exprbeta}. Notice that by \eqref{eq:eigen} $\mathfrak{J}$
and $T$ are exactly diagonal from $X_{1,\d}$ to $Y_{1,\d}$ and from
$X_{2,\d}$ to $Y_{2,\d}$ respectively, while $\L_0$  is nearly
diagonal (see also \eqref{eq:systbjxj}).

\begin{lem}\label{l:invF0} Letting $\Pi_{Y_{3,\d}}$ denote the
orthogonal projection onto $Y_{3,\d}$, there exists a sequence $\e_k
\to 0$ such that $\L_0$ is invertible from $X_{3,\d}$ into
$Y_{3,\d}$ and  such that its inverse satisfies $\|(\Pi_{Y_{3,\d}}
\L_0)^{-1}\| \leq \frac{C}{\e_k}$ for some fixed constant $C$.
\end{lem}

\begin{pf} First of all we show that there exists $\e_k \to 0$ such
that $\Pi_{Y_{3,\d}} \L_0$ cannot have eigenvalues in $Y_{3,\d}$
smaller in absolute value than $C^{-1} \e_k$: after this, we
estimate the (stronger) $X_{3,\d}$ norm of its inverse.

To prove the claim we apply Kato's theorem (see \cite{ka}, page
445): the latter allows to compute the derivative of an eigenvalue
$\nu(\e)$ of $\Pi_{Y_{3,\d}} \L_0$ with respect to $\e$. The
(possibly multiple) value of this derivative is given by the
eigenvalues of $\Pi_{Y_{3,\d}} \pa_\e \L_0$, restricted to the
$\nu(\e)$-eigenspace of $\Pi_{Y_{3,\d}} \L_0$.

Suppose that $\b$ satisfies the eigenvalue equation $\Pi_{Y_{3,\d}}
\L_0 \b = \nu \b$, which is equivalent to
\begin{equation}\label{eq:hahaha}
    \int_0^L \L_0(\b, \xi, \underline{\b}, \underline{\xi}) = \nu \int_0^L (Q_{1,\a}
  \b \underline{\b} + Q_{2,\a} \xi \underline{\xi}) \qquad \hbox{ for all }
(\underline{\b}, \underline{\xi}) \hbox{ with } \underline{\b} \in Y_{3,\d}.
\end{equation}
Looking at the powers of $\e$ in $\L_0$, see \eqref{eq:L0}, we write
$\L_0 = \L_{0,0} + \e \L_{0,1} + \e^2 \L_{0,2}$: notice that
$\L_{0,0}$ is negative-definite and $\L_{0,2}$  positive-definite.
We also point out that, since $f'$ satisfies \eqref{eq:f'Cintr}, for
$f(\e s)/\e$ to be $L/\e$-periodic, when we vary $\e$ also
$\mathcal{A}$ needs to be adjusted. Precisely, since the total
variation of phase in \eqref{eq:prof} is
$$
  \mathcal{A} \int_0^{L/\e} h(\e s)^\s ds = \frac{\mathcal{A}}{\e} \int_0^L
  h(\ov{s}) d \ov{s} = const.,
$$
when differentiating with respect to $\e$ we find that $\frac{\pa
\mathcal{A}}{\pa \e} = \frac{\mathcal{A}}{\e}$.  Hence, applying
Kato's theorem we find
\begin{equation}\label{eq:mM}
    \frac{\pa \nu}{\pa \e} \in \left[  \min_{\b_1,\b_2 \neq 0}
  \Theta(\b_1,\b_2), \max_{\b_1,\b_2 \neq 0} \Theta(\b_1,\b_2) \right],
\end{equation}
where
$$
  \Theta(\b_1,\b_2) = \frac{\int_0^L
  (\L_{0,1} + 2 \e \L_{0,2})(\b_1, \xi_1, \b_2, \xi_2)}{\int_0^L (Q_{1,\a}
  \b_1 \b_2 + Q_{2,\a} \xi_1 \xi_2)} + \frac 1 \e \frac{\int_0^L (2 f' Q_{3,\a} \left( \e \b'_1
  \xi_2 - \e \xi'_1 \b_2 - k \a \b_1 \b_2 - k \a \xi_1 \xi_2 \right)}{\int_0^L (Q_{1,\a}
  \b_1 \b_2 + Q_{2,\a} \xi_1 \xi_2)},
$$
and where $(\b_1,\xi_1), (\b_2, \xi_2)$ are functions satisfying
\eqref{eq:hahaha}: using this, $\Theta(\b_1,\b_2)$ can be written as
\begin{eqnarray*}
  & &  \frac{\frac 1 \e \int_0^L (\L_0 - \L_{0,0})(\b_1, \xi_1, \b_2,
  \xi_2) + \e \int_0^L \L_{0,2}(\b_1, \xi_1, \b_2,
  \xi_2)}{\int_0^L (Q_{1,\a} \b_1 \b_2 + Q_{2,\a} \xi_1 \xi_2)}
\\ & + & \frac 1 \e \frac{\int_0^L (2 f' Q_{3,\a} \left( \e \b'_1
  \xi_2 - \e \xi'_1 \b_2 - k \a \b_1 \b_2 - k \a \xi_1 \xi_2 \right)}{\int_0^L (Q_{1,\a}
  \b_1 \b_2 + Q_{2,\a} \xi_1 \xi_2)}
\\
   & = & \frac{\nu}{\e} + \frac{\e \int_0^L (Q_{1,\a} \b'_1 \b'_2 + Q_{2,\a} \xi'_1 \xi'_2) +
  \frac 1 \e \int_0^L [k^2\a^2 (Q_{1,\a} \b_1 \b_2 + Q_{2,\a} \xi_1 \xi_2) + 2 f' k \a
  Q_{3,\a} (\b_1 \b_2 + \xi_1 \xi_2)] }{\int_0^L (Q_{1,\a} \b_1 \b_2 + Q_{2,\a} \xi_1 \xi_2)}
\\ & + & \frac 1 \e \frac{\int_0^L (2 f' Q_{3,\a} \left( \e \b'_1
  \xi_2 - \e \xi'_1 \b_2 - k \a \b_1 \b_2 - k \a \xi_1 \xi_2 \right)}{\int_0^L (Q_{1,\a}
  \b_1 \b_2 + Q_{2,\a} \xi_1 \xi_2)}.
\end{eqnarray*}
Applying \eqref{eq:asynul}, \eqref{eq:bbj} and $Q_{1,\a} + Q_{2,\a}
= 1$ (see \eqref{eq:Q1Q2Q3} and the lines after
\eqref{eq:systuvadiff}), the last expression simplifies as
$$
  \frac{\nu}{\e} + \frac{1}{\e} \frac{\int_0^L (2 \a^2 k^2 + 4 f' \a
  k Q_{3,\a}) \xi_1 \xi_2}{\int_0^L \xi_1 \xi_2} + O(\d^2) \frac{1}{\e}.
$$
Since the numerator is symmetric in $\xi_1, \xi_2$, the infimum of
the above ratio is realized by some $\xi_0$, so by \eqref{eq:mM} and
the latter formula we find
\begin{equation}\label{eq:bddve}
    \frac{\pa \nu}{\pa \e} \geq \frac{\nu}{\e} + \frac{1}{\e} \frac{\int_0^L
(2 \a^2 k^2 + 4 f' \a  k Q_{3,\a}) \xi_0^2}{\int_0^L \xi_0^2} + O(\d^2) \frac{1}{\e}
\geq \frac{1}{\e} \left[ \nu + \inf_{[0,L]} \left( 2 \a^2 k^2 + 4 f'
\a  k Q_{3,\a} \right) - C \d^2  \right].
\end{equation}
Notice that for $\nu$ and $\d$ sufficiently small, the coefficient
of $\frac{1}{\e}$ in the above formula is positive and uniformly
bounded away from zero. From \eqref{eq:systbjxj} and the asymptotics
in \eqref{eq:asynul} (which follows from the Weyl's formula), one
can show that $\Pi_{Y_{3,\d}} \L_0$ has a number of negative
eigenvalues of order $\frac{\d^2}{\e}$. This fact and
\eqref{eq:bddve} yield the desired claim, which can be obtained as
in \cite{malm2}, Proposition 4.5: since the argument is quite
similar, we omit the details.

\

\noindent The above claim provides invertibility of $\Pi_{Y_{3,\d}}
\L_0$ in $Y_{3,\d}$, and gives
\begin{equation}\label{eq:bdY3Y3}
    \|(\Pi_{Y_{3,\d}} \L_0)^{-1} \b\|_{Y_{3,\d}} \leq \frac{C}{\e}
\|\b\|_{Y_{3,\d}} \qquad \hbox{ for any } \b \in Y_{3,\d}:
\end{equation}
we want next to estimate the $X_{3,\d}$ norm of $(\Pi_{Y_{3,\d}}
\L_0)^{-1} \b$.  Let $\b =
\sum_{j=-\frac{\d^2}{\e}}^{\frac{\d^2}{\e}} b_j \b_j$ and suppose
$\hat{\b} = \sum_{j=-\frac{\d^2}{\e}}^{\frac{\d^2}{\e}} \hat{b}_j
\b_j$ is such that $\Pi_{Y_{3,\d}} \L_0 \hat{\b} = \b$, in the sense
that
$$
  \int_0^L \L_0(\hat{\b}, \hat{\xi}, \underline{\b}, \underline{\xi}) = \int_0^L (Q_{1,\a}
  \b \underline{\b} + Q_{2,\a} \xi \underline{\xi}) \qquad \hbox{ for all }
(\underline{\b}, \underline{\xi}) \hbox{ with } \underline{\b} \in Y_{3,\d}.
$$
If $\underline{\b} = \sum_{j=-\frac{\d^2}{\e}}^{\frac{\d^2}{\e}}
\underline{b}_j \b_j$, then by \eqref{eq:systbjxj} integrating one
finds
$$
  \sum_{j=-\frac{\d^2}{\e}}^{\frac{\d^2}{\e}} \nu_j \hat{b}_j \underline{b}_j
  + O \bigg( \bigg( \sum_{j=-\frac{\d^2}{\e}}^{\frac{\d^2}{\e}} (\nu_j^2 + \e)^2
  \hat{b}_j^2 \bigg)^{\frac 12} \bigg( \sum_{l=-\frac{\d^2}{\e}}^{\frac{\d^2}{\e}}
 \underline{b}_l^2 \bigg)^{\frac 12} \bigg) \leq C \bigg( \sum_{j=-\frac{\d^2}{\e}}^{\frac{\d^2}{\e}}
 \hat{b}_j^2 \bigg)^{\frac 12} \bigg( \sum_{l=-\frac{\d^2}{\e}}^{\frac{\d^2}{\e}}
 \underline{b}_l^2 \bigg)^{\frac 12}.
$$
Choosing $\underline{b}_j = \hat{b}_j$ for $j > 0$ and
$\underline{b}_j = - \hat{b}_j$ for $j < 0$, from the asymptotics of
$\nu_j$ in \eqref{eq:asynul} we obtain for $\tilde{C}_1 > 0$
sufficiently large that
$$
  \e \sum_{\tilde{C}_1 \leq j \leq \frac{\d^2}{\e}} |j| \; \hat{b}_j^2
  \leq C \sum_{j=-\frac{\d^2}{\e}}^{\frac{\d^2}{\e}}
 \hat{b}_j^2 \leq C \|\hat{\b}\|_{Y_{3,\d}}^2.
$$
By \eqref{eq:bdY3Y3} we have  $\|\hat{\b}\|_{Y_{3,\d}} \leq
\frac{C}{\e} \|\b\|_{Y_{3,\d}}$, so recalling \eqref{eq:normsharp}
we get $\|\hat{\b}\|_{X_{3,\d}}^2 := \|\hat{\b}\|_{\sharp}^2 \leq C
\|\hat{\b}\|_{Y_{3,\d}}^2 \leq \frac{C^2}{\e^2}
\|\b\|_{Y_{3,\d}}^2$, which yields the conclusion.
\end{pf}

\

\begin{pfn} {\sc of Theorem \ref{t:main}} Let us introduce the operators
$$
  G_l : X_{1,\d} \times X_{2,\d} \times X_{3,\d} \to Y_{l,\d}, \qquad l = 1, 2, 3,
$$
defined by duality as
$$
  (G_1(\Phi, f_2, \b), \underline{\Phi})_{Y_{1,\d}} = \int_{\tilde{D}_\e} \check{v} \;
  \overline{v_{\underline{\Phi}}} dV_{\tilde{g}_\e}; \qquad  (G_2(\Phi,  f_2, \b),
  \underline{f}_2)_{Y_{2,\d}} = \int_{\tilde{D}_\e} \check{v} \; \overline{v_{\underline{f}_2}}
   dV_{\tilde{g}_\e};
$$
$$
   (G_3(\Phi,  f_2, \b), \underline{\b})_{Y_{3,\d}} =
  \int_{\tilde{D}_\e} \check{v} \; \overline{v_{\underline{\b}}} dV_{\tilde{g}_\e},
$$
where $\check{v} = \check{v}(\Phi, f_2, \b)$ is the function
appearing in Proposition \ref{p:exhatphi}.

By Proposition \ref{p:sec2final}, equation \eqref{eq:new} (or
\eqref{eq:pe}) is solved if and only if $\check{v} = 0$: in the
above notations, this is equivalent to finding $(\Phi, f_2, \b)$
such that $G_l(\Phi, f_2, \b) = 0$ for every $l = 1, 2, 3$. If
$\e_k$ is the sequence given in Lemma \ref{l:invF0} then $\L_0$ is
invertible, and the condition $\check{v} = 0$ is equivalent to the
system (we set $\e = \e_k$)
\begin{equation}\label{eq:fethi2}
  \left\{
    \begin{array}{ll}
      \Phi = \mathfrak{G_1}(\Phi, f_2, \b) := - \frac 1 \e \tilde{\mathfrak{J}}^{-1}
      [G_1(\Phi, f_2, \b) - \e \tilde{\mathfrak{J}}(\Phi)] &  \\      \\
      f_2 - \check{f}_2 = \mathfrak{G_2}(\Phi, f_2, \b) := - \frac{1}{\e^2} \tilde{T}^{-1}
      \left[G_2(\Phi, f_2, \b) - \e^2 \tilde{T} f_2
      - \e^2 \int_{\R^{n-1}} F_{4,i,e} U(k(\ov{s})z)dz \right]   &  \\       \\
      \b = \mathfrak{G_3}(\Phi, f_2, \b) := - \e (\Pi_{Y_{3,\d}} \L_0)^{-1}
     \left[ G_3(\Phi, f_2, \b) - \frac 1 \e \Pi_{Y_{3,\d}} \L_0 \b \right],  &
    \end{array}
  \right.
  \end{equation}
where $\tilde{\mathfrak{J}} = - \frac{p-1}{2 \th} C_0 \mathfrak{J}$,
$\tilde{T} = C_0 T$ ($C_0 = \int_{\R^{n-1}} U(y)^2 dy$) and where
$\check{f}_2 = - \tilde{T}^{-1} (\int_{\R^{n-1}} F_{4,i,e}
U(k(\ov{s})z)dz)$. By \eqref{eq:estR1}-\eqref{eq:lipR2},
\eqref{eq:PP3}-\eqref{eq:Lippr3} and \eqref{eq:lipR3} one finds
$$
    \|\mathfrak{G_1}(0,0,0)\|_{X_{1,\d}} \leq C \e; \qquad \|\mathfrak{G_2}(0,0,0)\|_{X_{2,\d}}
    \leq C \d; \qquad \|\mathfrak{G_3}(0,0,0)\|_{X_{3,\d}} \leq C \d \e^2;
$$
moreover if $\Phi, f_2, \b$ satisfy the bounds \eqref{eq:bdc1c2},
\eqref{eq:bdbe3}, then
$$
  \|\mathfrak{G_1}(\Phi, f_2, \b) - \mathfrak{G_1}(\tilde{\Phi}, \tilde{f}_2, \tilde{\b})
  \|_{X_{1,\d}} \leq C(c_1,c_2,c_3) \left( \d \|\Phi - \tilde{\Phi}\|_{X_{1,\d}} +
   \e \|f_2 - \tilde{f}_2\|_{X_{2,\d}} + \frac \d \e \|\b - \tilde{\b}\|_{X_{3,\d}} \right);
$$
$$
  \|\mathfrak{G_2}(\Phi, f_2, \b) - \mathfrak{G_2}(\tilde{\Phi}, \tilde{f}_2, \tilde{\b})
  \|_{X_{2,\d}} \leq C(c_1,c_2,c_3) \left( \frac \d \e \|\Phi - \tilde{\Phi}\|_{X_{1,\d}} +
   \d \|f_2 - \tilde{f}_2\|_{X_{2,\d}} + \frac{\d}{\e^2} \|\b - \tilde{\b}\|_{X_{3,\d}} \right);
$$
$$
  \|\mathfrak{G_3}(\Phi, f_2, \b) - \mathfrak{G_3}(\tilde{\Phi}, \tilde{f}_2, \tilde{\b})
  \|_{X_{3,\d}} \leq C(c_1,c_2,c_3) \left( \d \e \|\Phi - \tilde{\Phi}\|_{X_{1,\d}} +
   \d \e^2 \|f_2 - \tilde{f}_2\|_{X_{2,\d}} + \d \|\b - \tilde{\b}\|_{X_{3,\d}} \right).
$$
We now consider the scaled norms $\e \| \cdot \|_{\hat{X}_{1,\d}} =
\| \cdot \|_{X_{1,\d}}$, $\d^{\frac 12} \| \cdot \|_{\hat{X}_{2,\d}}
= \| \cdot \|_{X_{2,\d}}$, $\e^2 \| \cdot \|_{\hat{X}_{3,\d}} = \|
\cdot \|_{X_{3,\d}}$: with this new notation the last formulas
become
\begin{equation}\label{eq:CCCC}
    \|\mathfrak{G_1}(0,0,0)\|_{\hat{X}_{1,\d}} \leq C; \qquad \|\mathfrak{G_2}(0,0,0)\|_{\hat{X}_{2,\d}}
    \leq C \d^{\frac 12 }; \qquad \|\mathfrak{G_3}(0,0,0)\|_{\hat{X}_{3,\d}} \leq C \d;
\end{equation}
$$
  \|\mathfrak{G_1}(\Phi, f_2, \b) - \mathfrak{G_1}(\tilde{\Phi}, \tilde{f}_2, \tilde{\b})
  \|_{\hat{X}_{1,\d}} \leq C(c_1,c_2,c_3) \left( \d \|\Phi - \tilde{\Phi}\|_{\hat{X}_{1,\d}} +
   \d^{\frac 12} \|f_2 - \tilde{f}_2\|_{\hat{X}_{2,\d}} + \d \|\b - \tilde{\b}\|_{\hat{X}_{3,\d}} \right);
$$
$$
  \|\mathfrak{G_2}(\Phi, f_2, \b) - \mathfrak{G_2}(\tilde{\Phi}, \tilde{f}_2, \tilde{\b})
  \|_{\hat{X}_{2,\d}} \leq C(c_1,c_2,c_3) \left( \d^{\frac 12} \|\Phi - \tilde{\Phi}\|_{\hat{X}_{1,\d}} +
   \d \|f_2 - \tilde{f}_2\|_{\hat{X}_{2,\d}} + \d^{\frac 12} \|\b - \tilde{\b}\|_{\hat{X}_{3,\d}} \right);
$$
$$
  \|\mathfrak{G_3}(\Phi, f_2, \b) - \mathfrak{G_3}(\tilde{\Phi}, \tilde{f}_2, \tilde{\b})
  \|_{\hat{X}_{3,\d}} \leq C(c_1,c_2,c_3) \left( \d \|\Phi - \tilde{\Phi}\|_{\hat{X}_{1,\d}} +
   \d^{\frac 32} \|f_2 - \tilde{f}_2\|_{\hat{X}_{2,\d}} + \d \|\b - \tilde{\b}\|_{\hat{X}_{3,\d}} \right).
$$
If $C$ is the constant appearing in \eqref{eq:CCCC}, from the last
four formulas we deduce that if $\d$ is sufficiently small then
$(\mathfrak{G_1},\mathfrak{G_2},\mathfrak{G_3})$ has a fixed point
in $\{ \| \cdot \|_{\hat{X}_{1,\d}} \leq 2 C \} \cap \{ \| \cdot -
\check{f}_2 \|_{\hat{X}_{2,\d}} \leq 2 C \d^{\frac 12} \} \cap \{ \|
\cdot \|_{\hat{X}_{3,\d}} \leq 2 C \d \}$. This, by the comments
before \eqref{eq:fethi2}, leads to a solution of \eqref{eq:pe} with
the desired asymptotics.
\end{pfn}

\

\begin{center}
{\bf Acknowledgments}
\end{center}

\noindent F. M. and A. M. are  supported by M.U.R.S.T within the
PRIN 2006 {\em Variational methods and nonlinear differential
equations}. F. M. is grateful to SISSA for the kind hospitality.

\

\end{document}